\documentclass{article}
\usepackage{amsmath,amssymb, amsthm, amsbsy,bm}
\usepackage[dvips]{graphicx}

\newtheorem{theorem}{Theorem}
\newtheorem{acknowledgement}[theorem]{Acknowledgement}

\newtheorem{corollary}[theorem]{Corollary}

\newtheorem{lemma}[theorem]{Lemma}

\newtheorem{remark}[theorem]{Remark}

\begin{document}

\title{Evolution of viscous vortex filaments and desingularization of the
Biot-Savart integral}
\author{Marco A. Fontelos \\
Instituto de Ciencias Matem\'{a}ticas (ICMAT, CSIC-UAM-UC3M-UCM), \\
C/ Nicol\'{a}s Cabrera 15, 28049 Madrid, Spain. \and Luis Vega \\
BCAM-UPV/EHU Bilbao, Spain.}
\maketitle

\begin{abstract}
We consider a viscous fluid with kinematic viscosity $\nu $ and initial data
consisting of a smooth closed vortex filament with circulation $\Gamma $. We
show that, for short enough time, the solution consists of a deformed
Lamb-Oseen vortex whose center (a filament) follows the binormal flow
dynamics plus leading order corrections that depend locally on the filament
curvature and the nonlocal interactions with distant parts of the filament.
In order to achieve this scale separation we require $\Gamma /\nu $ to be
sufficiently small.
\end{abstract}

\section{\protect\bigskip Introduction}

One of the most fundamental elements in fluid mechanics are the fluid
filaments. These are structures around which vorticity tends to concentrate.
Their formation, structure and evolution has been a central theme since the
early days of fluid dynamics and the main results bear the names of
Helmholtz, Lord Kelvin, Kirchhoff, etc (see the classical book by Saffman
\cite{S} for a general overview or the book by Darrigol \cite{D} for a
historical perspective). When a vortex filament is curved and its motion
affected by the velocity field created by itself, one could in principle
think of a closed-form equation to describe the evolution\ of the filament.
When the flow is inviscid, the velocity field is given by the Biot-Savart
law applied to the vorticity. This yields a singular integral and an
infinite velocity field at the filament. By replacing the zero-thickness
filament by a filament of radius $\varepsilon $, the singularity is removed
and the velocity is given by%
\begin{equation}
\mathbf{v}=-\frac{\Gamma }{4\pi }\kappa (s,t)(\log \varepsilon )\mathbf{b}%
(s,t)+O(1),  \label{velo}
\end{equation}%
where $\Gamma $ is the velocity circulation around the vortex filament (or
vortex filament strength), $\kappa (s,t)$ is the local curvature with $s$
being the arc-length parameter and $\mathbf{b}(s,t)$ is the binormal vector.
By neglecting the $O(1)$ terms in (\ref{velo}) one introduces the celebrated
localized induction approximation (LIA) and to the so-called binormal flow
dynamics for the vortex filament. This evolution problem was introduced by
Da Rios in 1906 \cite{DR} and has produced hundreds of studies on existence,
uniqueness and properties of the solutions in relation with the physical
situations that the equation intends to model. Obvious problems with the
binormal flow are the following: LIA is an ad hoc assumption that can only
be justified as leading order for a thin vortex tube in inviscid flow and it
is still unclear what is the role of vortex filaments and the binormal flow
dynamics in the case of viscous flows. Nevertheless, there is nowadays
overwhelming evidence that binormal flow is connected to observed physical
phenomena such as the motion of vortex rings or the oscillation of perturbed
vortex rings \cite{LN}, \cite{SL}. In the case of Euler flows, the existence
of vortex rings with a core of finite thickness and moving with constant
velocity was proved by Fraenkel \cite{F} and important progress has been
made recently to prove the existence of travelling helices \cite{W1} and
leapfrogging rings \cite{W2}, for instance. The more general question on
whether the binormal flow dynamics represents the evolution of a vortex
filament of thickness $\varepsilon $ for a general class of initial data has
been answered positively in \cite{JS2} under some hypothesis.

A vortex filament may be considered as a singular initial data for the
vorticity in the form
\begin{equation}
\boldsymbol{\omega }(\mathbf{x},0)=\frac{\Gamma }{2\pi }\delta _{G}\mathbf{t}%
(s,0),  \label{vor0}
\end{equation}%
where $\delta _{G}$ is the Dirac delta measure supported in the curve $G$, $%
\Gamma $ the vortex circulation and $\mathbf{t}(s,0)$ the tangent vector. A
natural question is then what is the later evolution of the filament
according to Navier-Stokes equations and whether there is a connection\ (or
not) with the binormal flow. In terms of the vorticity, Navier-Stokes system
for a fluid with kinematic viscosity $\nu $ is%
\begin{equation}
\frac{\partial \boldsymbol{\omega }}{\partial t}+\mathbf{v}\cdot \nabla
\boldsymbol{\omega }-\boldsymbol{\omega }\cdot \nabla \mathbf{v}-\nu \Delta
\boldsymbol{\omega }=0,  \label{vor1}
\end{equation}%
where $\mathbf{v}$ is the velocity field given by Biot-Savart integral of
the vorticity:%
\begin{equation}
\mathbf{v}(\mathbf{x},t)=K\ast \boldsymbol{\omega }(\mathbf{x},t)\equiv
\frac{1}{4\pi }\int \frac{\boldsymbol{\omega }(\mathbf{x}^{\prime },t)\times
(\mathbf{x}-\mathbf{x}^{\prime })}{\left\vert \mathbf{x}-\mathbf{x}^{\prime
}\right\vert ^{3}}d\mathbf{x}^{\prime }.  \label{vor2}
\end{equation}

Giga and Miyakawa, in a seminal paper \cite{GM}, prove the existence of weak
solutions with singular initial data. This includes initial data conditions
of Dirac delta functions supported in a smooth curve, i.e. a filament. The
assumptions on initial data are smallness in a Morrey space norm which, in
the case of a vortex filament, amounts to requiring that the filament
strength $\Gamma $ is sufficiently small. Although solutions are smooth at
later times, one can expect them to be concentrated around a curve at least
for short enough times. A natural question then is what is the evolution of
such a curve, being a natural candidate the solution of the binormal flow.
Recent articles address the issue of existence of solutions when the initial
data is not necessarily small. Gallay and Sverak \cite{GS} considered the
case of a vortex ring with arbitrary circulation $\Gamma $ and in the limit
of vanishing viscosity $\nu $. They show that the vortex ring has a $O((\nu
t)^{\frac{1}{2}})$ thickness and its center translates, for a long period of
time, in the direction orthogonal to the ring with a precise velocity.
Bedrossian, Germain and Harrop-Griffiths \cite{BGH} also remove the
smallness assumption on $\Gamma $ and show existence of global solutions
when the initial data is small perturbation of an infinite filament as well
as local in time solutions for arbitrary initial data.

In this paper we will provide a refined description of solutions to the
Navier-Stokes system with a smooth closed filament $\mathbf{x}_{0}(s,0)$ as
initial data. The filament will evolve with the binormal flow dynamics:%
\begin{equation}
\frac{d\mathbf{x}_{0}(s,t)}{dt}=-\frac{\Gamma }{4\pi }\kappa (s,t)\log (\nu
t)^{\frac{1}{2}}\mathbf{b}(s,t),  \label{bin}
\end{equation}%
where $\Gamma $ is the vortex filament circulation, $\kappa (s,t)$ the local
curvature and $\mathbf{b}(s,t)$ the binormal vector. We will denote by $%
G_{t} $ the corresponding curve at time $t$. Under suitable smoothness
conditions on $\mathbf{x}_{0}(s,0)$ we can find a tube of radius $R$
(independent of $t$) sufficiently small around the filament $\mathbf{x}%
_{0}(s,t)$ and define at any given time coordinates $(s,\rho ,\theta )$ for
a point $\mathbf{x}$ where%
\begin{equation*}
\rho =\frac{dist(\mathbf{x,x}_{0}(s,t))}{(\nu t)^{\frac{1}{2}}},
\end{equation*}%
the coordinate $s$ is the arclength at which the minimum distance between $%
\mathbf{x}$ and $\mathbf{x}_{0}(s,t)$ is achieved and $\theta $ is such that%
\begin{equation*}
\mathbf{e}_{r}\cdot \mathbf{n}(s,t)=\cos \theta ,
\end{equation*}%
where%
\begin{equation*}
\mathbf{e}_{r}=\frac{\mathbf{x-x}_{0}(s,t)}{\left\vert \mathbf{x-x}%
_{0}(s,t)\right\vert }.
\end{equation*}%
We will prove the following result:

\begin{theorem}
\label{teorema}Let $\mathbf{x}_{0}(s,0)\in C^{3}$ and $\mathbf{x}_{0}(s,t)$
a solution to (\ref{bin}) without selfintersections. If $\frac{\Gamma }{\nu }
$ is sufficiently small, there exists $c_{0}>0$ such that for any given $t<%
\frac{c_{0}}{\nu }$ a solution $\boldsymbol{\omega }(x,t)$ to (\ref{vor1}), (%
\ref{vor2}) with initial data (\ref{vor0}) exists. For such solution, there
\ exists $R>0$ independent of $t$ and sufficiently small such that for any $%
\mathbf{x}$ with $dist(\mathbf{x},\mathbf{x}_{0}(s,t))<\frac{R}{2}$ the
vorticity can be written as%
\begin{eqnarray}
\boldsymbol{\omega }(x,t) &=&\frac{1}{(\nu t)}\frac{\Gamma }{4\pi }e^{-\frac{%
\rho ^{2}}{4}}\mathbf{x}_{0s}(s,t)+\frac{1}{(\nu t)^{\frac{1}{2}}}\frac{%
\Gamma \kappa }{8\pi }\rho e^{-\frac{\rho ^{2}}{4}}\left( \cos \theta
\right) \mathbf{x}_{0s}(s,t)  \notag \\
&&+\frac{1}{(\nu t)^{\frac{1}{2}}}\left( \Omega _{1}^{c(2)}(\rho )\left(
\cos \theta \right) +\Omega _{1}^{s(2)}(\rho )\left( \sin \theta \right)
\right) \mathbf{x}_{0s}(s,t)+\widetilde{\boldsymbol{\omega }}(\mathbf{x},t),
\label{expan}
\end{eqnarray}%
with%
\begin{equation*}
\left\vert \Omega _{1}^{s,c(2)}(\rho )\right\vert \leq C\frac{\Gamma ^{2}}{%
\nu }(\rho +\rho ^{2})e^{-\frac{\rho ^{2}}{4}},
\end{equation*}%
and%
\begin{equation*}
\left\Vert \widetilde{\boldsymbol{\omega }}\right\Vert _{L^{2}}^{2}(t)+\nu
\int_{0}^{t}\left\Vert \nabla \widetilde{\boldsymbol{\omega }}\right\Vert
_{L^{2}}^{2}(t^{\prime })dt^{\prime }\leq C\Gamma ^{2}(\nu t)\left\vert \log
(\nu t)\right\vert ^{2},
\end{equation*}%
where $C$ is a suitable constant independent of parameters.
\end{theorem}

This theorem presents the vorticity as the sum of a $O((\nu t)^{-1})$
contribution which is explicit near the moving filament and yields an $%
O((\nu t)^{-\frac{1}{2}})$ norm in $L^{2}$, a semiexplicit $O((\nu t)^{-%
\frac{1}{2}})$ contribution with an $O(1)$ norm in $L^{2}$ and a
contribution that is $O((\nu t)^{\frac{1}{2}}\left\vert \log (\nu
t)\right\vert )$ in $L^{2}\ $and is therefore much smaller than the other
two provided $(\nu t)$ sufficiently small.

Note that under the change $\boldsymbol{\omega }\rightarrow \Gamma \mathbf{%
\omega }$, $\mathbf{v}\rightarrow \Gamma \mathbf{v}$, $t\rightarrow t/\nu $,
equation (\ref{vor2}) is preserved while equation (\ref{vor1}) transforms
into%
\begin{equation*}
\frac{\partial \boldsymbol{\omega }}{\partial t}+\frac{\Gamma }{\nu }\left(
\mathbf{v}\cdot \nabla \boldsymbol{\omega }-\boldsymbol{\omega }\cdot \nabla
\mathbf{v}\right) -\Delta \boldsymbol{\omega }=0,
\end{equation*}%
and the initial data has unit vortex circulation. Hence, the estimates for
the vorticity in Theorem \ref{teorema}, once rescaled by $\Gamma $, depend
solely on the parameter $\frac{\Gamma }{\nu }$ and time appears always under
the product $(\nu t)$. This is also the reason for a timespan $O(\nu ^{-1})$
of the solution.

It is interesting to view the result of the theorem in relation to a
well-known physical situation: the small oscillations of a perturbed vortex
ring. It is well-known since the early works of Levi-Civita \cite{LC} that
small perturbations of a vortex ring under the binormal flow dynamics
undergo oscillatory motion. This is easily seen when taking the formulation
in terms of curvature and torsion:%
\begin{eqnarray}
\kappa _{t^{\prime }} &=&\frac{\Gamma }{\nu }\left( -\kappa \tau
_{s}-2\kappa _{s}\tau \right) ,  \label{eqk} \\
\tau _{t^{\prime }} &=&\frac{\Gamma }{\nu }\left( \left( \frac{\kappa _{ss}}{%
\kappa }-\tau ^{2}\right) _{s}+\kappa \kappa _{s}\right) ,  \label{eqt}
\end{eqnarray}%
with $dt^{\prime }=-\log (\nu t)d(\nu t)$ and linearising around a vortex
ring of radius $R$ so that%
\begin{eqnarray*}
\kappa _{t^{\prime }} &\simeq &-\frac{\Gamma }{\nu }\frac{\tau _{s}}{R}, \\
\tau _{t^{\prime }} &\simeq &\frac{\Gamma }{\nu }\left( R\kappa _{sss}+\frac{%
1}{R}\kappa _{s}\right) ,
\end{eqnarray*}%
that can be combined to provide%
\begin{equation*}
\kappa _{t^{\prime }t^{\prime }}=-\left( \frac{\Gamma }{\nu }\right)
^{2}\left( \kappa _{ssss}+\frac{1}{R^{2}}\kappa _{ss}\right) ,
\end{equation*}%
with solutions of the form $\kappa (s,t^{\prime })=e^{i\frac{\Gamma /\nu }{%
R^{2}}\sqrt{n^{4}-n^{2}}t^{\prime }}\cos \left( \frac{ns}{R}\right) $.
Hence, a whole period of oscillation is $T=O(R^{2}/\Gamma )$. According to
our theorem, the expansion is valid for times such that $(\nu t)^{\frac{1}{2}%
}\ll R$ provided $\frac{\Gamma }{\nu }$ is sufficiently small but order
unit. Therefore $t\ll \frac{R^{2}}{\nu }=O(R^{2}/\Gamma )$ and we can expect
the expansion (\ref{expan}) to be valid for a significant fraction of a
whole period.

The solutions constructed this way are small perturbations of circles and
therefore do not have self intersections. Another examples of curves without
self intersections are those obtained by Kida in \cite{K} and \cite{K2}. In
this case the evolution is given by a rigid motion. For closed curves Kida
obtains a three parameters family of examples. One of the parameters is the
scaling. The other two represent the number of rotations around a fixed axis
and the number of winds around a given circle. Our results can be applied to
these examples for a time which depends on the size of the parameters. In
future work we want to explore how we can use our argument to optimize this
time of existence.

Another way to obtain examples of curves without selfintersections is using
the so called Galilean transformations that leaves invariant the set of
solutions of the 1d cubic NLS. Since the work of Hasimoto \cite{H} is well
known the relation between the binormal flow and NLS. For example the
helices can be obtained applying this transformation to circles. This raises
the natural question about how our results extend to non-closed curves like
helices. This is related to the existence of a tubular neighbourhood such
that the change of variables to \textquotedblleft
cylindrical\textquotedblright\ coordinates can be implemented along the
evolution. For example, if the curve is a small perturbation of straight
line the method can be extended without difficulty. The case of helices is
more delicate if one wants to obtain good dependence with respect to the
curvature and the pitch.

A relevant situation of non-closed curves is the one given by the
self-similar solutions of the BF, see for example \cite{VG} for a full
classification of them together with issues concerning the existence of
self-intersections and \cite{BV} where the stability of such solutions is
addressed. Although due to the self similarity the tubular neighbourhood has
to depend on time, one can just work with the associated time independent
equation. This particular case will be studied in a forthcoming paper.

In the process of deducing the formula (\ref{expan}) we compute the $O(1)$
corrections to the law of motion of the filament beyond the $O(\log (\nu t))$
term given by (\ref{bin}). More precisely, one can write, up to $O(\Gamma
^{2}/\nu )$ error,
\begin{equation}
\frac{d\mathbf{x}_{0}}{dt}=\left( -\frac{\Gamma }{4\pi }\kappa \log (\nu t)^{%
\frac{1}{2}}+\frac{\Gamma \kappa }{8\pi }\gamma -\frac{\Gamma \kappa }{4\pi }%
\right) \mathbf{b}+\mathbf{v}^{\ast }(s),  \label{vlaw}
\end{equation}%
where%
\begin{equation}
\mathbf{v}^{\ast }(s,t)=\frac{\Gamma }{4\pi }\lim_{\varepsilon \rightarrow
0}\left( \int_{G_{t}\backslash \left[ -\varepsilon ,\varepsilon \right] }%
\frac{\mathbf{x}_{0s}(s^{\prime },t)\times (\mathbf{x}_{0}(s,t)-\mathbf{x}%
_{0}(s^{\prime },t))}{\left\vert \mathbf{x}_{0}(s,t)-\mathbf{x}%
_{0}(s^{\prime },t)\right\vert ^{3}}ds^{\prime }+\kappa (\log \varepsilon )%
\mathbf{b}(s,t)\right) ,  \label{vstar}
\end{equation}%
with the first term at the right hand side of (\ref{vstar}) bounded and
regular and $\gamma $ being Euler's constant. Note that the first two terms
in (\ref{expan}) imply a displacement $\kappa (\nu t)^{\frac{1}{2}}$of the
maximum of the gaussian vorticity profile in the negative normal direction.
In the case of a circular vortex this implies that the vortex ring's radius
tends to increase.

The paper is organized as follows. In section 2, we will present some
preliminary results on the regularity of vortex filaments and the existence
of tubular neighborhoods that will be used in the next sections. Section 3
will be devoted to discuss the structure of Navier-Stokes system in moving
frames. In sections 4 and 5 we will construct approximate solutions to
Navier-Stokes system as deformed Lamb-Oseen vortices. The construction of
the full solution is done in section 7 after obtaining some necessary
estimates in section 6 that are based on the previously constructed
approximate solution. Finally, in an appendix we deduce a novel formula for
the desingularization of the Biot-Savart integral that is used throughout
the paper.

\section{\protect\bigskip Preliminaries}

In this section we will discuss two of the main ingredients for our further
analysis: 1) the existence of a vortex filament following the binormal flow
dynamics with a high degree of regularity provided the initial data is
sufficiently regular and selfintersections do not occur in the time interval
$(0.T)$ and, 2) the existence of a tubular neighborhood around a regular
vortex filament. In the next lemma we will denote\ by $%
C((0.T),W(T^{1},S^{2}))$ the set of functions whose $W$-norm as functions
from $s\in T^{1}$ (that is, in a fixed period) into $S^{2}$ depends
continuously on time in the interval $(0.T).$

\begin{lemma}
Let $\mathbf{x}_{0}(s,0)$ be an initial curve such that $\mathbf{t}(s,0)=%
\frac{d\mathbf{x}_{0}(s,0)}{ds}\in H^{k}(T^{1},S^{2})$ with $k\geq 3$, then
\begin{eqnarray*}
\mathbf{t}(t,s) &\in &C(\mathbb{R},H^{k}(T^{1},S^{2})), \\
\left( \mathbf{n}(s,t),\mathbf{b}(s,t)\right) &\in
&C((0.T),C^{k-2}(T^{1},S^{2})), \\
\kappa (s,t),\tau (s,t) &\in &C((0.T),C^{k-3}(T^{1},\mathbb{R})).
\end{eqnarray*}
\end{lemma}

\textbf{Proof. }In \cite{JS} (see also \cite{SSB}) it was shown that $%
\mathbf{t}(s,t)\in C(%
\mathbb{R}
,H^{k}(T^{1},S^{2}))$ for initial data in $H^{k}(T^{1},S^{2})$. Hence, by
Sobolev embeddings, $\mathbf{x}(s,t)$ $\in C((0.T),C^{k}(T^{1},S^{2}))$ if $%
\mathbf{x}(s,0)\in H^{k+1}(T^{1},S^{2})$. In order to extend these
regularity results to curvature and torsion, we define the vectors \ $%
\left\{ \mathbf{e}_{1}(s,t),\mathbf{e}_{2}(s,t)\right\} $ in the normal
plane to $\mathbf{x}(s,t)$ and satisfying

\begin{eqnarray}
\mathbf{e}_{1,s} &=&-(\mathbf{t}_{s}\cdot \mathbf{e}_{1})\mathbf{t},
\label{e11} \\
\mathbf{e}_{2,s} &=&-(\mathbf{t}_{s}\cdot \mathbf{e}_{2})\mathbf{t},
\label{e12}
\end{eqnarray}%
where $\mathbf{t}=\mathbf{x}_{s}$ is the tangent vector. These are first
order ODEs defining the vectors $\mathbf{e}_{1},\mathbf{e}_{2}$ and hence
having the same regularity at $\mathbf{t}$, i.e. $\mathbf{e}_{1}(s,t),%
\mathbf{e}_{2}(s,t)\in C((0.T),C^{k-1}(T^{1},S^{2}))$. Now, since%
\begin{equation*}
\mathbf{t}_{s}=(\mathbf{t}_{s}\cdot \mathbf{e}_{1})\mathbf{e}_{1}+(\mathbf{t}%
_{s}\cdot \mathbf{e}_{2})\mathbf{e}_{2}=-\left( \mathbf{e}_{1,s}\cdot
\mathbf{t}\right) \mathbf{e}_{1}-\left( \mathbf{e}_{2,s}\cdot \mathbf{t}%
\right) \mathbf{e}_{2}\equiv \alpha \mathbf{e}_{1}+\beta \mathbf{e}_{2},
\end{equation*}%
where $\alpha ,\beta \in C((0.T),C^{k-2}(T^{1},%
\mathbb{R}
))$, we can write the complex%
\begin{equation*}
\alpha +i\beta =\kappa (s,t)e^{i\int_{0}^{z}\tau (s^{\prime },t)ds^{\prime
}},
\end{equation*}%
implying $\kappa \in C((0.T),C^{k-2}(T^{1},%
\mathbb{R}
))$ and $\tau \in C((0.T),C^{k-3}(T^{1},%
\mathbb{R}
))$.

Finally, we discuss the regularity of $\mathbf{n}(s,t),\mathbf{b}(s,t)$. The
Frenet-Serret system is

\begin{equation*}
\frac{d}{ds}\left(
\begin{array}{c}
\mathbf{t} \\
\mathbf{n} \\
\mathbf{b}%
\end{array}%
\right) =\left(
\begin{array}{ccc}
0 & \kappa & 0 \\
-\kappa & 0 & \tau \\
0 & -\tau & 0%
\end{array}%
\right) \left(
\begin{array}{c}
\mathbf{t} \\
\mathbf{n} \\
\mathbf{b}%
\end{array}%
\right) ,
\end{equation*}%
which, after taking an $s$-derivative, leads to%
\begin{eqnarray*}
\frac{d^{2}}{ds^{2}}\left(
\begin{array}{c}
\mathbf{t} \\
\mathbf{n} \\
\mathbf{b}%
\end{array}%
\right) &=&\left(
\begin{array}{ccc}
0 & \kappa & 0 \\
-\kappa & 0 & \tau \\
0 & -\tau & 0%
\end{array}%
\right) \frac{d}{ds}\left(
\begin{array}{c}
\mathbf{t} \\
\mathbf{n} \\
\mathbf{b}%
\end{array}%
\right) +\left(
\begin{array}{ccc}
0 & \kappa _{s} & 0 \\
-\kappa _{s} & 0 & \tau _{s} \\
0 & -\tau _{s} & 0%
\end{array}%
\right) \left(
\begin{array}{c}
\mathbf{t} \\
\mathbf{n} \\
\mathbf{b}%
\end{array}%
\right) \\
&=&\left[ \left(
\begin{array}{ccc}
0 & \kappa & 0 \\
-\kappa & 0 & \tau \\
0 & -\tau & 0%
\end{array}%
\right) ^{2}+\left(
\begin{array}{ccc}
0 & \kappa _{s} & 0 \\
-\kappa _{s} & 0 & \tau _{s} \\
0 & -\tau _{s} & 0%
\end{array}%
\right) \right] \left(
\begin{array}{c}
\mathbf{t} \\
\mathbf{n} \\
\mathbf{b}%
\end{array}%
\right) ,
\end{eqnarray*}%
and after successive derivatives:%
\begin{equation*}
\frac{d^{m}}{ds^{m}}\left(
\begin{array}{c}
\mathbf{t} \\
\mathbf{n} \\
\mathbf{b}%
\end{array}%
\right) =A\left[ \kappa ,\tau ,...,\kappa ^{(m-1)},\tau ^{(m-1)}\right]
\left(
\begin{array}{c}
\mathbf{t} \\
\mathbf{n} \\
\mathbf{b}%
\end{array}%
\right) .
\end{equation*}%
\bigskip Hence \bigskip $\left( \mathbf{n}(s,t),\mathbf{b}(s,t)\right) \in
C((0.T),C^{k-2}(T^{1},S^{2}))$ since $(\kappa (s,r),\tau (s,t))$ is $%
C((0.T),C^{k-3}(T^{1},%
\mathbb{R}
))$.

Next, we prove the following auxiliary Lemma showing that a vortex filament
can be extended as a vortex tube of a certain radius $R$:

\begin{lemma}
Given a vortex filament $\mathbf{x}_{0}(s,t)$ $\in
C((0.T),C^{2}(T^{1},S^{2}))$, there exists an $R>0$ with%
\begin{equation*}
R<\frac{1}{2\max_{s,t}\left\vert \kappa (s,t)\right\vert },
\end{equation*}%
such that at a given time $t$ and for any $\mathbf{x}$ with $dist(\mathbf{x,x%
}_{0}(s,t))<R$ there exists a unique $s_{0}$ such that%
\begin{equation*}
(\mathbf{x}-\mathbf{x}_{0}(s_{0},t))\bot \mathbf{t}(t,s_{0}),
\end{equation*}%
and therefore
\begin{equation*}
dist(\mathbf{x,x}_{0}(s,t))=dist(\mathbf{x,x}_{0}(s_{0},t)).
\end{equation*}
\end{lemma}

We take two points $\mathbf{x}_{1}$ and $\mathbf{x}_{2}$ that lie in the
normal planes to $\mathbf{x}_{0}(s_{1})$ and $\mathbf{x}_{0}(s_{2})$
respectively. Hence we can write
\begin{eqnarray*}
\mathbf{x}_{1} &=&\mathbf{x}_{0}(s_{1})+\mathbf{\rho }_{1}, \\
\mathbf{x}_{2} &=&\mathbf{x}_{0}(s_{2})+\mathbf{\rho }_{2},
\end{eqnarray*}%
where we can write, using polar coordinates in the normal planes,%
\begin{equation*}
\mathbf{\rho }_{i}=r_{i}\cos \theta _{i}\mathbf{e}_{1}(s_{i})+r_{i}\sin
\theta _{i}\mathbf{e}_{2}(s_{i}).
\end{equation*}%
Notice that the displacement $\delta \mathbf{\rho }=\mathbf{\rho }_{2}-%
\mathbf{\rho }_{1}$ is given, at leading order in $\delta r=r_{2}-r_{1}$, $%
\delta \theta =\theta _{2}-\theta _{1}$, $\delta \mathbf{e}_{1,2}(s)=\mathbf{%
e}_{1,2}(s_{2})-\mathbf{e}_{1,2}(s_{1})$ we have%
\begin{eqnarray*}
\delta \mathbf{\rho } &=&\left( \cos \theta _{1}\mathbf{e}_{1}(s_{1})+\sin
\theta _{1}\mathbf{e}_{2}(s_{1})\right) \delta r \\
&&+r_{1}\left( -\sin \theta _{1}\mathbf{e}_{1}(s_{1})+\cos \theta _{1}%
\mathbf{e}_{2}(s_{1})\right) \delta \theta \\
&&+r_{1}\cos \theta _{1}\delta \mathbf{e}_{1}(s)+r_{1}\sin \theta _{1}\delta
\mathbf{e}_{2}(s).
\end{eqnarray*}%
Moreover,
\begin{eqnarray*}
\delta \mathbf{e}_{1}(s) &=&-(\mathbf{t}_{s}\cdot \mathbf{e}_{1})\mathbf{t}%
\delta s+O((\delta s)^{2}), \\
\delta \mathbf{e}_{2}(s) &=&-(\mathbf{t}_{s}\cdot \mathbf{e}_{2})\mathbf{t}%
\delta s+O((\delta s)^{2}),
\end{eqnarray*}%
and hence%
\begin{equation*}
\delta \mathbf{\rho }=\delta r\frac{\mathbf{\rho }_{1}}{\left\vert \mathbf{%
\rho }_{1}\right\vert }+r_{1}\delta \theta \frac{\mathbf{\rho }_{1}^{\perp }%
}{\left\vert \mathbf{\rho }_{1}\right\vert }-(\mathbf{\rho }_{1}\cdot
\mathbf{n})\kappa (s_{1})\mathbf{t}(s_{1})\delta s,
\end{equation*}%
so that%
\begin{equation*}
\left\vert \delta \mathbf{x}\right\vert ^{2}=\left( \delta r\right)
^{2}+r_{1}^{2}(\delta \theta )^{2}+(1-r_{1}\cos \theta _{1}\kappa
(s_{1}))^{2}(\delta s)^{2}+O(\delta ^{3}),
\end{equation*}%
where the error term is cubic in the displacements. The minimum distance $%
\left\vert \delta \mathbf{x}\right\vert $ is then achieved when $\delta r=0$%
, $\delta \theta =0$ and then%
\begin{equation*}
\left\vert \delta \mathbf{x}\right\vert ^{2}=(1-r_{1}\cos \theta _{1}\kappa
(s_{1}))^{2}(\delta s)^{2}+O(\delta ^{3}),
\end{equation*}%
which is positive provided%
\begin{equation*}
r_{1}\kappa (s_{1})<1.
\end{equation*}%
Hence, we can find a $\delta s$ sufficiently small and%
\begin{equation*}
R<\frac{1}{2\max \left\vert \kappa \right\vert },
\end{equation*}%
so that the discs of radius $R$ normal to the filament at $s_{1}$ and $%
s_{1}+\delta s$ do not intersect. The argument can be repeated for each
portion between $s_{N}=s_{1}+N\delta s$ and $s_{N+1}=s_{1}+(N+1)\delta s$ ($%
N\in
\mathbb{Z}
$) and construct a tube of radius $R$ around the filament where each point $%
\mathbf{x}$ is biunivocally defined by a set of coordinates $(s,r,\theta )$.
This concludes the proof of the Lemma.

\begin{figure}[t]
\includegraphics[width=1.2\textwidth]{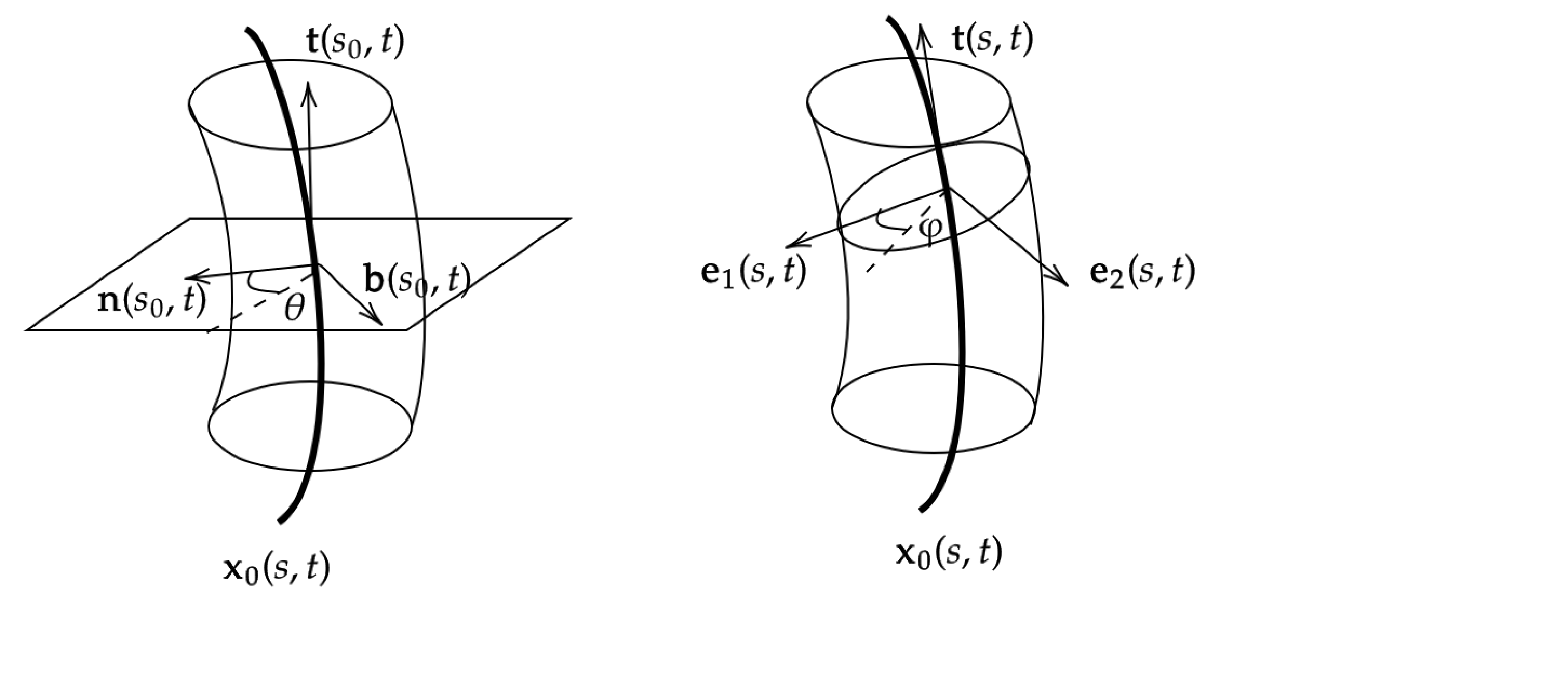}
\caption{Sketch of the local Frenet-Serret frame (left) and the filament
parallel (or natural) frame (right).}
\label{fig1}
\end{figure}

By the previous Lemma, for any $\mathbf{x}$ inside the tube or radius $R$
there exists a unique value of the arclength parameter $s_{0}$ such that%
\begin{equation*}
\mathbf{x}=\mathbf{x}_{0}(s_{0},t)+\boldsymbol{\rho },
\end{equation*}%
where $\boldsymbol{\rho }$ lies in the plane normal to $\mathbf{x}_{0}(s,t)$
at $s=s_{0}$. The vector $\boldsymbol{\rho }$ can be written in different
systems of coordinates of the normal plane: i) the Frenet-Serret frame (see
figure \ref{fig1})), so that%
\begin{equation}
\mathbf{x}=\mathbf{x}_{0}(s_{0},t)+x_{n}\mathbf{n}(s_{0},t)+x_{b}\mathbf{b}%
(s_{0},t),  \label{lfs}
\end{equation}%
with $(\mathbf{n}(s_{0},t),\mathbf{b}(s_{0},t))$ the normal and binormal
vectors to $\mathbf{x}_{0}(s,t)$ at $s=s_{0}$; ii) one can rewrite the
coordinates $(x_{n},x_{b})$ in a polar system such that%
\begin{eqnarray}
x_{n} &=&r\cos \theta ,  \label{xn} \\
x_{b} &=&r\sin \theta ,  \label{xb}
\end{eqnarray}%
and then
\begin{equation*}
\mathbf{x}=\mathbf{x}_{0}(s_{0},t)+r\mathbf{e}_{r}(s_{0},t),
\end{equation*}%
with%
\begin{eqnarray*}
\mathbf{e}_{r} &=&\left( \sin \theta \right) \mathbf{b}+\left( \cos \theta
\right) \mathbf{n}, \\
\mathbf{e}_{\theta } &=&\left( \cos \theta \right) \mathbf{b}-\left( \sin
\theta \right) \mathbf{n},
\end{eqnarray*}%
and, finally iii) one can use the so-called parallel or natural frame (see
figure \ref{fig1}) where
\begin{equation}
\mathbf{x}=\mathbf{x}_{0}(s_{0},t)+x_{1}\mathbf{e}_{1}(s_{0},t)+x_{2}\mathbf{%
e}_{2}(s_{0},t),  \label{xpara}
\end{equation}%
where $\left\{ \mathbf{e}_{1}(s,t),\mathbf{e}_{2}(s,t)\right\} $ were
defined above by means of (\ref{e11}), (\ref{e12}). Hence, the points $%
\mathbf{x}$ inside the tube of radius $R$ can be univocally labelled with
coordinates $(x_{n},x_{b},s)$ (Frenet-Serret frame), $(r,\theta ,s)$ (polar
coordinates in Frenet-Serret frame) or $(x_{1},x_{2},s)$ (parallel frame).

\bigskip Inside a tube of radius $R$ we can define a vector field $\mathbf{J}
$ in the form%
\begin{equation*}
\mathbf{J}=A_{1}(s,r,\theta )\mathbf{t}(s)+A_{2}(s,r,\theta )\mathbf{n}%
(s)+A_{3}(s,r,\theta )\mathbf{b}(s),
\end{equation*}%
and extend it to $%
\mathbb{R}
^{3}$ by using a cutoff function $\eta _{R}(dist(\mathbf{x},\mathbf{x}%
_{0}(s,t))$ so that%
\begin{eqnarray}
\eta _{R}(dist(\mathbf{x},\mathbf{x}_{0}(s,t)) &\in &C^{\infty },
\label{co0} \\
\eta _{R}(dist(\mathbf{x},\mathbf{x}_{0}(s)) &=&0\text{ outside the tube of
radius }R,  \label{co2a} \\
\eta _{R}(dist(\mathbf{x},\mathbf{x}_{0}(s)) &=&1\text{ inside the tube of
radius }R/2.  \label{co2}
\end{eqnarray}%
The vector field
\begin{equation*}
\mathbf{J}_{\eta }=\eta _{R}(dist(\mathbf{x},\mathbf{x}_{0}(s)))\mathbf{J},
\end{equation*}%
is identical to $\mathbf{J}$ inside the tube of radius $R/2$ and has the
same regularity (as a function in $%
\mathbb{R}
^{3}$) as $\mathbf{J}$.

We can also define moving tubes around a filament $\mathbf{x}_{0}(s,t)$ (for
$t\in \left[ 0,T\right] $) with uniform radius $R$ such that%
\begin{equation*}
R<\frac{1}{2\max_{s,t}\left\vert \kappa (s,t)\right\vert },
\end{equation*}%
and the corresponding moving cutoff function $\eta _{R}(dist(\mathbf{x},%
\mathbf{x}_{0}(s,t)))$.

\section{The Frenet-Serret frame and the parallel frame}

Navier-Stokes equations are not invariant under coordinate changes such as
those described in the previous section (the Frenet-Serret frame and the
parallel frame coordinates). The fact that reference frames do translate and
rotate in a nonuniform way lead to the appearance of the so-called
"fictitious forces" in the new reference frame such as Coriolis and
centrifugal forces. Nevertheless, for a vortex filament, we will see in
section \ref{secf} that at leading order those forces are negligible, so
that Navier-Stokes system still describes (at leading order in $t$) the
velocity and vorticity fields. Central to the analysis will be to estimate
the velocity at which vectors $(\mathbf{n},\mathbf{b},\mathbf{t})$ and $(%
\mathbf{e}_{1},\mathbf{e}_{2})$ translate and rotate. This is the goal of
the present section.

In order to understand how these fictitious forces appear, let us consider
first a local (in $s_{0}$) Frenet-Serret frame (see figure \ref{fig1}). The
coordinates in that frame will be $x^{\prime }=(x_{n},x_{b},z)$ such that%
\begin{equation}
\mathbf{x}=\mathbf{x}_{0}(s_{0},t)+x_{n}\mathbf{n}(s_{0},t)+x_{b}\mathbf{b}%
(s_{0},t)+z\mathbf{t}(s_{0},t).  \label{lfs2}
\end{equation}%
If we consider now a trajectory $\mathbf{x}(t)$ following a velocity field $%
\mathbf{v}(\mathbf{x},t)$ then%
\begin{eqnarray}
\mathbf{v}(\mathbf{x}(t),t) &=&\frac{d\mathbf{x}(t)}{dt}  \notag \\
&=&\frac{d\mathbf{x}_{0}(s_{0},t)}{dt}+\frac{dx_{n}(t)}{dt}\mathbf{n}%
(s_{0},t)+\frac{dx_{b}(t)}{dt}\mathbf{b}(s_{0},t)+\frac{dz(t)}{dt}\mathbf{t}%
(s_{0},t)  \notag \\
&&+x_{n}(t)\frac{d\mathbf{n}(s_{0},t)}{dt}+x_{b}(t)\frac{d\mathbf{b}(s_{0},t)%
}{dt}+z(t)\frac{d\mathbf{t}(s_{0},t)}{dt}  \notag \\
&=&\frac{d\mathbf{x}_{0}(s_{0},t)}{dt}+\mathbf{v}^{\prime }(\mathbf{x}%
^{\prime }(t),t)+\mathbf{\varpi }\times \mathbf{x}^{\prime }(t),
\label{velfs}
\end{eqnarray}%
which relates the velocities $\mathbf{v}(\mathbf{x},t)$ in the fixed frame
and $\mathbf{v}^{\prime }(\mathbf{x}^{\prime },t)$ in the local
Frenet-Serret frame. The vector $\mathbf{\varpi }$ is the instantaneous
angular velocity and its explicit expression and relation to $\frac{d}{dt}(%
\mathbf{n},\mathbf{b},\mathbf{t})$ will be described below (formulas (\ref%
{d0p}) and (\ref{d1})-(\ref{d3})). Concerning the various terms that appear
in the vorticity equation, we compute first the convective derivative of the
vorticity (using the notation $(\mathbf{u}_{1},\mathbf{u}_{2},\mathbf{u}%
_{3})=(\mathbf{n},\mathbf{b},\mathbf{t})$ and Einstein's summation
convention):
\begin{eqnarray*}
&&\frac{d\boldsymbol{\omega }(\mathbf{x}(t),t)}{dt} \\
&=&\frac{d\boldsymbol{\omega }^{\prime }(\mathbf{x}^{\prime }(t),t)}{dt}=%
\frac{\partial \omega _{i}}{\partial t}\mathbf{u}_{i}(t)+\omega _{i}\frac{d%
\mathbf{u}_{i}(t)}{dt}+\frac{d\mathbf{x}^{\prime }(t)}{dt}\cdot \nabla
\mathbf{\omega }^{\prime }(\mathbf{x}^{\prime }(t),t) \\
&=&\frac{\partial \omega _{i}}{\partial t}\mathbf{u}_{i}(t)+\omega _{i}\frac{%
d\mathbf{u}_{i}(t)}{dt}+\left( \mathbf{v}^{\prime }(\mathbf{x}^{\prime
}(t),t)\cdot \nabla \omega _{i}(\mathbf{x}^{\prime }(t),t)\right) \mathbf{u}%
_{i}(t) \\
&&+x_{i}(t)\frac{d\mathbf{u}_{i}(t)}{dt}\cdot \nabla \boldsymbol{\omega }%
^{\prime }(\mathbf{x}^{\prime }(t),t) \\
&=&\left( \frac{\partial \omega _{i}}{\partial t}+\mathbf{v}^{\prime }(%
\mathbf{x}^{\prime }(t),t)\cdot \nabla \omega _{i}(\mathbf{x}^{\prime
}(t),t)\right) \mathbf{u}_{i}(t)+\mathbf{\varpi }\times \boldsymbol{\omega }%
^{\prime }(\mathbf{x}^{\prime }(t),t) \\
&&+\left( \mathbf{\varpi }\times \mathbf{x}^{\prime }(t)\right) \cdot \nabla
\boldsymbol{\omega }^{\prime }(\mathbf{x}^{\prime }(t),t),
\end{eqnarray*}%
where the first term at the right hand side is the convective derivative in
the Frenet-Serret coordinate frame and the last two terms are corrections
that amount to "fictitious forces" for an observer in the moving frame. The
same can be said of the vortex stretching term:%
\begin{equation*}
\boldsymbol{\omega }(\mathbf{x},t)\cdot \nabla \mathbf{v}(\mathbf{x},t)=%
\mathbf{\omega }^{\prime }(\mathbf{x}^{\prime },t)\cdot \nabla \mathbf{v}%
^{\prime }(\mathbf{x}^{\prime },t)+\boldsymbol{\omega }^{\prime }(\mathbf{x}%
^{\prime },t)\cdot \nabla \left( \mathbf{\varpi }\times \mathbf{x}^{\prime
}\right) ,
\end{equation*}%
with an extra term at the right hand side while the viscous term remains
unchanged under the coordinate change since the Laplacian operator is
invariant under translations and rotations.

We compute and estimate next the vector $\mathbf{\varpi }$. Based on the
binormal flow equations%
\begin{equation*}
\frac{d\mathbf{x}}{dt}=-\frac{\kappa \Gamma }{4\pi }\left( \log (\nu t)^{%
\frac{1}{2}}\right) \mathbf{b},
\end{equation*}%
we introduce%
\begin{equation*}
t^{\prime }=\frac{\Gamma }{4\pi }\int_{0}^{t}\log (\nu \tau )^{\frac{1}{2}%
}d\tau ,
\end{equation*}%
and compute%
\begin{eqnarray*}
\frac{d\mathbf{t}}{dt^{\prime }} &=&\kappa _{s}\mathbf{b}+\kappa \mathbf{b}%
_{s}=\kappa _{s}\mathbf{b}-\kappa \tau \mathbf{n=-}\kappa _{s}\mathbf{n}%
\times \mathbf{t-}\kappa \tau \mathbf{b}\times \mathbf{t}, \\
\frac{d(\kappa \mathbf{n})}{dt^{\prime }} &=&\kappa _{ss}\mathbf{b}+\kappa
_{s}\mathbf{b}_{s}-(\kappa \tau )_{s}\mathbf{n}-\kappa \tau \left( -\kappa
\mathbf{t}+\tau \mathbf{b}\right) ,
\end{eqnarray*}%
so that%
\begin{equation*}
\frac{d\mathbf{n}}{dt^{\prime }}=-\frac{\kappa _{t^{\prime }}+\kappa
_{s}\tau +(\kappa \tau )_{s}}{\kappa }\mathbf{n+}\kappa \tau \mathbf{t}+%
\frac{\kappa _{ss}-\kappa \tau ^{2}}{\kappa }\mathbf{b}.
\end{equation*}%
Since%
\begin{equation*}
0=\mathbf{n}\cdot \frac{d\mathbf{n}}{dt^{\prime }}=-\frac{\kappa _{t^{\prime
}}+\kappa _{s}\tau +(\kappa \tau )_{s}}{\kappa },
\end{equation*}%
we have%
\begin{eqnarray*}
\frac{d\mathbf{n}}{dt^{\prime }} &=&\kappa \tau \mathbf{t}+\frac{\kappa
_{ss}-\kappa \tau ^{2}}{\kappa }\mathbf{b} \\
&=&\frac{\kappa _{ss}-\kappa \tau ^{2}}{\kappa }\mathbf{t}\times \mathbf{n-}%
\kappa \tau \mathbf{b}\times \mathbf{n}.
\end{eqnarray*}%
Notice that, by equation (\ref{eqt}),
\begin{equation}
\frac{\kappa _{ss}-\kappa \tau ^{2}}{\kappa }=\int \tau _{t^{\prime }}ds-%
\frac{1}{2}\kappa ^{2},  \label{d0}
\end{equation}%
which is bounded.

Note now%
\begin{eqnarray*}
\frac{d\mathbf{b}}{dt^{\prime }} &=&\mathbf{t}\times \frac{d\mathbf{n}}{%
dt^{\prime }}+\frac{d\mathbf{t}}{dt^{\prime }}\times \mathbf{n}=\frac{\kappa
_{ss}-\kappa \tau ^{2}}{\kappa }\mathbf{t}\times \mathbf{b}+\kappa _{s}%
\mathbf{b}\times \mathbf{n} \\
&=&\frac{\kappa _{ss}-\kappa \tau ^{2}}{\kappa }\mathbf{t}\times \mathbf{b}%
-\kappa _{s}\mathbf{n}\times \mathbf{b}.
\end{eqnarray*}%
Hence, by introducing%
\begin{equation}
\boldsymbol{\varpi }=\frac{\Gamma }{2\pi }\log (\nu t)\left[ \frac{\kappa
_{ss}-\kappa \tau ^{2}}{\kappa }\mathbf{t}-\kappa _{s}\mathbf{\mathbf{n}-}%
\kappa \tau \mathbf{b}\right] ,  \label{d0p}
\end{equation}%
we have%
\begin{eqnarray}
\frac{d\mathbf{t}}{dt} &=&\boldsymbol{\varpi }\times \mathbf{t},  \label{d1}
\\
\frac{d\mathbf{n}}{dt} &=&\boldsymbol{\varpi }\times \mathbf{n},  \label{d2}
\\
\frac{d\mathbf{b}}{dt} &=&\boldsymbol{\varpi }\times \mathbf{b},  \label{d3}
\end{eqnarray}%
where $\boldsymbol{\varpi }$ plays the same role as the Darboux vector but
for the binormal flow. Equations (\ref{d1})-(\ref{d3}) imply that the
Frenet-Serret frame rotates with instantaneous angular velocity $\boldsymbol{%
\varpi }$. Note that by (\ref{d0p}) we have
\begin{equation*}
\left\vert \boldsymbol{\varpi }\right\vert \leq C\Gamma \left\vert \log (\nu
t)\right\vert .
\end{equation*}

In the case of the parallel (or natural) frame coordinates $(x_{1},x_{2},s)$
(see figure \ref{fig1}) defined by means of (\ref{xpara}), the vectors $(%
\mathbf{e}_{1}(s,t),\mathbf{e}_{2}(s,t))$ form a basis for vectors in the
plane normal to $x_{0}(s,t)$ and relate to the normal and binormal vectors
by means of the relations%
\begin{eqnarray}
\mathbf{n}(s,t) &=&\cos \left( \theta _{0}(s,t)\right) \mathbf{e}%
_{1}(s,t)+\sin \left( \theta _{0}(s,t)\right) \mathbf{e}_{2}(s,t),
\label{fo1} \\
\mathbf{b}(s,t) &=&-\sin \left( \theta _{0}(s,t)\right) \mathbf{e}%
_{1}(s,t)+\cos \left( \theta _{0}(s,t)\right) \mathbf{e}_{2}(s,t),
\label{fo2}
\end{eqnarray}%
where%
\begin{equation*}
\theta _{0}(s,t)=\int_{s_{0}}^{s}\tau (s^{\prime },t)ds^{\prime },
\end{equation*}%
and $s_{0}$ correspond to an arbitrary point along the filament. We can also
combine $(x_{1},x_{2})$ in a system of polar coordinates so that%
\begin{eqnarray*}
x_{1} &=&r\cos \varphi , \\
x_{2} &=&r\sin \varphi .
\end{eqnarray*}%
Hence, one can rewrite%
\begin{eqnarray*}
x_{n} &=&r\cos \left( \varphi -\theta _{0}(s,t)\right) , \\
x_{b} &=&r\sin \left( \varphi -\theta _{0}(s,t)\right) .
\end{eqnarray*}%
In this way, the system $(r,\varphi ,s)$ forms an orthogonal system with
unit vectors%
\begin{equation}
(\mathbf{e}_{r},\mathbf{e}_{\varphi },\mathbf{e}_{s})=(\cos \varphi \mathbf{e%
}_{1}+\sin \varphi \mathbf{e}_{2},-\sin \varphi \mathbf{e}_{1}+\cos \varphi
\mathbf{e}_{2},\mathbf{t}),  \label{erefes}
\end{equation}%
once we have untwisted the origin of the angular coordinate according to the
torsion of $\mathbf{x}_{0}(s,t)$. The associated scale factors are
\begin{equation}
h_{r}=1,h_{\varphi }=1,h_{s}=1-\kappa r\cos \left( \varphi -\theta
_{0}(s,t)\right) ,  \label{scale1}
\end{equation}%
and hence the volume element $dV_{i}$ relates to the volume element $dV_{fs}$
in Frenet-Serret system by%
\begin{equation}
dV_{fs}=\left( 1-\kappa r\cos \left( \varphi -\theta _{0}(s,t)\right)
\right) dV_{i}.  \label{scale2}
\end{equation}

Analogously to (\ref{velfs}) we can relate the velocity $\mathbf{v}(\mathbf{x%
},t)$ in the fixed frame and the velocity $\widetilde{\mathbf{v}}(\widetilde{%
\mathbf{x}},t)$ the parallel frame:%
\begin{eqnarray}
&&\mathbf{v}(\mathbf{x}(t),t)  \notag \\
&=&\mathbf{x}_{0,t}(s(t),t)+s^{\prime }(t)\mathbf{t}(s(t),t)+x_{1,t}(t)%
\mathbf{e}_{1}(s(t),t)+x_{2,t}(t)\mathbf{e}_{2}(s(t),t)  \notag \\
&&+x_{1}(t)\left( \mathbf{e}_{1,t}(s(t),t)+s^{\prime }(t)\mathbf{e}%
_{1,s}(s(t),t)\right) +x_{2}(t)\left( \mathbf{e}_{2,t}(s(t),t)+s^{\prime }(t)%
\mathbf{e}_{2,s}(s(t),t)\right)  \notag \\
&=&\mathbf{x}_{0,t}(s(t),t)+v_{s}(t)\mathbf{t}(s(t),t)+v_{1}(t)\mathbf{e}%
_{1}(s(t),t)+v_{2}(t)\mathbf{e}_{2}(s(t),t)  \notag \\
&&+x_{1}(t)\left( \mathbf{e}_{1,t}(s(t),t)+v_{s}(t)\mathbf{e}%
_{1,s}(s(t),t)\right) +x_{2}(t)\left( \mathbf{e}_{2,t}(s(t),t)+v_{s}(t)%
\mathbf{e}_{2,s}(s(t),t)\right)  \notag \\
&=&\frac{d\mathbf{x}_{0}(s,t)}{dt}+\widetilde{\mathbf{v}}(\widetilde{\mathbf{%
x}}(t),t)+x_{1}(t)\left( \mathbf{e}_{1,t}(s(t),t)+v_{s}(t)\mathbf{e}%
_{1,s}(s(t),t)\right)  \notag \\
&&+x_{2}(t)\left( \mathbf{e}_{2,t}(s(t),t)+v_{s}(t)\mathbf{e}%
_{2,s}(s(t),t)\right) ,  \label{vpf}
\end{eqnarray}%
where we note the appearance of time and $s$ derivatives of $\mathbf{e}%
_{1}(s,t)$ and $\mathbf{e}_{2}(s,t)$ as well as the velocity $v_{s}$ with is
the component of $\mathbf{v}(\mathbf{x},t)$ tangential to $\mathbf{x}%
_{0}(s,t)$. Those derivatives will also appear in the various terms involved
in the vorticity equation (although we postpone the calculation till section %
\ref{secf}, see (\ref{ff1}), (\ref{ff2})). \ Concerning the derivatives of
the vectors $(\mathbf{e}_{1}(s,t),\mathbf{e}_{2}(s,t))$, we can prove the
following Lemma:

\begin{lemma}
\label{lemmae}Given a vortex filament $\mathbf{x}_{0}(s,t)$ $\in
C((0.T),C^{2}(T^{1},S^{2}))$, the unit vectors $(\mathbf{e}_{1}(s,t),\mathbf{%
e}_{2}(s,t))$ satisfy%
\begin{eqnarray*}
\left\vert \mathbf{e}_{1,s}(s,t)\right\vert +\left\vert \mathbf{e}%
_{2,s}(s,t)\right\vert &\leq &C \\
\left\vert \mathbf{e}_{1,t}(s,t)\right\vert +\left\vert \mathbf{e}%
_{2,t}(s,t)\right\vert &\leq &C\Gamma \left\vert \log (\nu t)\right\vert
\end{eqnarray*}
\end{lemma}

\textbf{Proof.} The estimate for $\mathbf{e}_{i,s}(s,t)$ simply follows from
(\ref{e11}), (\ref{e12}) with $\mathbf{t}_{s}=\kappa \mathbf{n}$ and
boundedness of the curvature. The estimate of $\mathbf{e}_{1,t}(s,t)$
follows from inverting the relations (\ref{fo1}), (\ref{fo2}) to obtain%
\begin{eqnarray*}
\mathbf{e}_{1}(s,t) &=&\cos \left( \theta _{0}(s,t)\right) \mathbf{n}%
(s,t)-\sin \left( \theta _{0}(s,t)\right) \mathbf{b}(s,t), \\
\mathbf{e}_{2}(s,t) &=&\sin \left( \theta _{0}(s,t)\right) \mathbf{n}%
(s,t)+\cos \left( \theta _{0}(s,t)\right) \mathbf{b}(s,t),
\end{eqnarray*}%
and hence, for $\mathbf{e}_{1,t}(s,t)$,
\begin{eqnarray*}
\mathbf{e}_{1,t}(s,t) &=&\left( -\sin \left( \theta _{0}(s,t)\right) \mathbf{%
n}(s,t)-\cos \left( \theta _{0}(s,t)\right) \mathbf{b}(s,t)\right) \theta
_{0,t}(s,t) \\
&&+\cos \left( \theta _{0}(s,t)\right) \mathbf{n}_{t}(s,t)-\sin \left(
\theta _{0}(s,t)\right) \mathbf{b}_{t}(s,t) \\
&=&\left( -\sin \left( \theta _{0}(s,t)\right) \mathbf{n}(s,t)-\cos \left(
\theta _{0}(s,t)\right) \mathbf{b}(s,t)\right) \int_{s_{0}}^{s}\tau
_{t}(s^{\prime },t)ds^{\prime } \\
&&+\cos \left( \theta _{0}(s,t)\right) \boldsymbol{\varpi }\times \mathbf{n}%
-\sin \left( \theta _{0}(s,t)\right) \boldsymbol{\varpi }\times \mathbf{b} \\
&=&\left( -\sin \left( \theta _{0}(s,t)\right) \mathbf{n}(s,t)-\cos \left(
\theta _{0}(s,t)\right) \mathbf{b}(s,t)\right) \int_{s_{0}}^{s}\tau
_{t}(s^{\prime },t)ds^{\prime } \\
&&+\cos \left( \theta _{0}(s,t)\right) \left( \int \tau _{t^{\prime }}ds-%
\frac{1}{2}\kappa ^{2}\right) \mathbf{b}+\cos \left( \theta _{0}(s,t)\right)
\frac{\Gamma }{2\pi }\log (\nu t)\kappa \tau \mathbf{t} \\
&&+\sin \left( \theta _{0}(s,t)\right) \left( \int \tau _{t^{\prime }}ds-%
\frac{1}{2}\kappa ^{2}\right) \mathbf{n}+\sin \left( \theta _{0}(s,t)\right)
\frac{\Gamma }{2\pi }\log (\nu t)\kappa _{s}\mathbf{t} \\
&=&-\frac{1}{2}\kappa ^{2}\mathbf{e}_{2}(s,t)+\frac{\Gamma }{2\pi }\log (\nu
t)\left( \cos \left( \theta _{0}(s,t)\right) \kappa \tau +\sin \left( \theta
_{0}(s,t)\right) \kappa _{s}\right) \mathbf{t},
\end{eqnarray*}%
and using (\ref{d0})-(\ref{d3}) and boundedness of curvature and torsion and
their derivatives we conclude:%
\begin{equation*}
\left\vert \mathbf{e}_{1,t}(s,t)\right\vert \leq C\left\vert \log (\nu
t)\right\vert .
\end{equation*}%
The vector $\mathbf{e}_{2,t}(s,t)$ can be estimated similarly and this
concludes the proof of the lemma.

We describe next the general strategy for our analysis in the next sections.
We write%
\begin{eqnarray*}
\boldsymbol{\omega } &=&\boldsymbol{\omega }_{0}+\widetilde{\boldsymbol{%
\omega }}, \\
\mathbf{v} &=&\mathbf{v}_{0}+\widetilde{\mathbf{v}},
\end{eqnarray*}%
where $\boldsymbol{\omega }_{0}$ is an approximate solution to the vorticity
equation in the moving frame and $\mathbf{v}_{0}$ the velocity field
resulting from Biot-Savart integral applied to $\boldsymbol{\omega }_{0}$.
Hence
\begin{equation}
\frac{\partial \boldsymbol{\omega }_{0}}{\partial t}+\boldsymbol{v}_{0}\cdot
\nabla _{\mathbf{x}}\boldsymbol{\omega }_{0}-\boldsymbol{\omega }_{0}\cdot
\nabla _{\mathbf{x}}\mathbf{v}_{0}-\nu \Delta _{\mathbf{x}}\boldsymbol{%
\omega }_{0}=-\mathbf{F}(\mathbf{x},t),  \notag
\end{equation}%
where $\mathbf{F}(\mathbf{x},t)$ will be small (in a norm to be defined) if $%
\boldsymbol{\omega }_{0}$ is sufficiently close to a solution of the
vorticity equation in the moving frame. The correction term $\widetilde{%
\mathbf{\omega }}$ satisfies then the equation
\begin{equation}
\frac{\partial \widetilde{\boldsymbol{\omega }}}{\partial t}-L_{0}\widetilde{%
\boldsymbol{\omega }}=\mathbf{F}(\mathbf{x},t)+\widetilde{\boldsymbol{\omega
}}\cdot \nabla _{\mathbf{x}}\widetilde{\mathbf{v}}-\widetilde{\boldsymbol{v}}%
\cdot \nabla _{\mathbf{x}}\widetilde{\boldsymbol{\omega }},  \label{eqnomega}
\end{equation}%
where $L_{0}$ is a linear operator. We will make a suitable choice of $%
\boldsymbol{\omega }_{0}$ so that $\mathbf{F}(\mathbf{x},t)$ is sufficiently
small as to guarantee that $\widetilde{\boldsymbol{\omega }}$ solution to (%
\ref{eqnomega}) is indeed a small perturbation to $\boldsymbol{\omega }_{0}$
itself. The vorticity $\boldsymbol{\omega }_{0}$ will be directed along \
the tangent to $\mathbf{x}_{0}(s,t)$, i.e. $\boldsymbol{\omega }_{0}=\omega
_{0}\boldsymbol{t}$, with $\omega _{0}$ being such that%
\begin{eqnarray}
&&A\left[ \omega _{0},v_{0r},v_{0\theta }\right]  \notag \\
&\equiv &\frac{\partial \omega _{0}}{\partial t}-\nu \left( \frac{1}{r^{2}}%
\frac{\partial ^{2}\omega _{0}}{\partial \theta ^{2}}+\frac{1}{r}\frac{%
\partial }{\partial r}\left( r\frac{\partial \omega _{0}}{\partial r}\right)
-\kappa (s,t)\cos \theta \frac{\partial \omega _{0}}{\partial r}\right)
\notag \\
&&+v_{0r}\frac{\partial \omega _{0}}{\partial r}+\frac{v_{0\theta }}{r}\frac{%
\partial \omega _{0}}{\partial \theta }-\kappa (s,t)\sin \theta \omega
_{0}v_{0\theta },  \label{defa}
\end{eqnarray}%
vanishes at leading order of $t$. In the definition of $A\left[ \omega
_{0},v_{0r},v_{0\theta }\right] $ in (\ref{defa}), $s$ may be viewed as a
parameter and $\omega _{0}$ is a function restricted to the normal $(\mathbf{%
n},\mathbf{b})$ plane that depends on polar coordinates $(r,\theta )$ in
that plane and only depends parametrically on $s$. The velocity field $%
(v_{0r},v_{0\theta })$ in $A\left[ \omega _{0},v_{0r},v_{0\theta }\right] $
results from applying Biot-Savart law to $\boldsymbol{\omega }_{0}$ and
keeping leading order (in $t$) terms. In this sense, $A\left[ \omega
_{0},v_{0r},v_{0\theta }\right] $ may be viewed as a $2D$ Navier-Stokes
operator defined on normal planes that includes information on the geometry
of $\mathbf{x}_{0}(s,t)$ through the curvature $\kappa (s,t)$. In the
construction of an approximate solution of $A=0$ in the next two sections we
will use Frenet-Serret local frames, while parallel frame coordinates will
fundamental to show that such approximate solution to $A\left[ \omega
_{0},v_{0r},v_{0\theta }\right] =0$ is also an approximate solution to the
vorticity formulation to Navier-Stokes equations.

\section{Lamb-Oseen vortex filament}

In this and the next section we will construct approximate solutions to
Navier-Stokes system. If we place ourselves in the local Frenet-Serret
system centered at a point $\mathbf{x}_{0}(s_{0},t)$ along the filament and
write cylindrical coordinates so that $\mathbf{x}^{\prime }=r\cos \theta
\mathbf{n}(s_{0},t)+r\sin \theta \mathbf{b}(s_{0},t)+z\mathbf{t}(s_{0},t)$
with $\mathbf{t}(s_{0},t)=\mathbf{e}_{z}$, the vorticity equation in that
frame is

\begin{eqnarray*}
&&\frac{\partial \omega _{z}}{\partial t}+v_{r}\frac{\partial \omega _{z}}{%
\partial r}+\frac{v_{\theta }}{r}\frac{\partial \omega _{z}}{\partial \theta
}+v_{z}\frac{\partial \omega _{z}}{\partial z} \\
&=&\omega _{r}\frac{\partial v_{z}}{\partial r}+\frac{\omega _{\theta }}{r}%
\frac{\partial v_{z}}{\partial \theta }+\omega _{z}\frac{\partial v_{z}}{%
\partial z}+\nu \Delta \omega _{z} \\
&&\frac{\partial \omega _{r}}{\partial t}+v_{r}\frac{\partial \omega _{r}}{%
\partial r}+\frac{v_{\theta }}{r}\frac{\partial \omega _{r}}{\partial \theta
}+v_{z}\frac{\partial \omega _{r}}{\partial z} \\
&=&\omega _{r}\frac{\partial v_{r}}{\partial r}+\frac{\omega _{\theta }}{r}%
\frac{\partial v_{r}}{\partial \theta }+\omega _{z}\frac{\partial v_{r}}{%
\partial z}+\nu \left( \Delta \omega _{r}-\frac{\omega _{r}}{r^{2}}-\frac{2}{%
r^{2}}\frac{\partial \omega _{\theta }}{\partial \theta }\right) \\
&&\frac{\partial \omega _{\theta }}{\partial t}+v_{r}\frac{\partial \omega
_{\theta }}{\partial r}+\frac{v_{\theta }}{r}\frac{\partial \omega _{\theta }%
}{\partial \theta }+v_{z}\frac{\partial \omega _{\theta }}{\partial z}-v_{r}%
\frac{\omega _{\theta }}{r} \\
&=&\omega _{r}\frac{\partial v_{\theta }}{\partial r}+\frac{\omega _{\theta }%
}{r}\frac{\partial v_{\theta }}{\partial \theta }+\omega _{z}\frac{\partial
v_{\theta }}{\partial z}-v_{\theta }\frac{\omega _{\theta }}{r}+\nu \left(
\Delta \omega _{\theta }-\frac{\omega _{\theta }}{r^{2}}+\frac{2}{r^{2}}%
\frac{\partial \omega _{r}}{\partial \theta }\right) ,
\end{eqnarray*}%
where%
\begin{equation*}
\Delta =\frac{1}{r}\frac{\partial }{\partial r}\left( r\frac{\partial }{%
\partial r}\right) +\frac{1}{r^{2}}\frac{\partial ^{2}}{\partial \theta ^{2}}%
+\frac{\partial ^{2}}{\partial z^{2}}.
\end{equation*}%
A $z-$independent solution is, of course, the Lamb-Oseen vortex $(\omega
_{r},\omega _{\theta },\omega _{z})=(0,0,\omega _{z}^{0})$ with:%
\begin{eqnarray}
\omega _{z}^{0} &=&\frac{1}{(\nu t)}\Omega _{0}\left( \frac{r}{(\nu t)^{%
\frac{1}{2}}}\right) ,  \label{grv} \\
v_{\theta }^{0} &=&\frac{1}{(\nu t)^{\frac{1}{2}}}V_{0}\left( \frac{r}{(\nu
t)^{\frac{1}{2}}}\right) ,  \label{gro}
\end{eqnarray}%
where the radial and angular unit vectors are defined by
\begin{eqnarray*}
\mathbf{e}_{\rho } &=&\left( \sin \theta \right) \mathbf{b}+\left( \cos
\theta \right) \mathbf{n}, \\
\mathbf{e}_{\theta } &=&\left( \cos \theta \right) \mathbf{b}-\left( \sin
\theta \right) \mathbf{n},
\end{eqnarray*}%
and%
\begin{eqnarray*}
\Omega _{0}\left( \rho \right) &=&\frac{\Gamma }{4\pi }e^{-\frac{\rho ^{2}}{4%
}}, \\
V_{0}\left( \rho \right) &=&\frac{\Gamma }{2\pi }\frac{1}{\rho }\left( 1-e^{-%
\frac{\rho ^{2}}{4}}\right) .
\end{eqnarray*}

For the sake of clarity, we will show next the expression given for the
velocity $V_{0}\left( \rho \right) $ together with the velocity associated
to the vorticity $(\omega _{r},\omega _{\theta },\omega _{z})=(0,0,\omega
_{z}^{0})$ with
\begin{equation*}
\omega _{z}^{0}=\frac{1}{(\nu t)^{\frac{1}{2}}}\Omega _{0}\left( \frac{r}{%
(\nu t)^{\frac{1}{2}}},\theta \right) ,
\end{equation*}%
where $\Omega _{0}\left( \rho ,\theta \right) =\frac{\Gamma }{4\pi }\rho e^{-%
\frac{\rho ^{2}}{4}}\cos \theta $ that will be used in this article.

\begin{lemma}
\label{lemma1}Let%
\begin{equation*}
\boldsymbol{\omega }^{0}=\frac{1}{(\nu t)}\Omega _{0}\left( \rho \right)
\mathbf{e}_{z},
\end{equation*}%
with%
\begin{equation}
\Omega _{0}\left( \rho \right) =\frac{\Gamma }{4\pi }e^{-\frac{\rho ^{2}}{4}%
}.  \label{re1}
\end{equation}%
Then, the associated velocity field that is bounded at the origin is given by%
\begin{equation*}
\mathbf{v}^{0}=\frac{1}{(\nu t)^{\frac{1}{2}}}V_{\theta }\left( \rho \right)
\mathbf{e}_{\theta },
\end{equation*}%
with%
\begin{equation}
V_{\theta }\left( \rho \right) =\frac{\Gamma }{2\pi }\frac{1}{\rho }\left(
1-e^{-\frac{1}{4}\rho ^{2}}\right) .  \label{re2}
\end{equation}%
Let%
\begin{equation*}
\boldsymbol{\omega }^{0}=\frac{1}{(\nu t)^{\frac{1}{2}}}\Omega _{0}\left(
\rho ,\theta \right) \mathbf{e}_{z},
\end{equation*}%
with
\begin{equation}
\Omega _{0}\left( \rho ,\theta \right) =\frac{\Gamma }{4\pi }\rho e^{-\frac{%
\rho ^{2}}{4}}\cos \theta .  \label{re3}
\end{equation}%
Then, the associated velocity field that is bounded at the origin is given by%
\begin{equation*}
\mathbf{v}^{0}=V_{r}\left( \rho ,\theta \right) \mathbf{e}_{r}+V_{\theta
}\left( \rho ,\theta \right) \mathbf{e}_{\theta },
\end{equation*}%
with%
\begin{eqnarray}
V_{r}\left( \rho ,\theta \right) &=&\frac{\Gamma }{4\pi }\left( \frac{1}{%
2\rho ^{2}}\left( 8e^{-\frac{1}{4}\rho ^{2}}+2\rho ^{2}e^{-\frac{1}{4}\rho
^{2}}-8\right) -e^{-\frac{1}{4}\rho ^{2}}\right) \sin \theta ,\   \label{re4}
\\
V_{\theta }\left( \rho ,\theta \right) &=&-\frac{\Gamma }{4\pi }\frac{1}{%
\rho ^{2}}\left( 4e^{-\frac{1}{4}\rho ^{2}}+2\rho ^{2}e^{-\frac{1}{4}\rho
^{2}}-4\right) \cos \theta ,  \label{re5}
\end{eqnarray}%
and, as $\rho \rightarrow 0$, one has a nonzero velocity along the binormal
direction given by%
\begin{equation}
V_{b}=-\frac{\Gamma }{4\pi },\ \ \text{as }\rho \rightarrow 0.
\label{vbzero}
\end{equation}
\end{lemma}

\bigskip \textbf{Proof.} \ In order to find the velocity field associated to
the vorticity $\Omega _{0}\left( \rho \right) $ we use the definition of
vorticity and the condition of incompressibility:
\begin{eqnarray*}
\frac{1}{\rho }\frac{\partial }{\partial \rho }\left( \rho V_{\theta }\left(
\rho ,\theta \right) \right) -\frac{1}{\rho }\frac{\partial V_{r}\left( \rho
,\theta \right) }{\partial \theta } &=&\Omega _{0}\left( \rho ,\theta
\right) , \\
\frac{1}{\rho }\frac{\partial }{\partial \rho }\left( \rho V_{r}\right) +%
\frac{1}{\rho }\frac{\partial V_{\theta }}{\partial \theta } &=&0.
\end{eqnarray*}%
Hence we can write, in terms of the stream function $\psi $,
\begin{equation*}
V_{r}=-\frac{1}{\rho }\frac{\partial \psi }{\partial \theta },\ V_{\theta }=%
\frac{\partial \psi }{\partial \rho },
\end{equation*}%
so that%
\begin{equation*}
\Delta \psi =\Omega _{0}\left( \rho ,\theta \right) ,
\end{equation*}%
and then, if $\Omega _{0}\left( \rho \right) =\frac{\Gamma }{4\pi }e^{-\frac{%
\rho ^{2}}{4}}$,%
\begin{equation*}
\psi =\frac{\Gamma }{4\pi }g(\rho ),
\end{equation*}%
where $g(\rho )$ satisfies \
\begin{equation*}
\frac{1}{\rho }\frac{d}{d\rho }\left( \rho \frac{dg(\rho )}{d\rho }\right)
=e^{-\frac{\rho ^{2}}{4}},
\end{equation*}%
$\allowbreak $and then%
\begin{equation*}
\frac{dg(\rho )}{d\rho }=\frac{2}{\rho }\left( 1-e^{-\frac{1}{4}\rho
^{2}}\right) ,
\end{equation*}%
so that%
\begin{equation*}
V_{\theta }\left( \rho \right) =\frac{\Gamma }{2\pi \rho }\left( 1-e^{-\frac{%
1}{4}\rho ^{2}}\right) .
\end{equation*}

In the case that $\Omega _{0}\left( \rho \right) =\frac{\Gamma }{4\pi }\rho
e^{-\frac{\rho ^{2}}{4}}\cos \theta $,%
\begin{equation*}
\psi \left( \rho ,\theta \right) =\frac{\Gamma }{4\pi }g(\rho )\cos \theta ,
\end{equation*}%
with%
\begin{equation*}
\frac{1}{\rho }\frac{d}{d\rho }\left( \rho \frac{dg(\rho )}{d\rho }\right) -%
\frac{1}{\rho ^{2}}g(\rho )=\rho e^{-\frac{\rho ^{2}}{4}}.
\end{equation*}

Hence, solving the equation above%
\begin{eqnarray*}
g(\rho ) &=&-\rho \int_{\rho }^{\infty }\rho ^{\prime -3}\left(
\int_{0}^{\rho ^{\prime }}\rho ^{\prime \prime 3}e^{-\frac{\rho ``^{2}}{4}%
}d\rho ^{\prime \prime }\right) d\rho ^{\prime } \\
&=&\frac{1}{2\rho }\left( 8e^{-\frac{1}{4}\rho ^{2}}+2\rho ^{2}e^{-\frac{1}{4%
}\rho ^{2}}-8\right) -\rho e^{-\frac{1}{4}\rho ^{2}},
\end{eqnarray*}%
with%
\begin{equation*}
g(\rho )=-\rho +O(\rho ^{3})\text{ as }\rho \rightarrow 0,\ g(\rho )=O(\rho
^{-1})\text{ as }\rho \rightarrow \infty ,
\end{equation*}%
and we have
\begin{eqnarray*}
\psi &=&-\frac{\Gamma }{4\pi }\rho \cos \theta ,\text{as }\rho \rightarrow 0,
\\
V_{r} &=&-\frac{\Gamma }{4\pi }\sin \theta ,\ V_{\theta }=-\frac{\Gamma }{%
4\pi }\cos \theta ,\ \ \ \ \text{as }\rho \rightarrow 0.
\end{eqnarray*}%
Notice that this implies a nonzero contribution to the binormal velocity at
the vortex filament:%
\begin{equation*}
V_{b}=-\frac{\Gamma }{4\pi },\ \ \text{as }\rho \rightarrow 0.
\end{equation*}%
More generally, for any $\rho $,%
\begin{eqnarray*}
V_{r} &=&\frac{\Gamma }{4\pi }\left( \frac{1}{2\rho ^{2}}\left( 8e^{-\frac{1%
}{4}\rho ^{2}}+2\rho ^{2}e^{-\frac{1}{4}\rho ^{2}}-8\right) -e^{-\frac{1}{4}%
\rho ^{2}}\right) \sin \theta ,\  \\
V_{\theta } &=&-\frac{\Gamma }{4\pi }\frac{1}{\rho ^{2}}\left( 4e^{-\frac{1}{%
4}\rho ^{2}}+2\rho ^{2}e^{-\frac{1}{4}\rho ^{2}}-4\right) \cos \theta ,
\end{eqnarray*}%
and this concludes the proof of the Lemma.

Hence, we introduce an approximate solution consisting of a Lamb-Oseen
vortex along the filament $\mathbf{x}_{0}(s,t)$ moving with the binormal
flow:%
\begin{equation}
\boldsymbol{\omega }_{LO}=\frac{1}{(\nu t)}\Omega _{0}\left( \rho \right)
\mathbf{t},  \label{lo1}
\end{equation}%
where%
\begin{equation}
\rho =\frac{dist(\mathbf{x},\mathbf{x}_{0})}{(\nu t)^{\frac{1}{2}}},\ \
\label{rho}
\end{equation}%
The distance of a given point $\mathbf{x}\in
\mathbb{R}
^{3}$ to the filament is uniquely well defined in a sufficiently small
neighborhood of the filament. We introduce then the vortex tube%
\begin{equation}
\boldsymbol{\omega }^{0}(\mathbf{x},t)=\frac{1}{(\nu t)}\Omega _{0}\left(
\rho \right) \eta _{R}(dist(\mathbf{x},\mathbf{x}_{0}))\mathbf{t},
\label{lo1p}
\end{equation}%
with $\eta _{R}$ a cutoff function defined by (\ref{co0})-(\ref{co2}), and
prove the following Lemma concerning the induced velocity:

\begin{lemma}
Let $\boldsymbol{\omega }^{0}(\mathbf{x},t)$ given by (\ref{lo1p}) be a
vorticity field with $\rho $ given by (\ref{rho}). Assume in addition that
the filament $\mathbf{x}_{0}(s,t)$ is three times differentiable in $s$
uniformly in $t$. Then, the induced velocity field given by Biot-Savart law
is, in the tube of radius $R$ around the filament, given by%
\begin{eqnarray}
\mathbf{v}^{0}(\rho ,\theta ,s,t) &=&\frac{1}{(\nu t)^{\frac{1}{2}}}%
V_{0}(\rho )\mathbf{e}_{\theta }(s,t)  \notag \\
&&-\frac{\Gamma }{4\pi }\kappa (s,t)\left( \log (\nu t)^{\frac{1}{2}}+1-\log
2+\frac{1}{4\pi }(F(\rho )+H(\rho ))\right) \mathbf{b}(s,t)  \notag \\
&&-\kappa (s,t)\mathbf{V}(\rho ,\theta )+\mathbf{v}^{\ast }(s,t)+O((\nu t)^{%
\frac{1}{2}}),  \label{v0}
\end{eqnarray}%
with $F(\rho )$ a smooth function such that $F(s)\sim 4\pi \log s,\ $as $%
s\rightarrow \infty $ and $F(0)=2\pi (2\log 2-\gamma )$ ($\gamma $ being
Euler's constant),$\ H(\rho )$ a smooth function such that $H(0)=2\pi $ and $%
O(\rho ^{-2})$ as $\rho \rightarrow \infty $,
\begin{equation*}
\mathbf{V}(\rho ,\theta )=V_{r}\left( \rho ,\theta \right) \mathbf{e}%
_{r}+V_{\theta }\left( \rho ,\theta \right) \mathbf{e}_{\theta },
\end{equation*}%
with $V_{r}\left( \rho ,\theta \right) $, $V_{\theta }\left( \rho ,\theta
\right) $ given by (\ref{re4}), (\ref{re5}) respectively, and%
\begin{equation}
\mathbf{v}^{\ast }(s,t)=\frac{\Gamma }{4\pi }\lim_{\varepsilon \rightarrow
0}\left( \int_{\left\vert s^{\prime }-s\right\vert >\varepsilon }\mathbf{t}%
(s^{\prime },t)\times \frac{\mathbf{x}(s,t)-\mathbf{x}^{\prime }(s^{\prime
},t)}{\left\vert \mathbf{x}(s,t)-\mathbf{x}^{\prime }(s^{\prime
},t)\right\vert ^{3}}ds^{\prime }+\kappa (s,t)\mathbf{b}(s,t)\log
\varepsilon \right) .  \label{vb}
\end{equation}%
At the center of the filament we have%
\begin{eqnarray}
\mathbf{v}^{0}(0,s,t) &=&\frac{\Gamma }{4\pi }\kappa (s,t)\left( \log (\nu
t)^{-\frac{1}{2}}+\log 2+\frac{1}{2}(\gamma -2\log 2)-\frac{1}{2}\right)
\mathbf{b}(s,t)  \notag \\
&&+\mathbf{v}^{\ast }(s,t)+O((\nu t)^{\frac{1}{2}}).  \label{va}
\end{eqnarray}
\end{lemma}

\textbf{Proof.} In order to compute the velocity field created by the
vorticity $\boldsymbol{\omega }^{0}$ we write%
\begin{eqnarray*}
\boldsymbol{v}^{0}(\mathbf{x},t) &=&\frac{1}{4\pi }\int \boldsymbol{\omega }%
^{0}(\mathbf{x}^{\prime },t)\times \frac{\mathbf{x}-\mathbf{x}^{\prime }}{%
\left\vert \mathbf{x}-\mathbf{x}^{\prime }\right\vert ^{3}}d\mathbf{x}%
^{\prime } \\
&=&\frac{1}{4\pi }\int_{\left\vert s^{\prime }-s_{0}\right\vert >\varepsilon
}+\frac{1}{4\pi }\int_{\left\vert s^{\prime }-s_{0}\right\vert \leq
\varepsilon }\equiv I_{1}+I_{2}.
\end{eqnarray*}%
First we note, since%
\begin{equation*}
\int_{0}^{2\pi }\int_{0}^{\infty }\frac{1}{(\nu t)}\Omega _{0}\left( \rho
\right) rdrd\theta =\Gamma ,
\end{equation*}%
that%
\begin{equation*}
I_{1}=\frac{\Gamma }{4\pi }\int_{\left\vert s^{\prime }-s_{0}\right\vert
>\varepsilon }\mathbf{t}(s^{\prime },t)\times \frac{\mathbf{x}-\mathbf{x}%
^{\prime }(s^{\prime },t)}{\left\vert \mathbf{x}-\mathbf{x}^{\prime
}(s^{\prime },t)\right\vert ^{3}}ds^{\prime }+O((\nu t)^{\frac{1}{2}}),
\end{equation*}%
which, according to the results in the appendix,
\begin{eqnarray}
&&\frac{\Gamma }{4\pi }\int_{\left\vert s^{\prime }-s_{0}\right\vert
>\varepsilon }\mathbf{t}(s^{\prime },t)\times \frac{\mathbf{x}-\mathbf{x}%
^{\prime }(s^{\prime },t)}{\left\vert \mathbf{x}-\mathbf{x}^{\prime
}(s^{\prime },t)\right\vert ^{3}}ds^{\prime }  \notag \\
&=&\frac{\Gamma }{4\pi }\left( \kappa (s_{0},t)\mathbf{n}(s_{0},t)\times
\mathbf{t}(s_{0},t)\right) \log \varepsilon +O(1)  \notag \\
&=&-\frac{\Gamma }{4\pi }\kappa (s_{0},t)\mathbf{b}(s_{0},t)\log \varepsilon
+O(1).  \label{I1}
\end{eqnarray}%
Hence we can write%
\begin{eqnarray}
I_{1} &=&-\frac{\Gamma }{4\pi }\kappa (s_{0},t)\mathbf{b}(s_{0},t)\log
\varepsilon  \notag \\
&&+\left[ \frac{\Gamma }{4\pi }\int_{\left\vert s^{\prime }-s_{0}\right\vert
>\varepsilon }\frac{\mathbf{t}(s^{\prime },t)\times \left( \mathbf{x}-%
\mathbf{x}^{\prime }(s^{\prime },t)\right) }{\left\vert \mathbf{x}-\mathbf{x}%
^{\prime }(s^{\prime },t)\right\vert ^{3}}ds^{\prime }+\frac{\Gamma }{4\pi }%
\kappa (s_{0},t)\mathbf{b}(s_{0},t)\log \varepsilon \right]  \notag \\
&&+O((\nu t)^{\frac{1}{2}}),  \label{i1c}
\end{eqnarray}%
where the first term at the right hand side of (\ref{i1c}) will cancel out
with a similar term coming from $I_{2}$ and the second term a the right hand
side of (\ref{i1c}) will provide, in the limit $\varepsilon \rightarrow 0$,
with the term (\ref{vb}).

We estimate next $I_{2}$. Locally near $\mathbf{x}_{0}(s_{0})$ (we drop the $%
t$ dependence for the sake of simplicity), we can estimate the Biot-Savart
kernel in the following fashion:%
\begin{equation*}
\frac{\mathbf{x}-\mathbf{x}^{\prime }}{\left\vert \mathbf{x}-\mathbf{x}%
^{\prime }\right\vert ^{3}}=\frac{r\mathbf{e}_{r}-r^{\prime }\mathbf{e}%
_{r^{\prime }}-\mathbf{r}_{0,s}(s_{0})(s^{\prime }-s_{0})-\frac{1}{2}\kappa
\mathbf{n}(s_{0})(s_{0}-s^{\prime })^{2}+O((s_{0}-s^{\prime })^{3})}{\left(
\left\vert r\mathbf{e}_{r}-r^{\prime }\mathbf{e}_{r^{\prime }}\right\vert
^{2}+(s^{\prime }-s_{0})^{2}-\kappa (s_{0}-s^{\prime })^{2}(\rho \mathbf{e}%
_{r}-\rho ^{\prime }\mathbf{e}_{r^{\prime }})\cdot \mathbf{n}(s_{0})\right)
^{\frac{3}{2}}},
\end{equation*}%
so that%
\begin{eqnarray}
&&\boldsymbol{\omega }^{0}(\mathbf{x}^{\prime },t)\times \frac{\mathbf{x}-%
\mathbf{x}^{\prime }}{\left\vert \mathbf{x}-\mathbf{x}^{\prime }\right\vert
^{3}}  \notag \\
&=&\frac{\Gamma e^{-\frac{r^{\prime 2}}{4\nu t}}}{4\pi \nu t}\left( \mathbf{t%
}_{0}(s_{0},t)+\kappa \mathbf{n}_{0}(s_{0},t)(s^{\prime }-s_{0})+\ldots
\right)  \notag \\
&&\times \frac{r\mathbf{e}_{r}-r^{\prime }\mathbf{e}_{r^{\prime }}-\mathbf{r}%
_{0,s}(s_{0})(s^{\prime }-s_{0})-\frac{1}{2}\kappa \mathbf{n}%
(s_{0})(s_{0}-s^{\prime })^{2}+...}{\left( \left\vert r\mathbf{e}%
_{r}-r^{\prime }\mathbf{e}_{r^{\prime }}\right\vert ^{2}+(s^{\prime
}-s_{0})^{2}-\kappa (s_{0}-s^{\prime })^{2}(r\mathbf{e}_{r}-r^{\prime }%
\mathbf{e}_{r^{\prime }})\cdot \mathbf{n}(s_{0})+\ldots \right) ^{\frac{3}{2}%
}}  \notag \\
&=&\frac{\Gamma e^{-\frac{r^{\prime 2}}{4\nu t}}\left[ \mathbf{t}%
_{0}(s_{0},t)\times \left( r\mathbf{e}_{r}-r^{\prime }\mathbf{e}_{r^{\prime
}}\right) +\kappa \mathbf{n}_{0}(s_{0},t)(s^{\prime }-s_{0})\times \left( r%
\mathbf{e}_{r}-r^{\prime }\mathbf{e}_{r^{\prime }}\right) +\frac{1}{2}\kappa
\mathbf{b}(s^{\prime }-s_{0})^{2}\right] }{4\pi \nu t\left( \left\vert r%
\mathbf{e}_{r}-r^{\prime }\mathbf{e}_{r^{\prime }}\right\vert
^{2}+(s^{\prime }-s_{0})^{2}\right) ^{\frac{3}{2}}}  \notag \\
&&+\frac{3}{2}\frac{\Gamma e^{-\frac{r^{\prime 2}}{4\nu t}}\left[ \mathbf{t}%
_{0}(s_{0},t)\times \left( r\mathbf{e}_{r}-r^{\prime }\mathbf{e}_{r^{\prime
}}\right) \right] \kappa (s_{0}-s^{\prime })^{2}(r\mathbf{e}_{r}-r^{\prime }%
\mathbf{e}_{r^{\prime }})\cdot \mathbf{n}(s_{0})}{4\pi \nu t\left(
\left\vert r\mathbf{e}_{r}-r^{\prime }\mathbf{e}_{r^{\prime }}\right\vert
^{2}+(s^{\prime }-s_{0})^{2}\right) ^{\frac{5}{2}}}+\ldots ,  \label{rhs}
\end{eqnarray}%
and hence, using (\ref{scale1}), (\ref{scale2}),%
\begin{eqnarray*}
I_{2} &=&\frac{1}{4\pi }\int_{\left\vert s^{\prime }-s_{0}\right\vert \leq
\varepsilon }\boldsymbol{\omega }_{0}(\mathbf{x}^{\prime },t)\times \frac{%
\mathbf{x}-\mathbf{x}^{\prime }}{\left\vert \mathbf{x}-\mathbf{x}^{\prime
}\right\vert ^{3}}\rho ^{\prime }h_{s^{\prime }}d\rho ^{\prime }ds^{\prime
}d\theta ^{\prime } \\
&=&\frac{1}{(\nu t)^{\frac{1}{2}}}V_{0}(r(\nu t)^{-\frac{1}{2}})\mathbf{e}%
_{\theta }+\frac{1}{8\pi }\kappa I\mathbf{b}(s_{0},t)+\frac{1}{4\pi }\kappa
\mathbf{W}(r(\nu t)^{-\frac{1}{2}},\theta )\mathbf{-}\kappa \mathbf{V}(r(\nu
t)^{-\frac{1}{2}},\theta ) \\
&&+O((\nu t)^{\frac{1}{2}})
\end{eqnarray*}%
which is a superposition of angular rotation, binormal velocity, a $O(1)$
velocity $\mathbf{W}$ due to the last term at the right-hand side of (\ref%
{rhs}) and a $O(1)$ velocity $\mathbf{V}$ due to the correction to the
volume introduced by the scale factor $h_{s^{\prime }}$ (see (\ref{scale1}),
(\ref{scale2})). The function $V_{0}(\rho )$ is defined by (\ref{gro}). The
value of $I$ is%
\begin{eqnarray}
I &\equiv &\int_{\left\vert s^{\prime }-s_{0}\right\vert \leq \varepsilon }%
\frac{\Gamma e^{-\frac{r^{\prime 2}}{4\nu t}}(s^{\prime }-s_{0})^{2}}{4\pi
\nu t\left( \left\vert r\mathbf{e}_{r}-r^{\prime }\mathbf{e}_{r^{\prime
}}\right\vert ^{2}+(s^{\prime }-s_{0})^{2}\right) ^{\frac{3}{2}}}r^{\prime
}dr^{\prime }ds^{\prime }d\theta ^{\prime }  \notag \\
&=&\int_{0}^{\infty }\int_{0}^{2\pi }\int_{-\frac{_{\varepsilon }}{(\nu t)^{%
\frac{1}{2}}\left\vert r(\nu t)^{-\frac{1}{2}}\mathbf{e}_{r}-\rho ^{\prime }%
\mathbf{e}_{r^{\prime }}\right\vert }}^{\frac{_{\varepsilon }}{(\nu t)^{%
\frac{1}{2}}\left\vert r(\nu t)^{-\frac{1}{2}}\mathbf{e}_{r}-\rho ^{\prime }%
\mathbf{e}_{r^{\prime }}\right\vert }}\frac{\Gamma e^{-\frac{\rho ^{\prime 2}%
}{4}}\sigma ^{2}}{4\pi \left( 1+\sigma ^{2}\right) ^{\frac{3}{2}}}\rho
^{\prime }d\rho ^{\prime }d\sigma d\theta ^{\prime }  \notag
\end{eqnarray}%
\begin{eqnarray}
&=&-\frac{\Gamma }{2\pi }\int_{0}^{2\pi }\int_{0}^{\infty }e^{-\frac{\rho
^{\prime 2}}{4}}\left( \log \frac{(\nu t)^{\frac{1}{2}}}{\varepsilon }+\log
\left\vert \frac{r}{(\nu t)^{\frac{1}{2}}}\mathbf{e}_{r}-\rho ^{\prime }%
\mathbf{e}_{r^{\prime }}\right\vert -\log 2+1\right) \rho ^{\prime }d\rho
^{\prime }d\theta ^{\prime }  \notag \\
&&+O((\nu t)^{\frac{1}{2}})  \notag \\
&=&-2\Gamma \log \frac{(\nu t)^{\frac{1}{2}}}{\varepsilon }-\frac{\Gamma }{%
2\pi }\int_{0}^{2\pi }\int_{0}^{\infty }e^{-\frac{\rho ^{\prime 2}}{4}}\log
\left\vert r(\nu t)^{-\frac{1}{2}}\mathbf{e}_{r}-\rho ^{\prime }\mathbf{e}%
_{r^{\prime }}\right\vert \rho ^{\prime }d\rho ^{\prime }d\theta ^{\prime }
\notag \\
&&+2\Gamma (\log 2-1)+O((\nu t)^{\frac{1}{2}})  \notag \\
&=&-2\Gamma \log \frac{(\nu t)^{\frac{1}{2}}}{\varepsilon }-\frac{\Gamma }{%
2\pi }F(r(\nu t)^{-\frac{1}{2}})+2\Gamma (\log 2-1)+O((\nu t)^{\frac{1}{2}}),
\label{I}
\end{eqnarray}%
where%
\begin{equation*}
F(s)=\int_{0}^{2\pi }\int_{0}^{\infty }e^{-\frac{\rho ^{\prime 2}}{4}}\log
\left\vert s\mathbf{e}_{r}-\rho ^{\prime }\mathbf{e}_{r^{\prime
}}\right\vert \rho ^{\prime }d\rho ^{\prime }d\theta ^{\prime }
\end{equation*}%
with%
\begin{eqnarray}
F(s) &\sim &4\pi \log \left\vert s\right\vert ,\ \text{as }s\rightarrow
\infty  \notag \\
F(0) &=&2\pi (2\log 2-\gamma )  \label{F0}
\end{eqnarray}%
Concerning $\mathbf{W}(\rho ,\theta )$, note that it is directed in the
binormal direction (due to the scalar product $(r\mathbf{e}_{r}-r^{\prime }%
\mathbf{e}_{r^{\prime }})\cdot \mathbf{n}(s_{0})$ in the integral only the
normal component of $r\mathbf{e}_{r}-r^{\prime }\mathbf{e}_{r^{\prime }}$
counts, but this quantity appears in $\mathbf{t}_{0}(s_{0},t)\times \left( r%
\mathbf{e}_{r}-r^{\prime }\mathbf{e}_{r^{\prime }}\right) $ which is
therefore a vector along $\mathbf{b}(s_{0},t)$) \ and therefore one can
write $\mathbf{W}(\rho ,\theta )=W_{r}\left( \rho ,\theta \right) \mathbf{e}%
_{r}+W_{\theta }\left( \rho ,\theta \right) \mathbf{e}_{\theta }$ with $%
W_{r}\left( \rho ,\theta \right) =-\frac{\Gamma }{4\pi }H(\rho )\sin \theta $
where%
\begin{equation}
H(\rho )=-\frac{3}{2}\left( \int \frac{e^{-\frac{\rho ^{\prime 2}}{4}}\sigma
^{\prime 2}\left\vert \rho \mathbf{e}_{\rho }-\rho ^{\prime }\mathbf{e}%
_{\rho ^{\prime }}\right\vert ^{\mathbf{2}}\cos ^{2}\theta ^{\prime }}{%
\left( \left\vert \rho \mathbf{e}_{\rho }-\rho ^{\prime }\mathbf{e}_{\rho
^{\prime }}\right\vert ^{\mathbf{2}}+\sigma ^{\prime 2}\right) ^{\frac{5}{2}}%
}d\sigma ^{\prime }d\theta ^{\prime }d\rho ^{\prime }\right) ,  \label{hrho}
\end{equation}%
and easily verify $H(0)=-2\pi $ and $\left\vert H(\rho )\right\vert =O(\rho
^{-2})$ as $\rho \rightarrow \infty $ as well as%
\begin{eqnarray}
\mathbf{W}(0,\theta ) &=&\frac{3}{2}\left( \int \frac{\Gamma e^{-\frac{\rho
^{\prime 2}}{4}}\sigma ^{\prime 2}\rho ^{\prime \mathbf{2}}\cos ^{2}\theta
^{\prime }}{4\pi \left( \rho ^{\prime \mathbf{2}}+\sigma ^{\prime 2}\right)
^{\frac{5}{2}}}d\sigma ^{\prime }d\theta ^{\prime }d\rho ^{\prime }\right)
\mathbf{b}(s_{0},t)  \notag \\
&=&\frac{\Gamma }{2}\mathbf{b}(s_{0},t).  \label{w0}
\end{eqnarray}%
Finally, $\mathbf{V}(\rho ,\theta )$ is the velocity field due to a
two-dimensional vorticity field $\Omega _{0}\left( \rho ,\theta \right)
=\Gamma e^{-\frac{\rho ^{2}}{4}}\rho \cos \theta $, which is (by Lemma \ref%
{lemma1}) given by $\mathbf{V}(\rho ,\theta )=V_{r}\left( \rho ,\theta
\right) \mathbf{e}_{r}+V_{\theta }\left( \rho ,\theta \right) \mathbf{e}%
_{\theta }$ with $V_{r}$, $V_{\theta }$ given by (\ref{re4}), (\ref{re5})
respectively. Note that%
\begin{equation*}
V_{r}\left( \rho ,\theta \right) =\frac{\Gamma }{(4\pi )^{2}}G(\rho )\sin
\theta ,
\end{equation*}%
with%
\begin{equation}
G(\rho )=4\pi \left( \frac{1}{2\rho ^{2}}\left( 8e^{-\frac{1}{4}\rho
^{2}}+2\rho ^{2}e^{-\frac{1}{4}\rho ^{2}}-8\right) -e^{-\frac{1}{4}\rho
^{2}}\right) ,  \label{grho}
\end{equation}%
and hence $G(\rho )=O(\rho ^{-2})$ as $\rho \rightarrow \infty $ and%
\begin{equation}
\mathbf{V}(0,\theta )=\frac{\Gamma }{4\pi }\mathbf{b}(s_{0},t).  \label{V0}
\end{equation}%
We also notice that the $\log (\varepsilon )$ contribution in (\ref{I})
cancels with the $\log (\varepsilon )$ term in (\ref{I1}).

At $\rho =0$ we have, by putting together (\ref{I}) with (\ref{F0}), (\ref%
{w0}), (\ref{V0}) and the contribution from (\ref{I1}),%
\begin{eqnarray*}
\mathbf{v}^{0}(\rho &=&0) \\
&=&\frac{\Gamma }{4\pi }\kappa (s,t)\left( \log (\nu t)^{-\frac{1}{2}}+\log
2+\frac{1}{2}(\gamma -2\log 2)-\frac{1}{2}\right) \mathbf{b}(s,t) \\
&&+\frac{\Gamma }{4\pi }\lim_{\varepsilon \rightarrow 0}\left(
\int_{\left\vert s^{\prime }-s\right\vert >\varepsilon }\mathbf{t}(s^{\prime
},t)\times \frac{\mathbf{x}(s,t)-\mathbf{x}^{\prime }(s^{\prime },t)}{%
\left\vert \mathbf{x}(s,t)-\mathbf{x}^{\prime }(s^{\prime },t)\right\vert
^{3}}ds^{\prime }+\kappa (s,t)\mathbf{b}(s,t)\log \varepsilon \right) \\
&&+O((\nu t)^{\frac{1}{2}}),
\end{eqnarray*}%
which is, at leading order, the binormal flow velocity. We conclude then
that the velocity field created by the vorticity $\boldsymbol{\omega }_{0}$
around the filament is given by formula (\ref{v0}) and this concludes the
proof of the Lemma.

\begin{remark}
Formula (\ref{va}) provides a method to compute the initial velocity of a
vortex ring of radius $\mathit{R}$. One can easily compute
\begin{equation*}
\mathbf{v}^{\ast }(s,t)=\frac{\Gamma }{4\pi \mathit{R}}\log (4\mathit{R})%
\mathbf{b}(s,t),
\end{equation*}%
providing the formula%
\begin{equation*}
\mathbf{v}^{0}=\frac{\Gamma }{4\pi \mathit{R}}\left( \log \left( \frac{8%
\mathit{R}}{(\nu t)^{\frac{1}{2}}}\right) +\frac{1}{2}(\gamma -2\log 2)-%
\frac{1}{2}+O((\nu t)^{\frac{1}{2}}\log (\nu t))\right) \mathbf{b},
\end{equation*}%
which can also be computed from the explicit formulas for the velocity field
of a vortex ring provided by Lamb (cf. \S\ 161 in \cite{L}, see also \cite{F}%
).
\end{remark}

Next, we notice that a Lamb-Oseen vortex for a straight filament is such
that the nonlinear terms in Navier-Stokes system (the stretching $\mathbf{%
\omega }_{0}\cdot \nabla \mathbf{v}_{0}$ and the convection $\mathbf{v}%
_{0}\cdot \nabla \boldsymbol{\omega }_{0}$ terms) identically cancel . This,
of course, is not the case for (\ref{lo1}) but we will see that
cancellations take place at leading order in $t$ and this will be crucial in
the next section in order to compute lower order corrections to the
vorticity. Specifically, we can prove the following Lemma:

\begin{lemma}
\bigskip Let $\boldsymbol{\omega }^{0}$ be given by (\ref{lo1p}) and $%
\mathbf{v}^{0}$ given by (\ref{v0}). Then, inside the vortex tube of radius $%
R$,%
\begin{eqnarray*}
\boldsymbol{\omega }^{0}\cdot \nabla \mathbf{v}^{0} &=&-\frac{1}{(\nu t)^{%
\frac{3}{2}}}\kappa \Omega _{0}\left( \rho \right) V_{0}(\rho )\sin \theta
\mathbf{t}+O((\nu t)^{-1}) \\
\mathbf{v}^{0}\cdot \nabla \boldsymbol{\omega }^{0} &=&-\frac{1}{(\nu t)^{%
\frac{3}{2}}}\frac{\Gamma }{(4\pi )^{2}}\kappa (4\pi (1-\log 2)+F(\rho
)+G(\rho )+H(\rho ))\frac{d\Omega _{0}\left( \rho \right) }{d\rho }\sin
\theta \mathbf{t} \\
&&+\frac{1}{(\nu t)^{\frac{3}{2}}}(v_{0b}^{\ast }\sin \theta +v_{0n}^{\ast
}\cos \theta )\frac{d\Omega _{0}\left( \rho \right) }{d\rho }\mathbf{t}%
+O((\nu t)^{-1}) \\
\nu \frac{\partial ^{2}\boldsymbol{\omega }^{0}}{\partial z^{2}} &=&-\frac{%
\nu }{(\nu t)^{\frac{3}{2}}}\kappa \Omega _{0}^{\prime }\left( \rho \right)
\cos \theta \mathbf{t}+O((\nu t)^{-1}).
\end{eqnarray*}
\end{lemma}

\textbf{Proof. }The vortex stretching term is, at leading order inside the
vortex tube
\begin{equation*}
\boldsymbol{\omega }^{0}\cdot \nabla \mathbf{v}^{0}=\frac{1}{(\nu t)}\Omega
_{0}\left( \rho \right) \mathbf{t\cdot \nabla v}_{0}=\frac{1}{(\nu t)}\Omega
_{0}\left( \rho \right) \frac{d\mathbf{v}_{0}}{ds},
\end{equation*}%
and since
\begin{eqnarray*}
\frac{d\mathbf{v}^{0}}{ds} &=&\frac{1}{(\nu t)^{\frac{1}{2}}}V_{0}(\rho )%
\frac{d}{ds}\left( \cos (\varphi -\theta _{0}(s))\mathbf{b}-\sin (\varphi
-\theta _{0}(s))\mathbf{n}\right) \\
&=&\frac{1}{(\nu t)^{\frac{1}{2}}}V_{0}(\rho )\left( -\tau \cos (\varphi
-\theta _{0}(s))\mathbf{n}+\kappa \sin (\varphi -\theta _{0}(s))\mathbf{t}%
-\tau \sin (\varphi -\theta _{0}(s))\mathbf{b}\right) \\
&&+\frac{1}{(\nu t)^{\frac{1}{2}}}V_{0}(\rho )\left( \theta _{0}^{\prime
}(s)\sin (\varphi -\theta _{0}(s))\mathbf{b}+\theta _{0}^{\prime }(s)\cos
(\varphi -\theta _{0}(s))\mathbf{n}\right) \\
&=&-\frac{1}{(\nu t)^{\frac{1}{2}}}\kappa V_{0}(\rho )\sin \theta \mathbf{t},
\end{eqnarray*}%
with%
\begin{equation*}
\theta _{0}^{\prime }(s)=\tau (s),
\end{equation*}%
we have%
\begin{equation}
\boldsymbol{\omega }^{0}\cdot \nabla \mathbf{v}^{0}=-\frac{1}{(\nu t)^{\frac{%
3}{2}}}\kappa \Omega _{0}\left( \rho \right) V_{0}(\rho )\sin \theta \mathbf{%
t},  \label{c1b}
\end{equation}%
yielding the important conclusion that (\ref{lo1}) produces a $O((\nu t)^{-%
\frac{3}{2}})$ vortex stretching along the $\mathbf{t}$ direction.

The fact that $\boldsymbol{\omega }^{0}$ depends solely on $\rho $ implies
that for any $O(1)$ regular vector field $\mathbf{v}^{\ast }(s,t)=v_{b}(s,t)%
\mathbf{b}(s,t)+v_{n}(s,t)\mathbf{n}(s,t)$ one also has%
\begin{equation}
\left( \mathbf{v}^{\ast }\cdot \nabla \right) \boldsymbol{\omega }^{0}=\frac{%
1}{(\nu t)^{\frac{3}{2}}}(v_{b}^{\ast }\sin \theta +v_{n}^{\ast }\cos \theta
)\frac{d\Omega _{0}\left( \rho \right) }{d\rho }\mathbf{t}+O((\nu t)^{-1}).
\label{c3}
\end{equation}%
This in particular applies to the velocity field (\ref{vb}) created by the
filament beyond the leading binormal flow term. As for the stretching $%
\boldsymbol{\omega }^{0}\cdot \nabla \mathbf{v}^{\ast }$ produced by $%
\mathbf{v}^{\ast }$, it will be $O((\nu t)^{-1})$ due to the fact that $%
\mathbf{\omega }^{0}$ is $O((\nu t)^{-1})$. Therefore, the contribution of (%
\ref{lo1}) to the convective nonlinearity is%
\begin{eqnarray}
\mathbf{v}^{0}\cdot \nabla \boldsymbol{\omega }^{0} &=&-\frac{1}{(\nu t)^{%
\frac{3}{2}}}\frac{\Gamma }{(4\pi )^{2}}\kappa (4\pi (1-\log 2)+F(\rho
)+G(\rho )+H(\rho ))\frac{d\Omega _{0}\left( \rho \right) }{d\rho }\sin
\theta \mathbf{t}  \notag \\
&&+\frac{1}{(\nu t)^{\frac{3}{2}}}(v_{0b}^{\ast }\sin \theta +v_{0n}^{\ast
}\cos \theta )\frac{d\Omega _{0}\left( \rho \right) }{d\rho }\mathbf{t+}%
O((\nu t)^{-1}),  \label{c2b}
\end{eqnarray}%
again a $O((\nu t)^{-\frac{3}{2}})$ contribution along the $\mathbf{t}$
direction. Finally, the curvature of the filament induces corrections at
lower order due to the second $z$ derivative of the vorticity $\mathbf{%
\omega }_{0}$ in (\ref{lo1}) in the Laplace operator. Since%
\begin{equation*}
\rho =\frac{\left( \mathbf{x}-\mathbf{x}_{0}(s,t)\right) \cdot \mathbf{e}%
_{\rho }}{(\nu t)^{\frac{1}{2}}},
\end{equation*}%
we obtain%
\begin{equation*}
\frac{d\rho }{ds}=-\frac{\mathbf{x}_{0,s}(s,t)\cdot \mathbf{e}_{\rho }}{(\nu
t)^{\frac{1}{2}}},
\end{equation*}%
which vanishes as $s\rightarrow s_{0}$, while%
\begin{equation*}
\frac{d^{2}\rho }{ds^{2}}=-\frac{\mathbf{x}_{0,ss}(s,t)\cdot \mathbf{e}%
_{\rho }}{(\nu t)^{\frac{1}{2}}}=-\kappa \frac{(\mathbf{e}_{\rho }\cdot
\mathbf{n})}{(\nu t)^{\frac{1}{2}}}=-\kappa \frac{\cos \theta }{(\nu t)^{%
\frac{1}{2}}},
\end{equation*}%
and bearing in mind that $dz=(1+O(s^{2}))ds$ we conclude
\begin{eqnarray}
\nu \frac{\partial ^{2}\boldsymbol{\omega }^{0}}{\partial z^{2}} &=&\frac{1}{%
(\nu t)}\frac{d^{2}\rho }{ds^{2}}\Omega _{0}^{\prime }\left( \rho \right)
\mathbf{t}+O((\nu t)^{-1})  \notag \\
&=&-\frac{\nu }{(\nu t)^{\frac{3}{2}}}\frac{d}{ds}(\mathbf{e}_{r}\cdot
\mathbf{x}_{0,s})\Omega _{0}^{\prime }\left( \rho \right) \mathbf{t}+O((\nu
t)^{-1})  \notag \\
&=&-\frac{\nu }{(\nu t)^{\frac{3}{2}}}\kappa \Omega _{0}^{\prime }\left(
\rho \right) \cos \theta \mathbf{t}+O((\nu t)^{-1}),  \label{c4}
\end{eqnarray}%
which is also a $O((\nu t)^{-\frac{3}{2}})$ contribution. This concludes the
proof of the Lemma.

We notice that the cancellation of $\boldsymbol{\omega }_{0}$ outside the
vortex tube implies the cancellation of the nonlinear terms $\boldsymbol{%
\omega }^{0}\cdot \nabla \mathbf{v}^{0}$ and $\mathbf{v}^{0}\cdot \nabla
\mathbf{\omega }^{0}$ outside the tube. We finally remark that in terms of
the vortex strength $\Gamma $, since $\Omega _{0}\left( \rho \right) $ is $%
O(\Gamma )$, the correction (\ref{c4}) is $O(\nu \Gamma )$ while the
corrections (\ref{c1b}), (\ref{c2b}) are all $O(\Gamma ^{2})$. We will
assume, in the next sections, that $\Gamma /\nu $ is sufficiently small so
that any lower order contribution to the vorticity is dominated by the (\ref%
{c4}) correction. Notice also that the fact that (\ref{c1b})-(\ref{c3}) are $%
O((\nu t)^{-\frac{3}{2}})$ are along the tangential direction suggests the
introduction of a correction $\boldsymbol{\omega }^{1}$ to the vorticity
which is $O((\nu t)^{-\frac{1}{2}})$ and along the tangential direction and,
correspondingly $\mathbf{v}^{1}$ which is a $O(1)$ correction to the
velocity.

\section{First order correction}

We introduce a first order correction to the Lamb-Oseen vortex filament (\ref%
{lo1}) in the form
\begin{equation*}
\omega _{z}^{1}=\frac{1}{(\nu t)^{\frac{1}{2}}}\Omega _{1}\left( \frac{r}{%
(\nu t)^{\frac{1}{2}}},\theta \right) ,
\end{equation*}%
and split it as%
\begin{equation*}
\Omega _{1}\left( \rho ,\theta \right) =\Omega _{1}^{(1)}\left( \rho ,\theta
\right) +\Omega _{1}^{(2)}\left( \rho ,\theta \right) ,
\end{equation*}%
with $\Omega _{1}^{(1)}\left( \rho ,\theta \right) $ being the vorticity
correction associated to the local curvature of the filament, that is the
solution to
\begin{equation}
-\frac{1}{2}\Omega _{1}^{(1)}-\frac{1}{2}\rho \frac{\partial \Omega
_{1}^{(1)}}{\partial \rho }-\left( \frac{1}{\rho }\frac{\partial }{\partial
\rho }\left( \rho \frac{\partial \Omega _{1}^{(1)}}{\partial \rho }\right) +%
\frac{1}{\rho ^{2}}\frac{\partial ^{2}\Omega _{1}^{(1)}}{\partial \theta ^{2}%
}\right) =-\kappa \frac{d\Omega _{0}(\rho )}{d\rho }\cos \theta ,
\label{correc1}
\end{equation}%
where we have introduced (\ref{c4}) at the right hand side.

\begin{lemma}
\label{lem1}The solution to (\ref{correc1}) bounded at the origin is%
\begin{equation*}
\Omega _{1}^{(1)}(\rho ,\theta )=\frac{1}{2}\kappa \frac{\Gamma }{4\pi }\rho
e^{-\frac{\rho ^{2}}{4}}\cos \theta ,
\end{equation*}%
and the associated velocity field is
\begin{eqnarray*}
V_{1r}^{(1)} &=&\frac{1}{2}\kappa \frac{\Gamma }{4\pi }\left( \frac{1}{2\rho
^{2}}\left( 8e^{-\frac{1}{4}\rho ^{2}}+2\rho ^{2}e^{-\frac{1}{4}\rho
^{2}}-8\right) -e^{-\frac{1}{4}\rho ^{2}}\right) \sin \theta ,\  \\
V_{1\theta }^{(1)} &=&-\frac{1}{2}\kappa \frac{\Gamma }{4\pi }\frac{1}{\rho
^{2}}\left( 4e^{-\frac{1}{4}\rho ^{2}}+2\rho ^{2}e^{-\frac{1}{4}\rho
^{2}}-4\right) \cos \theta ,
\end{eqnarray*}%
which is bounded and regular at the origin.
\end{lemma}

\bigskip It is simple to find the solution to (\ref{correc1}) in the form
\begin{equation*}
\Omega _{1}^{(1)}(\rho ,\theta )=\frac{1}{2}\kappa \frac{\Gamma }{4\pi }\rho
e^{-\frac{\rho ^{2}}{4}}\cos \theta .
\end{equation*}%
In order to find the velocity field associated to the vorticity $\Omega
_{1}^{(1)}$ we use the definition of vorticity and the condition of
incompressibility:
\begin{eqnarray*}
\frac{1}{\rho }\frac{\partial }{\partial \rho }\left( \rho V_{1\theta
}^{(1)}\right) -\frac{1}{\rho }\frac{\partial V_{1r}^{(1)}}{\partial \theta }
&=&\frac{1}{2}\kappa \frac{\Gamma }{4\pi }\rho e^{-\frac{\rho ^{2}}{4}}\cos
\theta , \\
\frac{1}{\rho }\frac{\partial }{\partial \rho }\left( \rho
V_{1r}^{(1)}\right) +\frac{1}{\rho }\frac{\partial V_{1\theta }^{(1)}}{%
\partial \theta } &=&0,
\end{eqnarray*}%
which amounts to applying Lemma \ref{lemma1} and this concludes the proof of
the Lemma.

We can compute now $\Omega _{1}^{(2)}\left( \rho ,\theta \right) $ as the
solution to the equation

\begin{eqnarray}
&&-\frac{1}{2}\Omega _{1}^{(2)}-\frac{1}{2}\rho \frac{\partial \Omega
_{1}^{(2)}}{\partial \rho }-\left( \frac{1}{\rho }\frac{\partial }{\partial
\rho }\left( \rho \frac{\partial \Omega _{1}^{(2)}}{\partial \rho }\right) +%
\frac{1}{\rho ^{2}}\frac{\partial ^{2}\Omega _{1}^{(2)}}{\partial \theta ^{2}%
}\right)  \notag \\
&=&F_{z}-\frac{1}{\nu }\frac{V_{0}(\rho )}{\rho }\frac{\partial \Omega
_{1}^{(2)}}{\partial \theta }-\frac{1}{\nu }V_{1r}^{(2)}(\rho ,\theta )\frac{%
\partial \Omega _{0}(\rho )}{\partial \rho },  \label{omega}
\end{eqnarray}%
where $F_{z}$ contains all nonlinear terms (\ref{c1b})-(\ref{c3}) together
with the nonlinear terms from $\Omega _{1z}^{(1)}$ and those due to the
velocity field created by distant parts of the filament (\ref{c3}):%
\begin{eqnarray}
F_{z} &=&-\frac{\kappa }{\nu }\Omega _{0}\left( \rho \right) V_{0}(\rho
)\sin \theta  \notag \\
&&+\frac{\Gamma }{\nu }\kappa (4\pi (1-\log 2)+F(\rho )+G(\rho )+H(\rho ))%
\frac{d\Omega _{0}\left( \rho \right) }{d\rho }\sin \theta  \notag \\
&&-\frac{1}{\nu }\frac{V_{0}(\rho )}{\rho }\frac{\partial \Omega _{1}^{(1)}}{%
\partial \theta }-\frac{1}{\nu }V_{1r}^{(1)}(\rho ,\theta )\frac{\partial
\Omega _{0}(\rho )}{\partial \rho }  \notag \\
&&-\frac{1}{\nu }(v_{b}^{\ast }\sin \theta +v_{n}^{\ast }\cos \theta )\frac{%
d\Omega _{0}\left( \rho \right) }{d\rho },  \label{efe}
\end{eqnarray}%
where $(v_{n}^{\ast },v_{b}^{\ast })$ are the normal and binormal components
of $\mathbf{v}^{\ast }$ defined \ by (\ref{vb}).

\begin{lemma}
\label{lem2}\bigskip For $\frac{\Gamma }{\nu }<c_{0}$ with $c_{0}$
sufficiently small and independent of $\Gamma $ and $\nu $, there exists a
unique bounded solution to (\ref{omega}) with $F_{z}$ given by (\ref{efe})
and such that%
\begin{equation*}
\Omega _{1}^{(2)}(\rho ,\theta )=\Omega _{1}^{c(2)}(\rho )\cos \theta
+\Omega _{1}^{s(2)}(\rho )\sin \theta ,
\end{equation*}%
with%
\begin{equation}
\sup \frac{e^{\frac{\rho ^{2}}{4}}\left\vert \Omega _{1}^{c,s(2)}\right\vert
}{\rho +\rho ^{2}}\leq C\frac{\Gamma ^{2}}{\nu },  \label{omega1}
\end{equation}%
where the constant $C$ depends solely on $\kappa $, $v_{b}^{\ast }$ and $%
v_{n}^{\ast }$.
\end{lemma}

\textbf{Proof.} It is simple to show that $F_{z}$ given in (\ref{efe}) is of
the form%
\begin{equation*}
F_{z}(\rho ,\theta )=F_{z}^{c}(\rho )\cos \theta +F_{z}^{s}(\rho )\sin
\theta ,
\end{equation*}%
with%
\begin{equation}
\left\vert F_{z}^{c,s}(\rho )\right\vert \leq (\rho +\rho ^{M})e^{-\frac{%
\rho ^{2}}{4}},  \label{estf}
\end{equation}%
where $M\geq 2$. We consider now equation (\ref{omega}) with a right hand
side of the form%
\begin{equation*}
f(\rho ,\theta )=f^{c}(\rho )\cos \theta +f^{s}(\rho )\sin \theta ,
\end{equation*}%
and seek for solutions of the form
\begin{equation*}
\Omega _{1}^{(2)}(\rho ,\theta )=\Omega _{1}^{c(2)}(\rho )\cos \theta
+\Omega _{1}^{s(2)}(\rho )\sin \theta ,
\end{equation*}%
so that $y(\rho )\equiv \Omega _{1}^{c,s(2)}(\rho )$ satisfies the ODE%
\begin{equation}
-\frac{1}{2}y(\rho )-\frac{1}{2}\rho \frac{dy(\rho )}{d\rho }-\left( \frac{1%
}{\rho }\frac{d}{d\rho }\left( \rho \frac{dy(\rho )}{d\rho }\right) -\frac{%
y(\rho )}{\rho ^{2}}\right) =f(\rho ),  \label{ode}
\end{equation}%
where $f(\rho )$ stands for $f^{c}(\rho )$ or $f^{s}(\rho )$. Following the
variation of constants method we write $y(\rho )=y_{0}(\rho )z(\rho )$ with
\begin{equation*}
-\frac{1}{2}y_{0}(\rho )-\frac{1}{2}\rho \frac{dy_{0}(\rho )}{d\rho }-\left(
\frac{1}{\rho }\frac{d}{d\rho }\left( \rho \frac{dy_{0}(\rho )}{d\rho }%
\right) -\frac{y_{0}(\rho )}{\rho ^{2}}\right) =0,
\end{equation*}%
that is%
\begin{equation*}
y_{0}(\rho )=\frac{1}{\rho }e^{-\frac{\rho ^{2}}{4}}.
\end{equation*}%
Hence%
\begin{equation*}
\left( \frac{1}{2}\rho +\frac{1}{\rho }+2\frac{y_{0}^{\prime }(\rho )}{%
y_{0}(\rho )}\right) z^{\prime }(\rho )+z^{\prime \prime }(\rho )=\frac{%
f(\rho )}{y_{0}(\rho )},
\end{equation*}%
or equivalently%
\begin{equation*}
\left( e^{\int \left( \frac{1}{2}\rho +\frac{1}{\rho }+2\frac{y_{0}^{\prime
}(\rho )}{y_{0}(\rho )}\right) }z^{\prime }(\rho )\right) ^{\prime }=e^{\int
\left( \frac{1}{2}\rho +\frac{1}{\rho }+2\frac{y_{0}^{\prime }(\rho )}{%
y_{0}(\rho )}\right) }\frac{f(\rho )}{y_{0}(\rho )},
\end{equation*}%
so that%
\begin{equation*}
\left( e^{\frac{\rho ^{2}}{4}}\rho y_{i}^{2}(\rho )z^{\prime }(\rho )\right)
^{\prime }=e^{\frac{\rho ^{2}}{4}}\rho y_{i}(\rho )f(\rho ),
\end{equation*}%
and then%
\begin{equation}
z^{\prime }(\rho )=-\frac{e^{-\frac{\rho ^{2}}{4}}}{\rho y_{0}^{2}(\rho )}%
\int_{\rho }^{\infty }e^{\frac{\rho ^{\prime 2}}{4}}\rho ^{\prime
}y_{0}(\rho ^{\prime })f(\rho ^{\prime })d\rho ^{\prime },  \label{zp}
\end{equation}%
concluding%
\begin{equation}
z(\rho )=-\int_{0}^{\rho }\left( \frac{e^{-\frac{\rho ^{\prime 2}}{4}}}{\rho
^{\prime }y_{0}^{2}(\rho ^{\prime })}\int_{\rho ^{\prime }}^{\infty }e^{%
\frac{\rho ^{\prime \prime 2}}{4}}\rho ^{\prime \prime }y_{0}(\rho ^{\prime
\prime })f(\rho ^{\prime \prime })d\rho ^{\prime \prime }\right) d\rho
^{\prime }.  \label{z}
\end{equation}%
Since

\begin{equation}
\left\vert f(\rho )\right\vert \leq (\rho +\rho ^{M})e^{-\frac{\rho ^{2}}{4}%
},  \label{est1}
\end{equation}%
we have%
\begin{equation*}
\frac{e^{-\frac{\rho ^{2}}{4}}}{\rho y_{0}^{2}(\rho )}\int_{\rho }^{\infty
}e^{\frac{\rho ^{\prime 2}}{4}}\rho ^{\prime }y_{0}(\rho ^{\prime
})\left\vert f(\rho ^{\prime })\right\vert d\rho ^{\prime }\leq C(\rho +\rho
^{M}),
\end{equation*}%
and therefore%
\begin{equation}
\left\vert y(\rho )\right\vert =\left\vert y_{0}(\rho )z(\rho )\right\vert
\leq C(\rho +\rho ^{M})e^{-\frac{\rho ^{2}}{4}}.  \label{est2}
\end{equation}%
Coming back to equation (\ref{omega}), notice the presence of the term $%
V_{1r}^{(2)}(\rho ,\theta )$ at the right hand side. This term is related to
the vorticity $\Omega _{1}^{(2)}(\rho ,\theta )$ by the relations
\begin{eqnarray*}
\Delta \psi &=&\Omega _{1}^{(2)}(\rho ,\theta ), \\
V_{1r}^{(2)}(\rho ,\theta ) &=&-\frac{1}{\rho }\frac{\partial \psi }{%
\partial \theta },
\end{eqnarray*}%
so that, if%
\begin{equation*}
\Omega _{1}^{(2)}(\rho ,\theta )=\Omega _{1}^{c(2)}(\rho )\cos \theta
+\Omega _{1}^{s(2)}(\rho )\sin \theta ,
\end{equation*}%
then%
\begin{equation*}
\psi (\rho ,\theta )=\psi ^{c}(\rho )\cos \theta +\psi ^{s}(\rho )\sin
\theta ,
\end{equation*}%
and hence%
\begin{equation*}
\frac{1}{\rho }\frac{d}{d\rho }\left( \rho \frac{d\psi ^{c,s}(\rho )}{d\rho }%
\right) -\frac{\psi ^{c,s}(\rho )}{\rho ^{2}}=\Omega _{1z}^{c,s(2)}(\rho ).
\end{equation*}%
We write%
\begin{equation*}
\psi ^{c,s}(\rho )=\rho \varphi ^{c,s}(\rho ),
\end{equation*}%
so that%
\begin{equation*}
\frac{1}{\rho }(\rho ^{3}\varphi ^{c,s\prime }(\rho ))^{\prime }=\Omega
_{1}^{c,s(2)}(\rho ),
\end{equation*}%
and then%
\begin{equation*}
\varphi ^{c,s\prime }(\rho )=-\int_{\rho }^{\infty }\rho ^{\prime -3}\left(
\int_{0}^{\rho ^{\prime }}\rho ^{\prime \prime }\Omega _{1}^{c,s(2)}(\rho
^{\prime \prime })d\rho ^{\prime \prime }\right) d\rho ^{\prime }
\end{equation*}%
leading to%
\begin{eqnarray*}
V_{1r}^{(2)}(\rho ,\theta ) &=&\left( \int_{\rho }^{\infty }\rho ^{\prime
-3}\left( \int_{0}^{\rho ^{\prime }}\rho ^{\prime \prime }\Omega
_{1}^{c(2)}(\rho ^{\prime \prime })d\rho ^{\prime \prime }\right) d\rho
^{\prime }\right) \sin \theta \\
&&-\left( \int_{\rho }^{\infty }\rho ^{\prime -3}\left( \int_{0}^{\rho
^{\prime }}\rho ^{\prime \prime }\Omega _{1}^{s(2)}(\rho ^{\prime \prime
})d\rho ^{\prime \prime }\right) d\rho ^{\prime }\right) \cos \theta \\
&\equiv &V_{1r}^{s(2)}(\rho )\sin \theta +V_{1r}^{c(2)}(\rho )\cos \theta .
\end{eqnarray*}%
Similarly,%
\begin{eqnarray*}
V_{1\theta }^{(2)}(\rho ,\theta ) &=&\frac{d}{d\rho }\left( \rho \int_{\rho
}^{\infty }\rho ^{\prime -3}\left( \int_{0}^{\rho ^{\prime }}\rho ^{\prime
\prime }\Omega _{1}^{c(2)}(\rho ^{\prime \prime })d\rho ^{\prime \prime
}\right) d\rho ^{\prime }\right) \cos \theta \\
&&+\frac{d}{d\rho }\left( \rho \int_{\rho }^{\infty }\rho ^{\prime -3}\left(
\int_{0}^{\rho ^{\prime }}\rho ^{\prime \prime }\Omega _{1}^{s(2)}(\rho
^{\prime \prime })d\rho ^{\prime \prime }\right) d\rho ^{\prime }\right)
\sin \theta \\
&\equiv &V_{1\theta }^{c(2)}(\rho )\cos \theta +V_{1\theta }^{s(2)}(\rho
)\sin \theta .
\end{eqnarray*}%
One can then easily get%
\begin{equation*}
\left\vert \int_{0}^{\rho ^{\prime }}\rho ^{\prime \prime }\Omega
_{1}^{c,s(2)}(\rho ^{\prime \prime })d\rho ^{\prime \prime }\right\vert \leq
\frac{C\rho ^{\prime 3}}{(1+\rho ^{\prime })^{3}}\sup \frac{e^{\frac{\rho
^{2}}{4}}\left\vert \Omega _{1}^{c,s(2)}\right\vert }{\rho +\rho ^{2}},
\end{equation*}%
and hence%
\begin{equation*}
\left\vert V_{1r}^{c,s(2)}\right\vert \leq C\frac{1}{1+\rho }\sup \frac{e^{%
\frac{\rho ^{2}}{4}}\left\vert \Omega _{1}^{c,s(2)}\right\vert }{\rho +\rho
^{2}}.
\end{equation*}

Finally, in order to solve (\ref{omega}) we use a fixed point argument.
Consider the equation%
\begin{eqnarray}
&&-\frac{1}{2}\Omega _{1}^{(2)}-\frac{1}{2}\rho \frac{\partial \Omega
_{1}^{(2)}}{\partial \rho }-\left( \frac{1}{\rho }\frac{\partial }{\partial
\rho }\left( \rho \frac{\partial \Omega _{1}^{(2)}}{\partial \rho }\right) +%
\frac{1}{\rho ^{2}}\frac{\partial ^{2}\Omega _{1}^{(2)}}{\partial \theta ^{2}%
}\right)  \notag \\
&=&F_{z}-\frac{1}{\nu }\frac{V_{0}(\rho )}{\rho }\frac{\partial \widetilde{%
\Omega }_{1}^{(2)}}{\partial \theta }-\frac{1}{\nu }\widetilde{V}%
_{1r}^{(2)}(\rho ,\theta )\frac{\partial \Omega _{0}(\rho )}{\partial \rho },
\label{omegap}
\end{eqnarray}%
where $\widetilde{\Omega }_{1}^{(2)}$ is a linear combination of sine and
cosine functions in $\theta $ with coefficients in the unit ball under the
topology
\begin{equation*}
\left\Vert \widetilde{\Omega }_{1}^{c,s(2)}\right\Vert =\sup \frac{e^{\frac{%
\rho ^{2}}{4}}\left\vert \widetilde{\Omega }_{1}^{c,s(2)}\right\vert }{\rho
+\rho ^{M}}.
\end{equation*}%
We define the mapping $T$ that maps $\widetilde{\Omega }_{1}^{c,s(2)}$ into
the solution $\Omega _{1}^{c,s(2)}$ associated to (\ref{omegap}). Using (\ref%
{est1}), (\ref{est2}) we can show the estimate
\begin{equation*}
\left\Vert \Omega _{1}^{c(2)}\right\Vert +\left\Vert \Omega
_{1}^{s(2)}\right\Vert \leq C\frac{\Gamma ^{2}}{\nu }+C\frac{\Gamma }{\nu }%
\left( \left\Vert \widetilde{\Omega ^{c}}_{1}^{(2)}\right\Vert +\left\Vert
\widetilde{\Omega ^{s}}_{1}^{(2)}\right\Vert \right) ,
\end{equation*}%
so that $T$ maps the closed unit ball into itself provided $\frac{\Gamma }{%
\nu }$ is sufficiently small. By linearity it is also a contractive mapping
for $\Gamma $ sufficiently small. By Banach's fixed point theorem there is a
unique solution to problem (\ref{omega}). We can obtain further regularity
for $\Omega _{1}^{c(2)}$ and $\Omega _{1}^{s(2)}$. By denoting $y(\rho
)=\Omega _{1}^{c(2)}$ or $y(\rho )=\Omega _{1}^{s(2)}$ it turns out that $%
y(\rho )$ satisfies (\ref{ode}) with $f(\rho )$ satisfying (\ref{est1}) and
hence $y(\rho )=y_{0}(\rho )z(\rho )$ with $z(\rho )$ given by (\ref{z}).
Straightforward calculation yields $y^{\prime }(\rho )=O(1)$ and $y^{\prime
\prime }(\rho )=O(\rho )$ as $\rho \rightarrow 0$ as well as the estimates%
\begin{equation}
\sup \frac{e^{\frac{\rho ^{2}}{4}}\left\vert \frac{d\widetilde{\Omega }%
_{1}^{c,s(2)}(\rho )}{d\rho }\right\vert }{1+\rho ^{M+1}}\leq C\frac{\Gamma
^{2}}{\nu },\ \sup \frac{e^{\frac{\rho ^{2}}{4}}\left\vert \frac{d^{2}%
\widetilde{\Omega }_{1}^{c,s(2)}(\rho )}{d\rho ^{2}}\right\vert }{\rho +\rho
^{M+2}}\leq C\frac{\Gamma ^{2}}{\nu },  \label{omega2}
\end{equation}%
and this concludes the proof of the Lemma.

Note that $\Omega _{1}^{c}$ and $\Omega _{1}^{s}$ depend parametrically on $%
s $ through the curvature $\kappa $ and the velocities $v_{b}^{\ast }$,$%
v_{n}^{\ast }$ so that one can write $\Omega _{1}^{c}(\rho ,s)$ and $\Omega
_{1}^{s}(\rho ,s)$ respectively. We can compute the $s$-derivative of
equation (\ref{omega}) with $F_{z,s}$ computed from (\ref{efe}) and hence
involving $\kappa _{s}$, $v_{b,s}^{\ast }$,$v_{n,s}^{\ast }$ to obtain the
same equation (\ref{omega}) for $\Omega _{1,s}^{c,s(2)}(\rho ,s)$ (with new
terms in the force term $F_{z,s}$ but such that estimate (\ref{estf}) still
applies) and hence identical estimates (\ref{omega1}) and (\ref{omega2})
with $C$ depending on $\kappa _{s}$, $v_{b,s}^{\ast }$ and $v_{n,s}^{\ast }$%
. Identical argument also applies to $\Omega _{1,ss}^{c,s(2)}(\rho ,s)$ with
$C$ depending on $\kappa _{ss}$, $v_{b,ss}^{\ast }$ and $v_{n,ss}^{\ast }.$

We finally remark that, by construction in lemmas \ref{lem1}, \ref{lem2}, if
we define%
\begin{equation*}
\omega _{0}=\frac{1}{(\nu t)}\Omega _{0}(\rho )+\frac{1}{(\nu t)^{\frac{1}{2}%
}}\left( \Omega _{1}^{c}(\rho ,s)\cos \theta +\Omega _{1}^{s}(\rho ,s)\sin
\theta \right) ,
\end{equation*}%
and%
\begin{equation*}
\boldsymbol{\omega }_{0}=\eta _{R}\omega _{0}\mathbf{t},
\end{equation*}%
as well as the velocity field
\begin{equation*}
\mathbf{v}_{0}(\mathbf{x},t)=K\ast \boldsymbol{\omega }_{0}(\mathbf{x},t)+%
\frac{\Gamma }{4\pi }\kappa (s,t)\log (\nu t)^{\frac{1}{2}}\mathbf{b}(s,t),
\end{equation*}%
with $\eta _{R}$ a cutoff function defined by (\ref{co0})-(\ref{co2}) and $K$
being the Biot-Savart kernel, then the quantity $A\left[ \omega
_{0},v_{0r},v_{0\theta }\right] $ (where $v_{0r},v_{0\theta }$ are the $r$
and $\theta $ components of $\mathbf{v}_{0}(\mathbf{x},t)$ respectively)
defined by (\ref{defa}) vanishes at the $O((\nu t)^{-2})$ and $O((\nu t)^{-%
\frac{3}{2}})$ orders inside the tube of radius $R$ and satisfies%
\begin{equation}
\left\vert A\left[ \omega _{0},v_{0r},v_{0\theta }\right] \right\vert \leq
\frac{C(\Gamma ^{2}+\Gamma \nu )}{(\nu t)}(1+\rho ^{4})e^{-\frac{\rho ^{2}}{4%
}}.  \label{estia}
\end{equation}

\section{Estimates of $\mathbf{F}(\mathbf{x},t)\label{secf}$}

In this section we will estimate the error $\mathbf{F}(\mathbf{x},t)$
produced when introducing the approximate solution constructed in previous
sections in Navier-Stokes system, i.e. the right hand side of (\ref{eqnomega}%
). We will find more convenient to compute the various terms in the
Navier-Stokes system in a system of coordinates $(x_{1},x_{2},s)$ adapted to
the filament geometry (the parallel frame) and related to the position
vector in the fixed reference frame $\mathbf{x}$ by the relation deduced in (%
\ref{xpara}):
\begin{equation*}
\mathbf{x=x}_{0}(s,t)+x_{1}\mathbf{e}_{1}(s,t)+x_{2}\mathbf{e}_{2}(s,t).
\end{equation*}%
We will denote $\mathbf{e}_{s}$ the vector orthogonal to $(\mathbf{e}_{1},%
\mathbf{e}_{2})$ i.e. the tangent vector $\mathbf{t}$. If we introduce%
\begin{equation*}
\widetilde{\mathbf{x}}=\mathbf{x-x}_{0}(s,t),
\end{equation*}%
and the vector fields $\widetilde{\boldsymbol{\omega }}_{0}(\widetilde{%
\mathbf{x}},t)=\boldsymbol{\omega }_{0}(\mathbf{x},t)$ and $\widetilde{%
\mathbf{v}}_{0}(\widetilde{\mathbf{x}},t)$ the velocity vector field as
measured in the parallel frame, i.e. by%
\begin{eqnarray*}
\mathbf{v}_{0}(\mathbf{x},t) &=&\mathbf{x}_{0,t}(s,t)+\widetilde{\mathbf{v}}%
_{0}(\widetilde{\mathbf{x}},t) \\
&&+x_{1}\left( \mathbf{e}_{1,t}(s,t)+v_{s}\mathbf{e}_{1,s}(s,t)\right)
+x_{2}\left( \mathbf{e}_{2,t}(s,t)+v_{s}\mathbf{e}_{2,s}(s,t)\right) ,
\end{eqnarray*}%
where $v_{s}$ is the component of the velocity in the tangential direction,
we can prove the following Lemma:

\begin{lemma}
Let%
\begin{equation*}
\boldsymbol{\omega }_{0}(\mathbf{x},t)=\eta _{R}\omega _{0}\mathbf{e}_{s},
\end{equation*}%
with
\begin{equation*}
\omega _{0}=\frac{1}{(\nu t)}\Omega _{0}(\rho )+\frac{1}{(\nu t)^{\frac{1}{2}%
}}\left( \Omega _{1}^{c}(\rho ,s)\cos \theta +\Omega _{1}^{s}(\rho ,s)\sin
\theta \right) ,
\end{equation*}%
and%
\begin{equation*}
\mathbf{v}_{0}(\mathbf{x},t)=K\ast \boldsymbol{\omega }_{0}(\mathbf{x},t).
\end{equation*}%
Then
\begin{eqnarray*}
&&\frac{\partial \boldsymbol{\omega }_{0}(\mathbf{x},t)}{\partial t}+\mathbf{%
v}_{0}(\mathbf{x},t)\cdot \nabla \boldsymbol{\omega }_{0}(\mathbf{x},t)-%
\boldsymbol{\omega }_{0}(\mathbf{x},t)\cdot \nabla \mathbf{v}_{0}(\mathbf{x}%
,t) \\
&=&\frac{\partial \widetilde{\boldsymbol{\omega }}_{0}(\widetilde{\mathbf{x}}%
,t)}{\partial t}+\widetilde{\mathbf{v}}_{0}(\widetilde{\mathbf{x}},t)\cdot
\nabla \widetilde{\boldsymbol{\omega }}_{0}(\widetilde{\mathbf{x}},t)-%
\widetilde{\mathbf{\omega }}_{0}(\widetilde{\mathbf{x}},t)\cdot \nabla
\widetilde{\mathbf{v}}_{0}(\widetilde{\mathbf{x}},t) \\
&&+\mathbf{D}(\widetilde{\mathbf{x}},t),
\end{eqnarray*}%
with%
\begin{equation*}
\left\vert \mathbf{D}(\widetilde{\mathbf{x}},t)\right\vert \leq \frac{%
C\Gamma ^{2}\left\vert \log (\nu t)\right\vert }{(\nu t)}(1+\rho ^{4})e^{-%
\frac{\rho ^{2}}{4}}.
\end{equation*}
\end{lemma}

\textbf{Proof}. We denote $\mathbf{x}(t)$ the trajectory of a fluid particle
in the fixed frame and $(x_{1}(t),x_{2}(t),s(t))$ the trajectory written in
the $(\mathbf{e}_{1}(s,t),\mathbf{e}_{2}(s,t),\mathbf{t}(s,t))$ frame. Note
that
\begin{equation*}
\mathbf{x}(t)=\mathbf{x}_{0}(s(t),t)+x_{1}(t)\mathbf{e}_{1}(s(t),t)+x_{2}(t)%
\mathbf{e}_{2}(s(t),t),
\end{equation*}%
so that%
\begin{equation*}
\mathbf{x}(t)-\mathbf{x}_{0}(s(t),t)=x_{1}(t)\mathbf{e}_{1}(s(t),t)+x_{2}(t)%
\mathbf{e}_{2}(s(t),t).
\end{equation*}%
Since%
\begin{equation*}
\boldsymbol{\omega }_{0}(\mathbf{x}(t),t)=\widetilde{\boldsymbol{\omega }}%
_{0}(\mathbf{x}(t)-\mathbf{x}_{0}(s(t),t),t),
\end{equation*}%
then by (\ref{vpf}) we have%
\begin{eqnarray*}
&&\frac{d\boldsymbol{\omega }_{0}(\mathbf{x}(t),t)}{dt} \\
&=&\frac{\partial \widetilde{\boldsymbol{\omega }}_{0}(\widetilde{\mathbf{x}}%
(t),t)}{\partial t}+\left( x_{1,t}(t)\mathbf{e}_{1}(s(t),t)+x_{2,t}(t)%
\mathbf{e}_{2}(s(t),t)\right) \cdot \nabla \widetilde{\boldsymbol{\omega }}%
_{0}(\widetilde{\mathbf{x}}(t),t) \\
&&+(x_{1}(t)\left( \mathbf{e}_{1,t}(s(t),t)+v_{s}(t)\mathbf{e}%
_{1,s}(s(t),t)\right) )\cdot \nabla \widetilde{\boldsymbol{\omega }}_{0}(%
\widetilde{\mathbf{x}}(t),t) \\
&&+(x_{2}(t)\left( \mathbf{e}_{2,t}(s(t),t)+v_{s}(t)\mathbf{e}%
_{2,s}(s(t),t)\right) )\cdot \nabla \widetilde{\boldsymbol{\omega }}_{0}(%
\widetilde{\mathbf{x}}(t),t).
\end{eqnarray*}%
Using the following relation between the velocities in both frames:%
\begin{eqnarray}
\mathbf{v}_{0}(\mathbf{x},t) &=&\mathbf{x}_{0,t}(s,t)+\widetilde{\mathbf{v}}%
_{0}(\widetilde{\mathbf{x}},t)  \notag \\
&&+x_{1}\left( \mathbf{e}_{1,t}(s,t)+v_{s}\mathbf{e}_{1,s}(s,t)\right)
+x_{2}\left( \mathbf{e}_{2,t}(s,t)+v_{0,s}\mathbf{e}_{2,s}(s,t)\right) ,
\label{relvel}
\end{eqnarray}%
we have%
\begin{eqnarray}
&&\frac{\partial \boldsymbol{\omega }_{0}(\mathbf{x},t)}{\partial t}+\mathbf{%
v}_{0}(\mathbf{x},t)\cdot \nabla \boldsymbol{\omega }_{0}(\mathbf{x},t)
\notag \\
&=&\frac{\partial \widetilde{\boldsymbol{\omega }}_{0}(\widetilde{\mathbf{x}}%
,t)}{\partial t}+\widetilde{\mathbf{v}}_{0}(\widetilde{\mathbf{x}},t)\cdot
\nabla \widetilde{\boldsymbol{\omega }}_{0}(\widetilde{\mathbf{x}},t)  \notag
\\
&&+(-v_{0,s}\mathbf{t}(s,t)+x_{1}\left( \mathbf{e}_{1,t}+v_{s}\mathbf{e}%
_{1,s}\right) +x_{2}(\mathbf{e}_{2,t}+v_{0,s}\mathbf{e}_{2,s}))\cdot \nabla
\widetilde{\boldsymbol{\omega }}_{0}(\widetilde{\mathbf{x}},t)  \label{ff1}
\end{eqnarray}

Concerning the stretching term $\boldsymbol{\omega }_{0}\cdot \nabla \mathbf{%
v}_{0}$, since the velocities in both frames are related by (\ref{relvel}),
we have%
\begin{eqnarray}
\boldsymbol{\omega }_{0}(\mathbf{x},t)\cdot \nabla \mathbf{v}_{0}(\mathbf{x}%
,t) &=&\widetilde{\boldsymbol{\omega }}_{0}(\widetilde{\mathbf{x}},t)\cdot
\nabla \widetilde{\mathbf{v}}_{0}(\widetilde{\mathbf{x}},t)  \notag \\
&&+\widetilde{\boldsymbol{\omega }}_{0}(\widetilde{\mathbf{x}},t)\cdot
\nabla \left( \mathbf{x}_{0,t}+x_{1}\left( \mathbf{e}_{1,t}+v_{0,s}\mathbf{e}%
_{1,s}\right) +x_{2}(\mathbf{e}_{2,t}+v_{0,s}\mathbf{e}_{2,s})\right) .
\label{ff2}
\end{eqnarray}%
We note finally that the last terms at the right hand side of both (\ref{ff1}%
) and (\ref{ff2}), denoted by $\mathbf{D}_{1}$ and $\mathbf{D}_{2}$
respectively, satisfy (by using lemma \ref{lemmae} to estimate derivatives
of $\mathbf{e}_{i}(s,t)$),%
\begin{equation*}
\left\vert \mathbf{D}_{1}\right\vert +\left\vert \mathbf{D}_{2}\right\vert
\leq \frac{C\Gamma ^{2}\left\vert \log (\nu t)\right\vert }{(\nu t)}(1+\rho
^{4})e^{-\frac{\rho ^{2}}{4}}.
\end{equation*}

\begin{lemma}
\label{lemmaa}Let%
\begin{equation}
\omega _{0}=\frac{1}{(\nu t)}\Omega _{0}(\rho )+\frac{1}{(\nu t)^{\frac{1}{3}%
}}\left( \Omega _{1}^{c}(\rho ,s)\cos \theta +\Omega _{1}^{s}(\rho ,s)\sin
\theta \right) .  \label{gg1}
\end{equation}%
Then the functions
\begin{equation*}
\boldsymbol{\omega }_{0}=\eta _{R}\omega _{0}\mathbf{e}_{s},
\end{equation*}%
and
\begin{equation*}
\mathbf{v}_{0}(\mathbf{x},t)=K\ast \boldsymbol{\omega }_{0}(\mathbf{x},t),
\end{equation*}%
with $K$ the Biot-Savart kernel and $s=\arg \min_{s}(dist(\mathbf{x},\mathbf{%
x}_{0}(s,t)))$, are such that%
\begin{equation}
\left\vert \frac{\partial \boldsymbol{\omega }_{0}}{\partial t}+\mathbf{v}%
_{0}\cdot \nabla \boldsymbol{\omega }_{0}-\boldsymbol{\omega }_{0}\cdot
\nabla \mathbf{v}_{0}-\nu \Delta \boldsymbol{\omega }_{0}\right\vert \leq
\frac{C(\Gamma ^{2}+\Gamma \nu )\left\vert \log (\nu t)\right\vert }{(\nu t)}%
(1+\rho ^{4})e^{-\frac{\rho ^{2}}{4}}.  \label{ef1}
\end{equation}%
\bigskip
\end{lemma}

\textbf{Proof.} First, it follows from the previous lemma that estimate (\ref%
{ef1}) holds provided%
\begin{equation*}
\left\vert \frac{\partial \widetilde{\boldsymbol{\omega }}_{0}}{\partial t}+%
\widetilde{\mathbf{v}}_{0}\cdot \nabla \widetilde{\boldsymbol{\omega }}_{0}-%
\widetilde{\boldsymbol{\omega }}_{0}\cdot \nabla \widetilde{\mathbf{v}}%
_{0}-\nu \Delta \widetilde{\boldsymbol{\omega }}_{0}\right\vert \leq \frac{%
C(\Gamma ^{2}+\Gamma \nu )\left\vert \log (\nu t)\right\vert }{(\nu t)}%
(1+\rho ^{4})e^{-\frac{\rho ^{2}}{4}}.
\end{equation*}%
We consider a small $s$-interval around $s_{0}$, place ourselves in a frame
whose center moves with velocity $\mathbf{x}_{0t}(s_{0},t)$ and use the
parallel frame orthonormal system written in polar coordinates (see (\ref%
{erefes}) for the definition of the corresponding unit vectors $\mathbf{e}%
_{r},\mathbf{e}_{\varphi }$ in terms of $\mathbf{e}_{1},\mathbf{e}_{2}$) and
denote for the sake of simplicity $\ (\boldsymbol{e}_{1},\boldsymbol{e}_{2},%
\boldsymbol{e}_{3})=(\mathbf{e}_{r},\mathbf{e}_{\varphi },\mathbf{e}_{s})$
where $\mathbf{e}_{s}=\mathbf{t}$. The polar coordinates associated to $(%
\boldsymbol{e}_{1},\boldsymbol{e}_{2})$ are then denoted as $(r,\varphi )$
and the coordinate associated to $\boldsymbol{e}_{3}$ is the arclength
parameter $s$. The scale factors are (see (\ref{scale1})) $h_{1}=1$, $%
h_{2}=r $, $h_{3}=1-\kappa (s)r\cos (\varphi -\theta _{0}(s))$ and $h_{ij}$
will denote the $j$-derivative of $h_{i}$.\ We also write, for the sake of
simplicity, $\theta =\varphi -\theta _{0}(s)$. We will use indistinctly the
notation $\widetilde{\boldsymbol{\omega }}_{0}=(\omega _{1},\omega
_{2},\omega _{3})$ or $(\omega _{r},\omega _{\varphi },\omega _{s})$ (and
likewise for the velocity field). By hypothesis, we have in this system of
coordinates%
\begin{equation*}
\widetilde{\boldsymbol{\omega }}_{0}=(\omega _{1},\omega _{2},\omega
_{3})=(0,0,\omega _{0}).
\end{equation*}%
We will compute, in this system, the various terms in Navier-Stokes system.
First, the stretching terms
\begin{eqnarray*}
(\widetilde{\boldsymbol{\omega }}_{0}\mathbf{\cdot \nabla })(v_{1}%
\boldsymbol{e}_{1}) &=&\left( \frac{\omega _{1}}{h_{1}}\frac{\partial v_{1}}{%
\partial x_{1}}+\frac{\omega _{2}}{h_{2}}\frac{\partial v_{1}}{\partial x_{2}%
}+\frac{\omega _{3}}{h_{3}}\frac{\partial v_{1}}{\partial x_{3}}\right)
\boldsymbol{e}_{1} \\
&&+\left( \frac{v_{1}\omega _{2}h_{21}}{h_{1}h_{2}}-\frac{v_{1}\omega
_{1}h_{12}}{h_{1}h_{2}}\right) \boldsymbol{e}_{2}+\left( \frac{v_{1}\omega
_{3}h_{31}}{h_{1}h_{3}}-\frac{v_{1}\omega _{1}h_{13}}{h_{1}h_{3}}\right)
\boldsymbol{e}_{3},
\end{eqnarray*}%
where%
\begin{eqnarray*}
h_{21} &=&1, \\
h_{31} &=&-\kappa \cos \theta ,\ h_{32}=\kappa r\sin \theta ,
\end{eqnarray*}%
and, in particular, when $\widetilde{\boldsymbol{\omega }}=\omega _{3}%
\boldsymbol{e}_{3}$, we have the nonzero terms
\begin{eqnarray}
\frac{\omega _{3}}{h_{3}}\frac{\partial v_{1}}{\partial x_{3}}\boldsymbol{e}%
_{1} &=&\frac{\omega _{s}}{1-\kappa r\cos \theta }\frac{\partial v_{r}}{%
\partial s}\mathbf{e}_{r},  \label{f1} \\
\frac{v_{1}\omega _{3}h_{31}}{h_{1}h_{3}}\boldsymbol{e}_{3} &=&-\frac{\kappa
\cos \theta }{1-\kappa r\cos \theta }v_{r}\omega _{s}\mathbf{e}_{s}.
\label{f2}
\end{eqnarray}%
Secondly%
\begin{eqnarray*}
(\widetilde{\boldsymbol{\omega }}_{0}\mathbf{\cdot \nabla })(v_{2}%
\boldsymbol{e}_{2}) &=&\left( \frac{\omega _{1}}{h_{1}}\frac{\partial v_{2}}{%
\partial x_{1}}+\frac{\omega _{2}}{h_{2}}\frac{\partial v_{2}}{\partial x_{2}%
}+\frac{\omega _{3}}{h_{3}}\frac{\partial v_{2}}{\partial x_{3}}\right)
\boldsymbol{e}_{2} \\
&&+\left( \frac{v_{2}\omega _{1}h_{12}}{h_{1}h_{2}}-\frac{v_{2}\omega
_{2}h_{21}}{h_{1}h_{2}}\right) \boldsymbol{e}_{1}+\left( \frac{v_{2}\omega
_{3}h_{32}}{h_{2}h_{3}}-\frac{v_{2}\omega _{2}h_{23}}{h_{2}h_{3}}\right)
\boldsymbol{e}_{3},
\end{eqnarray*}%
where we have the nonzero terms%
\begin{eqnarray}
\frac{\omega _{3}}{h_{3}}\frac{\partial v_{2}}{\partial x_{3}}\boldsymbol{e}%
_{2} &=&\frac{\omega _{s}}{1-\kappa r\cos \theta }\frac{\partial v_{\varphi }%
}{\partial s}\mathbf{e}_{\varphi },  \label{f3} \\
\frac{v_{2}\omega _{3}h_{32}}{h_{2}h_{3}}\boldsymbol{e}_{3} &=&\frac{\kappa
\sin \theta }{1-\kappa r\cos \theta }v_{\varphi }\omega _{s}\mathbf{e}_{s},
\label{f4}
\end{eqnarray}%
and finally%
\begin{eqnarray*}
(\widetilde{\boldsymbol{\omega }}_{0}\mathbf{\cdot \nabla })(v_{3}%
\boldsymbol{e}_{3}) &=&\left( \frac{\omega _{1}}{h_{1}}\frac{\partial v_{3}}{%
\partial x_{1}}+\frac{\omega _{2}}{h_{2}}\frac{\partial v_{3}}{\partial x_{2}%
}+\frac{\omega _{3}}{h_{3}}\frac{\partial v_{3}}{\partial x_{3}}\right)
\boldsymbol{e}_{3} \\
&&+\left( \frac{v_{3}\omega _{1}h_{13}}{h_{3}h_{1}}-\frac{v_{3}\omega
_{3}h_{31}}{h_{1}h_{3}}\right) \boldsymbol{e}_{1}+\left( \frac{v_{3}\omega
_{2}h_{23}}{h_{2}h_{3}}-\frac{v_{3}\omega _{3}h_{32}}{h_{2}h_{3}}\right)
\boldsymbol{e}_{2},
\end{eqnarray*}%
where we have the nonzero terms%
\begin{eqnarray}
\frac{\omega _{3}}{h_{3}}\frac{\partial v_{3}}{\partial x_{3}}\boldsymbol{e}%
_{3} &=&\frac{\omega _{s}}{1-\kappa r\cos \theta }\frac{\partial v_{s}}{%
\partial s}\mathbf{e}_{s},  \label{f5} \\
\frac{v_{3}\omega _{3}h_{31}}{h_{1}h_{3}}\boldsymbol{e}_{1} &=&-\frac{\kappa
\cos \theta }{1-\kappa r\cos \theta }v_{s}\omega _{s}\mathbf{e}_{r},
\label{f6} \\
\frac{v_{3}\omega _{3}h_{32}}{h_{2}h_{3}}\boldsymbol{e}_{2} &=&\frac{\kappa
\sin \theta }{1-\kappa r\cos \theta }v_{s}\omega _{s}\mathbf{e}_{\varphi }.
\label{f7}
\end{eqnarray}%
Analogously, we can compute the convective terms $(\widetilde{\mathbf{v}}_{0}%
\mathbf{\cdot \nabla })(\omega _{3}\boldsymbol{e}_{3})$ by simply exchanging
$v$ and $\omega $ in the formulas above and obtain the nonzero terms%
\begin{eqnarray}
(\widetilde{\mathbf{v}}_{0}\mathbf{\cdot \nabla })(\omega _{3}\boldsymbol{e}%
_{3}) &=&\left( v_{r}\frac{\partial \omega _{s}}{\partial r}+\frac{%
v_{\varphi }}{r}\frac{\partial \omega _{s}}{\partial \varphi }+\frac{v_{s}}{%
1-\kappa r\cos \theta }\frac{\partial \omega _{s}}{\partial s}\right)
\mathbf{e}_{s}  \notag \\
&&+\frac{\kappa \cos \theta }{1-\kappa r\cos \theta }v_{s}\omega _{s}\mathbf{%
e}_{r}-\frac{\kappa \sin \theta }{1-\kappa r\cos \theta }v_{s}\omega _{s}%
\mathbf{e}_{\varphi }.  \label{f8}
\end{eqnarray}

We compute now the viscous terms%
\begin{eqnarray*}
&&\Delta \widetilde{\boldsymbol{\omega }}_{0} \\
&=&\frac{1}{h_{1}}\frac{\partial }{\partial x_{1}}\left[ \frac{1}{%
h_{1}h_{2}h_{3}}\left[ \frac{\partial }{\partial x_{1}}(h_{2}h_{3}\omega
_{1})+\frac{\partial }{\partial x_{2}}(h_{3}h_{1}\omega _{2})+\frac{\partial
}{\partial x_{3}}(h_{1}h_{2}\omega _{3})\right] \right] \boldsymbol{e}_{1} \\
&&-\frac{1}{h_{2}h_{3}}\left[ \frac{\partial }{\partial x_{2}}\left( \frac{%
h_{3}}{h_{1}h_{2}}\left[ \frac{\partial (h_{2}\omega _{2})}{\partial x_{1}}-%
\frac{\partial (h_{1}\omega _{1})}{\partial x_{2}}\right] \right) -\frac{%
\partial }{\partial x_{3}}\left( \frac{h_{2}}{h_{1}h_{3}}\left[ \frac{%
\partial (h_{1}\omega _{1})}{\partial x_{3}}-\frac{\partial (h_{3}\omega
_{3})}{\partial x_{1}}\right] \right) \right] \boldsymbol{e}_{1} \\
&&+\frac{1}{h_{2}}\frac{\partial }{\partial x_{2}}\left[ \frac{1}{%
h_{1}h_{2}h_{3}}\left[ \frac{\partial }{\partial x_{1}}(h_{2}h_{3}\omega
_{1})+\frac{\partial }{\partial x_{2}}(h_{3}h_{1}\omega _{2})+\frac{\partial
}{\partial x_{3}}(h_{1}h_{2}\omega _{3})\right] \right] \boldsymbol{e}_{2} \\
&&-\frac{1}{h_{1}h_{3}}\left[ \frac{\partial }{\partial x_{3}}\left( \frac{%
h_{1}}{h_{2}h_{3}}\left[ \frac{\partial (h_{3}\omega _{3})}{\partial x_{2}}-%
\frac{\partial (h_{2}\omega _{2})}{\partial x_{3}}\right] \right) -\frac{%
\partial }{\partial x_{1}}\left( \frac{h_{3}}{h_{1}h_{2}}\left[ \frac{%
\partial (h_{2}\omega _{2})}{\partial x_{1}}-\frac{\partial (h_{1}\omega
_{1})}{\partial x_{2}}\right] \right) \right] \boldsymbol{e}_{2} \\
&&+\frac{1}{h_{3}}\frac{\partial }{\partial x_{3}}\left[ \frac{1}{%
h_{1}h_{2}h_{3}}\left[ \frac{\partial }{\partial x_{1}}(h_{2}h_{3}\omega
_{1})+\frac{\partial }{\partial x_{2}}(h_{3}h_{1}\omega _{2})+\frac{\partial
}{\partial x_{3}}(h_{1}h_{2}\omega _{3})\right] \right] \boldsymbol{e}_{3} \\
&&-\frac{1}{h_{1}h_{2}}\left[ \frac{\partial }{\partial x_{1}}\left( \frac{%
h_{2}}{h_{1}h_{3}}\left[ \frac{\partial (h_{1}\omega _{1})}{\partial x_{2}}-%
\frac{\partial (h_{3}\omega _{3})}{\partial x_{1}}\right] \right) -\frac{%
\partial }{\partial x_{2}}\left( \frac{h_{1}}{h_{2}h_{3}}\left[ \frac{%
\partial (h_{3}\omega _{3})}{\partial x_{2}}-\frac{\partial (h_{2}\omega
_{2})}{\partial x_{3}}\right] \right) \right] \boldsymbol{e}_{3},
\end{eqnarray*}%
and compute the nonzero terms in the $\boldsymbol{e}_{1}=\mathbf{e}_{r}$
component as
\begin{eqnarray*}
B_{r} &\equiv &\frac{1}{h_{1}}\frac{\partial }{\partial x_{1}}\left( \frac{1%
}{h_{1}h_{2}h_{3}}\frac{\partial }{\partial x_{3}}(h_{1}h_{2}\omega
_{3})\right) -\frac{1}{h_{2}h_{3}}\frac{\partial }{\partial x_{3}}\left(
\frac{h_{2}}{h_{1}h_{3}}\frac{\partial }{\partial x_{1}}(h_{3}\omega
_{3})\right) \\
&=&\frac{\partial }{\partial r}\left( \frac{1}{(1-\kappa r\cos \theta )}%
\frac{\partial }{\partial s}\omega _{s}\right) -\frac{1}{(1-\kappa r\cos
\theta )}\frac{\partial }{\partial s}\left( \frac{1}{(1-\kappa r\cos \theta )%
}\frac{\partial }{\partial r}((1-\kappa r\cos \theta )\omega _{s})\right) \\
&=&\frac{\kappa \cos \theta }{(1-\kappa r\cos \theta )^{2}}\frac{\partial }{%
\partial s}\omega _{s}+\frac{1}{(1-\kappa r\cos \theta )}\frac{\partial }{%
\partial s}\left( \frac{\kappa \cos \theta }{(1-\kappa r\cos \theta )}\omega
_{s}\right)
\end{eqnarray*}%
as well as the nonzero terms in the $\boldsymbol{e}_{2}=\mathbf{e}_{\varphi
} $ component as
\begin{eqnarray*}
B_{\varphi } &\equiv &\frac{1}{h_{2}}\frac{\partial }{\partial x_{2}}\left(
\frac{1}{h_{1}h_{2}h_{3}}\frac{\partial }{\partial x_{3}}(h_{1}h_{2}\omega
_{3})\right) -\frac{1}{h_{1}h_{3}}\frac{\partial }{\partial x_{3}}\left(
\frac{h_{1}}{h_{2}h_{3}}\frac{\partial }{\partial x_{2}}(h_{3}\omega
_{3})\right) \\
&=&\frac{1}{r}\frac{\partial }{\partial \varphi }\left( \frac{1}{(1-\kappa
r\cos \theta )}\frac{\partial }{\partial s}\omega _{s}\right) -\frac{1}{%
r(1-\kappa r\cos \theta )}\frac{\partial }{\partial s}\left( \frac{1}{%
(1-\kappa r\cos \theta )}\frac{\partial }{\partial \varphi }((1-\kappa r\cos
\theta )\omega _{s})\right) \\
&=&-\frac{\kappa \sin \theta }{(1-\kappa r\cos \theta )^{2}}\frac{\partial }{%
\partial s}\omega _{s}-\frac{1}{(1-\kappa r\cos \theta )}\frac{\partial }{%
\partial s}\left( \frac{\kappa \sin \theta }{(1-\kappa r\cos \theta )}\omega
_{s}\right) .
\end{eqnarray*}%
Note that the nonzero terms in the $\boldsymbol{e}_{3}=\mathbf{e}_{s}$
component are
\begin{eqnarray}
&&\frac{1}{h_{3}}\frac{\partial }{\partial x_{3}}\left( \frac{1}{%
h_{1}h_{2}h_{3}}\frac{\partial }{\partial x_{3}}(h_{1}h_{2}\omega
_{3})\right) +\frac{1}{h_{1}h_{2}}\frac{\partial }{\partial x_{2}}\left(
\frac{h_{1}}{h_{2}h_{3}}\frac{\partial }{\partial x_{2}}(h_{3}\omega
_{3})\right)  \notag \\
&&+\frac{1}{h_{1}h_{2}}\frac{\partial }{\partial x_{1}}\left( \frac{h_{2}}{%
h_{1}h_{3}}\frac{\partial }{\partial x_{1}}(h_{3}\omega _{3})\right)  \notag
\\
&=&\frac{1}{(1-\kappa r\cos \theta )}\frac{\partial }{\partial s}\left(
\frac{1}{(1-\kappa r\cos \theta )}\frac{\partial }{\partial s}\omega
_{s}\right) +\frac{1}{r^{2}}\frac{\partial }{\partial \varphi }\left( \frac{1%
}{(1-\kappa r\cos \theta )}\frac{\partial }{\partial \varphi }((1-\kappa
r\cos \theta )\omega _{s})\right)  \notag \\
&&+\frac{1}{r}\frac{\partial }{\partial r}\left( \frac{r}{(1-\kappa r\cos
\theta )}\frac{\partial }{\partial r}((1-\kappa r\cos \theta )\omega
_{s})\right)  \notag \\
&=&\frac{1}{(1-\kappa r\cos \theta )}\frac{\partial }{\partial s}\left(
\frac{1}{(1-\kappa r\cos \theta )}\frac{\partial }{\partial s}\omega
_{s}\right) +\frac{1}{r^{2}}\frac{\partial ^{2}\omega _{s}}{\partial \varphi
^{2}}+\frac{1}{r}\frac{\partial }{\partial r}\left( r\frac{\partial \omega
_{s}}{\partial r}\right)  \notag \\
&&+\frac{\kappa \sin \theta }{r}\frac{\partial }{\partial \varphi }\left(
\frac{1}{(1-\kappa r\cos \theta )}\omega _{s}\right) -\kappa \cos \theta
\frac{\partial }{\partial r}\left( \frac{1}{(1-\kappa r\cos \theta )}\omega
_{s}\right) .  \label{f9}
\end{eqnarray}%
From the vorticity field in (\ref{gg1}) we can easily estimate
\begin{equation*}
\left\vert B_{r}\right\vert +\left\vert B_{\varphi }\right\vert \leq \frac{%
C\Gamma ^{2}/\nu }{(\nu t)}(1+\rho ^{4})e^{-\frac{\rho ^{2}}{4}}
\end{equation*}%
and similar estimates for the terms (\ref{f1}), (\ref{f3}), (\ref{f6}), (\ref%
{f7}) and the last two at the right hand side of (\ref{f8}), that we can
group together into a contribution $B$ such that%
\begin{equation*}
\left\vert B\right\vert \leq \frac{C\Gamma ^{2}}{(\nu t)}(1+\rho ^{4})e^{-%
\frac{\rho ^{2}}{4}}.
\end{equation*}%
As for the sum of the stretching, convective and viscous terms in the $%
\mathbf{e}_{s}$ direction, we have%
\begin{eqnarray*}
&&\left( v_{r}\frac{\partial \omega _{s}}{\partial r}+\frac{v_{\varphi }}{r}%
\frac{\partial \omega _{s}}{\partial \varphi }+\frac{v_{s}}{1-\kappa r\cos
\theta }\frac{\partial \omega _{s}}{\partial s}\right) -\left( \frac{\omega
_{s}}{1-\kappa r\cos \theta }\frac{\partial v_{s}}{\partial s}\right) \\
&&+\frac{\kappa \cos \theta }{1-\kappa r\cos \theta }\omega _{s}v_{r}-\frac{%
\kappa \sin \theta }{1-\kappa r\cos \theta }\omega _{s}v_{\varphi } \\
&&-\nu \left( \frac{1}{(1-\kappa r\cos \theta )}\frac{\partial }{\partial s}%
\left( \frac{1}{(1-\kappa r\cos \theta )}\frac{\partial }{\partial s}\omega
_{s}\right) +\frac{1}{r^{2}}\frac{\partial ^{2}\omega _{s}}{\partial \varphi
^{2}}+\frac{1}{r}\frac{\partial }{\partial r}\left( r\frac{\partial \omega
_{s}}{\partial r}\right) \right. \\
&&\left. +\frac{\kappa \sin \theta }{r}\frac{\partial }{\partial \varphi }%
\left( \frac{1}{(1-\kappa r\cos \theta )}\omega _{s}\right) -\kappa \cos
\theta \frac{\partial }{\partial r}\left( \frac{1}{(1-\kappa r\cos \theta )}%
\omega _{s}\right) \right) \\
&=&v_{r}\frac{\partial \omega _{s}}{\partial r}+\frac{v_{\varphi }}{r}\frac{%
\partial \omega _{s}}{\partial \varphi }-\kappa \sin \theta \omega
_{s}v_{\varphi } \\
&&-\nu \left( \frac{1}{r^{2}}\frac{\partial ^{2}\omega _{s}}{\partial
\varphi ^{2}}+\frac{1}{r}\frac{\partial }{\partial r}\left( r\frac{\partial
\omega _{s}}{\partial r}\right) -\kappa \cos \theta \frac{\partial \omega
_{s}}{\partial r}\right) +~~B_{s},
\end{eqnarray*}%
with%
\begin{equation*}
\left\vert B_{s}\right\vert \leq \frac{C\Gamma ^{2}}{(\nu t)}(1+\rho
^{4})e^{-\frac{\rho ^{2}}{4}}.
\end{equation*}

Concerning $\frac{\partial \widetilde{\boldsymbol{\omega }}}{\partial t}$,
we must take into account that both the position vector and the vorticity
itself are written in a coordinate frame that evolves in time. Since $%
\widetilde{\boldsymbol{\omega }}=\omega _{i}\boldsymbol{e}_{i}$ (using
Einstein's summation convention) we have
\begin{equation*}
\frac{\partial \widetilde{\boldsymbol{\omega }}_{0}}{\partial t}=\frac{%
\partial \omega _{i}}{\partial t}\boldsymbol{e}_{i}+\omega _{i}\frac{d%
\boldsymbol{e}_{i}}{dt}\equiv \frac{\partial \omega _{i}}{\partial t}%
\boldsymbol{e}_{i}+\mathbf{C},
\end{equation*}%
but we can easily estimate%
\begin{equation*}
\left\vert \mathbf{C}\right\vert =\left\vert \omega _{i}\frac{d\boldsymbol{e}%
_{i}}{dt}\right\vert \leq \frac{C\Gamma ^{2}\left\vert \log (\nu
t)\right\vert }{(\nu t)}(1+\rho ^{4})e^{-\frac{\rho ^{2}}{4}}.
\end{equation*}

The only terms in Navier-Stokes system that are left to estimate are then
collected in
\begin{eqnarray*}
&&\frac{\partial \omega _{s}}{\partial t}+v_{r}\frac{\partial \omega _{s}}{%
\partial r}+\frac{v_{\varphi }}{r}\frac{\partial \omega _{s}}{\partial
\varphi }-\kappa \sin \theta \omega _{s}v_{\theta } \\
&&-\nu \left( \frac{1}{r^{2}}\frac{\partial ^{2}\omega _{s}}{\partial
\varphi ^{2}}+\frac{1}{r}\frac{\partial }{\partial r}\left( r\frac{\partial
\omega _{s}}{\partial r}\right) +\frac{\kappa \sin \theta }{r}\frac{\partial
\omega _{s}}{\partial \varphi }-\kappa \cos \theta \frac{\partial \omega _{s}%
}{\partial r}\right) \\
&\equiv &A\left[ \omega _{s},v_{r},v_{\varphi }\right] ,
\end{eqnarray*}%
which corresponds exactly to the definition of the nonlinear operator $A$
defined in (\ref{defa}). By the procedure of construction of $\omega _{0}$,
we have (see (\ref{estia})):%
\begin{equation*}
\left\vert A\left[ \omega _{s},v_{r},v_{\varphi }\right] \right\vert \leq
\frac{C(\Gamma ^{2}+\Gamma \nu )}{(\nu t)}(1+\rho ^{4})e^{-\frac{\rho ^{2}}{4%
}}.
\end{equation*}

Finally, if%
\begin{equation*}
\widetilde{\boldsymbol{\omega }}_{0}=\eta _{R}\omega _{0}\mathbf{e}_{s},
\end{equation*}%
we have%
\begin{eqnarray*}
&&\widetilde{\boldsymbol{\omega }}_{0,t}+\widetilde{\mathbf{v}}\cdot \nabla
\widetilde{\boldsymbol{\omega }}_{0}-\widetilde{\boldsymbol{\omega }}%
_{0}\cdot \nabla \widetilde{\mathbf{v}}_{0}-\nu \Delta \widetilde{%
\boldsymbol{\omega }}_{0} \\
&&\eta _{R}\left[ (\omega _{0}\mathbf{e}_{s})_{t}+\widetilde{\mathbf{v}}%
\cdot \nabla (\omega _{0}\mathbf{e}_{s})-\left( \omega _{0}\mathbf{e}%
_{s}\right) \cdot \nabla \widetilde{\mathbf{v}}-\nu \Delta (\omega _{0}%
\mathbf{e}_{s})\right] +\eta _{R,t}(\omega _{0}\mathbf{e}_{s}) \\
&&+\left( \widetilde{\mathbf{v}}\cdot \nabla \eta _{R}\right) (\omega _{0}%
\mathbf{e}_{s})-2\nu \left( \nabla \eta _{R}\right) \cdot \nabla (\omega _{0}%
\mathbf{e}_{s})-(\omega _{0}\mathbf{e}_{s})\Delta \eta _{R},
\end{eqnarray*}%
but then%
\begin{eqnarray*}
&&\left\vert \eta _{R,t}(\omega _{0}\mathbf{e}_{s})+\left( \widetilde{%
\mathbf{v}}_{0}\cdot \nabla \eta _{R}\right) (\omega _{0}\mathbf{e}%
_{s})-2\nu \left( \nabla \eta _{R}\right) \cdot \nabla (\omega _{0}\mathbf{e}%
_{s})-\nu (\omega _{0}\mathbf{e}_{s})\Delta \eta _{R}\right\vert \\
&\leq &\frac{C(\Gamma ^{2}+\Gamma \nu )\left\vert \log (\nu t)\right\vert }{%
(\nu t)}(1+\rho ^{4})e^{-\frac{\rho ^{2}}{4}}.
\end{eqnarray*}%
The estimates extend to other small $s$-intervals that cover the whole
vortex filament.

As a corollary of the previous Lemmas, we have the following

\begin{corollary}
\label{lemaf}%
\begin{equation}
\left\Vert \mathbf{F}(\mathbf{x},t)\right\Vert _{L^{2}}\leq C(\Gamma
^{2}+\Gamma \nu )\frac{\left\vert \log (\nu t)\right\vert }{(\nu t)^{\frac{1%
}{2}}}.  \label{cotaf}
\end{equation}
\end{corollary}

Finally, we will need estimates for the function $\mathbf{G}$ resulting from
the application of Biot-Savart law to $\mathbf{F}$, that is%
\begin{equation}
\mathbf{G}(\mathbf{x},t)=\frac{1}{4\pi }\int \frac{\mathbf{F}(\mathbf{x}%
^{\prime },t)\times (\mathbf{x}-\mathbf{x}^{\prime })}{\left\vert \mathbf{x}-%
\mathbf{x}^{\prime }\right\vert ^{3}}d\mathbf{x}^{\prime }.  \label{bs}
\end{equation}%
We have the following Lemma:

\begin{lemma}
\label{lemag}The function $\mathbf{G}(\mathbf{x},t)$ defined by (\ref{bs})
satisfies%
\begin{equation}
\left\Vert \mathbf{G}(\mathbf{x},t)\right\Vert _{L^{2}}\leq C(\Gamma
^{2}+\Gamma \nu )\left\vert \log (\nu t)\right\vert ^{\frac{3}{2}}.
\label{cotag}
\end{equation}
\end{lemma}

\textbf{Proof.} The function $\mathbf{G}(\mathbf{x},t)$ is the velocity
field associated to the vorticity field $\mathbf{F}(\mathbf{x},t)$. In
regions far from the vortex filament, i.e. for $\left\vert \mathbf{x}%
\right\vert $ sufficiently large, we have%
\begin{equation*}
\mathbf{G}(\mathbf{x},t)\sim \frac{1}{4\pi \left\vert \mathbf{x}\right\vert
^{2}}\int \mathbf{F}(\mathbf{x}^{\prime },t)\times \frac{(\mathbf{x}-\mathbf{%
x}^{\prime })}{\left\vert \mathbf{x}-\mathbf{x}^{\prime }\right\vert }d%
\mathbf{x}^{\prime },
\end{equation*}%
and hence%
\begin{equation*}
\left\vert \mathbf{G}(\mathbf{x},t)\right\vert \leq \frac{C}{4\pi \left\vert
\mathbf{x}\right\vert ^{2}}\int \left\vert \mathbf{F}(\mathbf{x}^{\prime
},t)\right\vert d\mathbf{x}^{\prime }\leq \frac{C(\Gamma ^{2}+\Gamma \nu
)\left\vert \log (\nu t)\right\vert }{4\pi \left\vert \mathbf{x}\right\vert
^{2}}.
\end{equation*}%
On the other hand, in regions close to the filament the vorticity $\mathbf{F}%
(\mathbf{x},t)$ generates an angular velocity field $\mathbf{G}_{\theta }(%
\mathbf{x},t)$ which is $O(\left\vert \log (\nu t)\right\vert r^{-1})$ for $%
(\nu t)^{\frac{1}{2}}\leq r\ll 1$. This yields an $L^{2}$ norm of $\mathbf{G}%
_{\theta }(\mathbf{x},t)$ which is%
\begin{equation*}
\left\Vert \mathbf{G}_{\theta }(\mathbf{x},t)\right\Vert _{L^{2}}^{2}\leq
C\int_{(\nu t)^{\frac{1}{2}}}^{1}\left\vert \log (\nu t)\right\vert
^{2}r^{-2}rdr=O(\left\vert \log (\nu t)\right\vert ^{3}).
\end{equation*}%
As for the region $r<(\nu t)^{\frac{1}{2}}$, the function $\mathbf{G}(%
\mathbf{x},t)$ is $O(\left\vert \log (\nu t)\right\vert /(\nu t)^{\frac{1}{2}%
})$ while the total volume of the region is $O(\nu t)$ which implies that
its contribution to the $L^{2}$ norm is merely $O(\left\vert \log (\nu
t)\right\vert )$.

\section{Construction of the full solution}

\bigskip Navier-Stokes equation for the velocity of an incompressible fluid
can be written as%
\begin{equation*}
\mathbf{v}_{t}-\mathbf{v}\times \boldsymbol{\omega +}\frac{\nabla \left\vert
\mathbf{v}\right\vert ^{2}}{2}=-\nabla p+\nu \Delta \mathbf{v},
\end{equation*}%
while the vorticity equation as%
\begin{equation}
\boldsymbol{\omega }_{t}-\nabla \times \left( \mathbf{v}\times \boldsymbol{%
\omega }\right) =\nu \Delta \boldsymbol{\omega }.  \label{edp1}
\end{equation}%
Under the change $t=t^{\prime }/\nu $ and $\boldsymbol{\omega =}\Gamma
\boldsymbol{\omega }^{\prime }$,$\mathbf{v}\boldsymbol{=}\Gamma \mathbf{v}%
^{\prime }$ the equations for the vorticity transforms (omitting primes for
the sake of simplicity) into%
\begin{equation*}
\boldsymbol{\omega }_{t}-\frac{\Gamma }{\nu }\nabla \times \left( \mathbf{v}%
\times \boldsymbol{\omega }\right) =\Delta \boldsymbol{\omega }.
\end{equation*}

By writing%
\begin{eqnarray*}
\mathbf{v} &=&\mathbf{v}_{0}+\widetilde{\mathbf{v}}, \\
\boldsymbol{\omega } &=&\boldsymbol{\omega }_{0}+\widetilde{\boldsymbol{%
\omega }},
\end{eqnarray*}%
where $\boldsymbol{\omega }_{0}$ is the approximate vorticity constructed in
previous sections and $\mathbf{v}_{0}$ the corresponding velocity, i.e.
\begin{equation*}
\nabla \times \mathbf{v}_{0}=\boldsymbol{\omega }_{0}.
\end{equation*}%
The approximate vorticity $\boldsymbol{\omega }_{0}$ satisfies
\begin{equation*}
\boldsymbol{\omega }_{0,t}-\frac{\Gamma }{\nu }\nabla \times \left( \mathbf{v%
}_{0}\times \boldsymbol{\omega }_{0}\right) =\Delta \boldsymbol{\omega }_{0}-%
\mathbf{F},
\end{equation*}%
where, according to Corollary \ref{lemaf},
\begin{equation}
\left\Vert \mathbf{F}\right\Vert _{L^{2}}\leq C(\Gamma /\nu +1)\frac{%
\left\vert \log t\right\vert }{t^{\frac{1}{2}}}.  \label{af0}
\end{equation}

Hence, $\widetilde{\boldsymbol{\omega }}$ satisfies

\begin{equation}
\widetilde{\boldsymbol{\omega }}_{t}-\Delta \widetilde{\boldsymbol{\omega }}=%
\mathbf{F}+\frac{\Gamma }{\nu }\nabla \times \left( \mathbf{v}_{0}\times
\widetilde{\boldsymbol{\omega }}+\widetilde{\mathbf{v}}\times \boldsymbol{%
\omega }_{0}\right) +\frac{\Gamma }{\nu }\nabla \times \left( \widetilde{%
\mathbf{v}}\times \widetilde{\boldsymbol{\omega }}\right) .  \label{ecuaw}
\end{equation}

\begin{theorem}
\bigskip Let $\mathbf{F}$ satisfy (\ref{af0}). Then there exists a unique
solution $\widetilde{\boldsymbol{\omega }}$ to (\ref{ecuaw}) such that $%
\left\Vert \widetilde{\boldsymbol{\omega }}\right\Vert _{L^{2}}(t=0)=0$ as
well as satisfying%
\begin{equation}
\sup_{\left[ 0,T\right] }\left\Vert \widetilde{\boldsymbol{\omega }}%
\right\Vert _{L^{2}}^{2}+\int_{0}^{T}\left\Vert \nabla \widetilde{%
\boldsymbol{\omega }}\right\Vert _{L^{2}}^{2}dt^{\prime }\leq CT\left\vert
\log (T)\right\vert ^{2},  \label{est0}
\end{equation}%
for $\frac{\Gamma }{\nu }$ sufficiently small. $C$ a suitable constant that
only depends on initial data and $\nu T$ is sufficiently small.
\end{theorem}

\textbf{Proof. }We consider the modified linear equation
\begin{equation*}
\widetilde{\boldsymbol{\omega }}_{t}-\Delta \widetilde{\boldsymbol{\omega }}=%
\mathbf{F}+\frac{\Gamma }{\nu }\nabla \times \left( \mathbf{v}_{0}\times
\overline{\boldsymbol{\omega }}+\overline{\mathbf{v}}\times \boldsymbol{%
\omega }_{0}\right) +\frac{\Gamma }{\nu }\nabla \times \left( \overline{%
\mathbf{v}}\times \overline{\boldsymbol{\omega }}\right) ,
\end{equation*}%
multiply by $\widetilde{\boldsymbol{\omega }}$ and integrate by parts to
obtain%
\begin{eqnarray*}
&&\frac{1}{2}\frac{d}{dt}\left\Vert \widetilde{\boldsymbol{\omega }}%
\right\Vert ^{2}+\left\Vert \nabla \widetilde{\boldsymbol{\omega }}%
\right\Vert ^{2} \\
&=&\int \mathbf{F}\cdot \widetilde{\boldsymbol{\omega }}-\frac{\Gamma }{\nu }%
\int \left( \nabla \times \widetilde{\boldsymbol{\omega }}\right) \cdot
\left( \mathbf{v}_{0}\times \overline{\boldsymbol{\omega }}+\widetilde{%
\mathbf{v}}\times \boldsymbol{\omega }_{0}+\overline{\mathbf{v}}\times
\overline{\boldsymbol{\omega }}\right) \\
&\leq &\left\Vert \mathbf{F}\right\Vert _{L^{2}}\left\Vert \widetilde{%
\boldsymbol{\omega }}\right\Vert _{L^{2}} \\
&&+C\frac{\Gamma }{\nu }\left\Vert \nabla \widetilde{\boldsymbol{\omega }}%
\right\Vert _{L^{2}}\left( \left\Vert \mathbf{v}_{0}\times \overline{%
\boldsymbol{\omega }}\right\Vert _{L^{2}}+\left\Vert \overline{\mathbf{v}}%
\times \boldsymbol{\omega }_{0}\right\Vert _{L^{2}}+\left\Vert \overline{%
\mathbf{v}}\times \overline{\boldsymbol{\omega }}\right\Vert _{L^{2}}\right)
,
\end{eqnarray*}%
and hence%
\begin{eqnarray*}
&&\frac{1}{2}\left\Vert \widetilde{\boldsymbol{\omega }}\right\Vert ^{2}(t)+%
\frac{1}{2}\int_{0}^{t}\left\Vert \nabla \widetilde{\boldsymbol{\omega }}%
\right\Vert ^{2}(\tau )d\tau \\
&\leq &\frac{1}{4}\sup_{t}\left\Vert \widetilde{\boldsymbol{\omega }}%
\right\Vert _{L^{2}}^{2}+\left( \int_{0}^{t}\left\Vert \mathbf{F}\right\Vert
_{L^{2}}d\tau \right) ^{2} \\
&&+\frac{1}{2}\left( \frac{C\Gamma }{\nu }\right) ^{2}\int_{0}^{t}\left(
\left\Vert \mathbf{v}_{0}\times \overline{\boldsymbol{\omega }}\right\Vert
_{L^{2}}^{2}+\left\Vert \overline{\mathbf{v}}\times \boldsymbol{\omega }%
_{0}\right\Vert _{L^{2}}^{2}+\left\Vert \overline{\mathbf{v}}\times
\overline{\boldsymbol{\omega }}\right\Vert _{L^{2}}^{2}\right) d\tau \\
&\equiv &\frac{1}{4}\sup_{t}\left\Vert \widetilde{\boldsymbol{\omega }}%
\right\Vert _{L^{2}}^{2}+R(t).
\end{eqnarray*}%
Therefore, we have%
\begin{eqnarray}
\frac{1}{2}\sup_{t}\left\Vert \widetilde{\boldsymbol{\omega }}\right\Vert
^{2} &\leq &\frac{1}{4}\sup_{t}\left\Vert \widetilde{\boldsymbol{\omega }}%
\right\Vert _{L^{2}}^{2}+R(t),  \label{af1} \\
\frac{1}{2}\int_{0}^{t}\left\Vert \nabla \widetilde{\boldsymbol{\omega }}%
\right\Vert ^{2}(\tau )d\tau &\leq &\frac{1}{4}\sup_{t}\left\Vert \widetilde{%
\boldsymbol{\omega }}\right\Vert _{L^{2}}^{2}+R(t),  \label{af2}
\end{eqnarray}%
and adding to (\ref{af1}), inequality (\ref{af2}) multiplied by $\alpha $ we
find
\begin{equation*}
\frac{1-\alpha }{4}\sup_{t}\left\Vert \widetilde{\boldsymbol{\omega }}%
\right\Vert ^{2}+\frac{\alpha }{2}\int_{0}^{t}\left\Vert \nabla \widetilde{%
\boldsymbol{\omega }}\right\Vert ^{2}(\tau )d\tau \leq (1+\alpha )R(t),
\end{equation*}%
and choosing $\alpha =1/3$,%
\begin{equation*}
\sup_{t}\left\Vert \widetilde{\boldsymbol{\omega }}\right\Vert
^{2}+\int_{0}^{t}\left\Vert \nabla \widetilde{\boldsymbol{\omega }}%
\right\Vert ^{2}(\tau )d\tau \leq 8R(t).
\end{equation*}%
Now we estimate the terms involved in $R(t)$. First, by (\ref{af0}),%
\begin{equation*}
\int_{0}^{t}\left\Vert \mathbf{F}\right\Vert _{L^{2}}d\tau \leq Ct^{\frac{1}{%
2}}\left\vert \log t\right\vert .
\end{equation*}%
Secondly,%
\begin{equation}
\left\Vert \mathbf{v}_{0}\times \overline{\boldsymbol{\omega }}\right\Vert
_{L^{2}}^{2}\leq \frac{C}{t}\left\Vert \overline{\boldsymbol{\omega }}%
\right\Vert _{L^{2}}^{2}.  \label{ass1}
\end{equation}%
Third, since by Fefferman-Phong inequality (cf. \cite{Fe}), we have%
\begin{equation*}
\left( \int \left\vert \overline{\mathbf{v}}\right\vert ^{2}\left\vert
\boldsymbol{\omega }_{0}\right\vert ^{2}\right) ^{\frac{1}{2}}\leq
C\sup_{B_{r}}\left( r\left( \frac{1}{\left\vert B_{r}\right\vert }%
\int_{B_{r}}\left\vert \boldsymbol{\omega }_{0}\right\vert ^{2p}\right) ^{%
\frac{1}{2p}}\right) \left( \int \left\vert \nabla \overline{\mathbf{v}}%
\right\vert ^{2}\right) ^{\frac{1}{2}},\ 1<p<n/2,
\end{equation*}%
and estimating%
\begin{eqnarray*}
\sup_{B_{r}}\left( r\left( \frac{1}{\left\vert B_{r}\right\vert }%
\int_{B_{r}}\left\vert \boldsymbol{\omega }_{0}\right\vert ^{2p}\right) ^{%
\frac{1}{2p}}\right) &\leq &\frac{C}{t^{\frac{1}{2}}}, \\
\left( \int \left\vert \nabla \overline{\mathbf{v}}\right\vert ^{2}\right) ^{%
\frac{1}{2}} &\leq &C\left\Vert \overline{\boldsymbol{\omega }}\right\Vert
_{L^{2}},
\end{eqnarray*}%
we conclude%
\begin{equation}
\left\Vert \overline{\mathbf{v}}\times \boldsymbol{\omega }_{0}\right\Vert
_{L^{2}}^{2}\leq \frac{C}{t}\left\Vert \overline{\boldsymbol{\omega }}%
\right\Vert _{L^{2}}^{2}.  \label{ass2}
\end{equation}%
Finally%
\begin{equation*}
\left\Vert \overline{\mathbf{v}}\times \overline{\boldsymbol{\omega }}%
\right\Vert _{L^{2}}^{2}\leq \left\Vert \overline{\mathbf{v}}\right\Vert
_{L^{6}}^{2}\left\Vert \overline{\boldsymbol{\omega }}\right\Vert
_{L^{3}}^{2},
\end{equation*}%
and using Sobolev and interpolation inequalities%
\begin{eqnarray*}
\left\Vert \overline{\mathbf{v}}\right\Vert _{L^{6}} &\leq &C\left\Vert
\overline{\boldsymbol{\omega }}\right\Vert _{L^{3}}, \\
\left\Vert \overline{\boldsymbol{\omega }}\right\Vert _{L^{3}} &\leq
&C\left\Vert \overline{\boldsymbol{\omega }}\right\Vert _{L^{2}}^{\frac{1}{2}%
}\left\Vert \nabla \overline{\boldsymbol{\omega }}\right\Vert _{L^{2}}^{%
\frac{1}{2}},
\end{eqnarray*}%
we conclude%
\begin{equation}
\left\Vert \overline{\mathbf{v}}\times \overline{\boldsymbol{\omega }}%
\right\Vert _{L^{2}}^{2}\leq C\left\Vert \overline{\boldsymbol{\omega }}%
\right\Vert _{L^{2}}^{3}\left\Vert \nabla \overline{\boldsymbol{\omega }}%
\right\Vert _{L^{2}}.  \label{ass3}
\end{equation}%
Using (\ref{ass1}), (\ref{ass2}), (\ref{ass3}) we find%
\begin{eqnarray}
&&\sup_{t}\left\Vert \widetilde{\boldsymbol{\omega }}\right\Vert
^{2}(t)+\int_{0}^{t}\left\Vert \nabla \widetilde{\boldsymbol{\omega }}%
\right\Vert ^{2}(\tau )d\tau  \notag \\
&\leq &8\left( \int_{0}^{t}\left\Vert \mathbf{F}\right\Vert _{L^{2}}d\tau
\right) ^{2}+C\left( \frac{\Gamma }{\nu }\right) ^{2}\int_{0}^{t}\frac{%
\left\Vert \overline{\boldsymbol{\omega }}\right\Vert ^{2}(\tau )}{\tau }%
d\tau +C\left( \frac{\Gamma }{\nu }\right) ^{2}\int_{0}^{t}\left\Vert
\overline{\boldsymbol{\omega }}\right\Vert _{L^{2}}^{3}\left\Vert \nabla
\overline{\boldsymbol{\omega }}\right\Vert _{L^{2}}d\tau .  \label{ineq}
\end{eqnarray}%
We define now the norm%
\begin{equation*}
\left\Vert \widetilde{\boldsymbol{\omega }}\right\Vert _{X_{T}}^{2}=\sup_{%
\left[ 0,T\right] }\left( \frac{1}{t\left\vert \log t\right\vert ^{2}}%
\left\Vert \widetilde{\boldsymbol{\omega }}\right\Vert _{L^{2}}^{2}+\frac{1}{%
t\left\vert \log t\right\vert ^{2}}\int_{0}^{t}\left\Vert \nabla \widetilde{%
\boldsymbol{\omega }}\right\Vert _{L^{2}}^{2}dt^{\prime }\right) ,
\end{equation*}%
and the mapping $\mathbb{T}$ that assigns to $\overline{\boldsymbol{\omega }}
$ in a ball of radius $R$ in $X_{T}$ the function $\widetilde{\boldsymbol{%
\omega }}$ satisfying (\ref{ecuaw}). By inequality (\ref{ineq}) we have
\begin{equation}
\left\Vert \widetilde{\boldsymbol{\omega }}\right\Vert _{X_{T}}^{2}\leq
\frac{1}{2}R^{2}+C\left( \frac{\Gamma }{\nu }\right) ^{2}\left( \left\Vert
\overline{\boldsymbol{\omega }}\right\Vert _{X_{T}}^{2}+T^{\frac{1}{2}%
}\left\Vert \overline{\boldsymbol{\omega }}\right\Vert _{X_{T}}^{4}\right) ,
\label{estiwt}
\end{equation}%
where%
\begin{equation*}
\frac{1}{2}R^{2}=\sup_{\left[ 0,T\right] }\frac{8}{t\left\vert \log
t\right\vert ^{2}}\left( \int_{0}^{t}\left\Vert \mathbf{F}\right\Vert
_{L^{2}}d\tau \right) ^{2}.
\end{equation*}%
Hence,if $\frac{\Gamma }{\nu }$ is sufficiently small, the mapping $\mathbb{T%
}$ maps the ball of radius $R$ into itself. Banach's fixed point theorem
implies then existence of solutions to (\ref{ecuaw}). It is straightforward
to see that the mapping is also a contraction and hence the solution is
unique. This ends the proof of the theorem.

Note that, in terms of the variables and unknowns in equation (\ref{edp1}),
the estimate (\ref{est0}) reads

\begin{equation*}
\sup_{\left[ 0,T\right] }\left\Vert \widetilde{\boldsymbol{\omega }}%
\right\Vert _{L^{2}}^{2}+\nu \int_{0}^{T}\left\Vert \nabla \widetilde{%
\boldsymbol{\omega }}\right\Vert _{L^{2}}^{2}dt^{\prime }\leq C\Gamma
^{2}(\nu T)\left\vert \log (\nu T)\right\vert ^{2}.
\end{equation*}

\begin{remark}
We discuss now the time-span at which our construction is a valid
approximations or, equivalently, for how long a vortex filament moving with
the binormal flow is a good approximation to Navier-Stokes system. On one
hand, our solutions are valid provided $\frac{\Gamma }{\nu }$ is
sufficiently small, that is $\Gamma =O(\nu )$. On the other hand, the
thickness of the vortex filament is $O((\nu t)^{\frac{1}{2}})$, and it needs
to be small compared with the smallest radius of curvature of the filament
(that we assume to be $O(1)$). The dynamics of the filament moving with the
binormal flow involves a velocity $O(\Gamma \log (\nu t))$ so that the
displacement of the filament is $O(\Gamma t\log (\nu t))$. Since $\Gamma
=O(t)$ we find that the filament displaces a distance $O(1)$ during a time $%
O(\nu ^{-1})$ and in that period of time the filament becomes $O(1)$ thick.
\end{remark}

\begin{acknowledgement}
We thank Santiago Montaner and Christian Zillinger for very valuable
discussions at the early stages of the present work. MAF was supported by
projects TED2021-131530B-I00 and PID2020-113596GB-I0. LV is funded by MICINN
(Spain) projects Severo Ochoa CEX2021-001142, and PID2021-126813NB-I00 (ERDF
A way of making Europe), and by Eusko Jaurlaritza project IT1615-22 and BERC
program.
\end{acknowledgement}

\section{Appendix: $O(1)$ corrections to LIA\protect\bigskip}

Our aim in this section is to estimate the Biot-Savart integral outside an
interval $s\in \left[ -\varepsilon ,\varepsilon \right] $ around a point
along the filament (that we consider, without loss of generality,
corresponding to $s=0$ and $\mathbf{r}(0,t)=\mathbf{0}$). We will consider
for the sake of simplicity an infinite filament, but the result extends
trivially to a closed filament . Let's denote%
\begin{equation*}
\mathbf{I}=\int_{%
\mathbb{R}
\diagdown \left[ -\varepsilon ,\varepsilon \right] }\frac{\mathbf{r}%
_{s}(s,t)\times \mathbf{r}(s,t)}{\left\vert \mathbf{r}(s,t)\right\vert ^{3}}%
ds.
\end{equation*}%
By introducing%
\begin{equation*}
\frac{1}{\left\vert \mathbf{r}(s,t)\right\vert ^{3}}=\sigma _{ss}(s,t),
\end{equation*}%
we have%
\begin{eqnarray*}
\sigma _{s}(s,t) &=&-\int_{s}^{\infty }\frac{1}{\left\vert \mathbf{r}%
(s^{\prime },t)\right\vert ^{3}}ds^{\prime },\ \ s>0, \\
\sigma _{s}(s,t) &=&\int_{-\infty }^{s}\frac{1}{\left\vert \mathbf{r}%
(s^{\prime },t)\right\vert ^{3}}ds^{\prime },\ \ s<0, \\
\sigma _{s}(\varepsilon ,t) &\sim &-\frac{1}{2\varepsilon ^{2}},\ \sigma
_{s}(-\varepsilon ,t)\sim \frac{1}{2\varepsilon ^{2}}\ \ \text{as }%
\varepsilon \rightarrow 0^{+}.
\end{eqnarray*}

We can compute now
\begin{eqnarray*}
I_{1} &=&\int_{\varepsilon }^{\infty }\frac{\mathbf{r}_{s}(s,t)\times
\mathbf{r}(s,t)}{\left\vert \mathbf{r}(s,t)\right\vert ^{3}}ds=\frac{1}{2}%
\frac{\mathbf{r}_{s}(\varepsilon ,t)\times \mathbf{r}(\varepsilon ,t)}{%
\varepsilon ^{2}}-\int_{\varepsilon }^{\infty }\sigma _{s}(s)\left( \mathbf{r%
}_{ss}(s,t)\times \mathbf{r}(s,t)\right) ds \\
&=&\frac{1}{2}\frac{\mathbf{r}_{s}(\varepsilon ,t)\times \mathbf{r}%
(\varepsilon ,t)}{\varepsilon ^{2}}-\int_{\varepsilon }^{\infty }\sigma
_{s}(s)\left( \kappa (s,t)\mathbf{n}(s,t)\times \mathbf{r}(s,t)\right) ds,
\end{eqnarray*}%
where%
\begin{eqnarray*}
\int_{\varepsilon }^{\infty }(s\sigma _{s}(s))\frac{\kappa (s,t)\mathbf{n}%
(s,t)\times \mathbf{r}(s,t)}{s}ds &=&-\mu (\varepsilon )\frac{\kappa
(\varepsilon ,t)\mathbf{n}(\varepsilon ,t)\times \mathbf{r}(\varepsilon ,t)}{%
\varepsilon } \\
&&-\int_{\varepsilon }^{\infty }\mu (s)\frac{d}{ds}\frac{\kappa (s,t)\mathbf{%
n}(s,t)\times \mathbf{r}(s,t)}{s}ds,
\end{eqnarray*}%
with%
\begin{eqnarray*}
\mu _{s}(s) &=&s\sigma _{s}(s), \\
\mu (s) &\sim &-\frac{1}{2}\log s\text{, as }\left\vert s\right\vert
\rightarrow 0,
\end{eqnarray*}%
and%
\begin{eqnarray*}
I_{2} &=&\int_{-\infty }^{-\varepsilon }\frac{\mathbf{r}_{s}(s,t)\times
\mathbf{r}(s,t)}{\left\vert \mathbf{r}(s,t)\right\vert ^{3}}ds \\
&=&\frac{1}{2}\frac{\mathbf{r}_{s}(-\varepsilon ,t)\times \mathbf{r}%
(-\varepsilon ,t)}{\varepsilon ^{2}}-\int_{-\infty }^{-\varepsilon }\sigma
_{s}(s)\left( \mathbf{r}_{ss}(s,t)\times \mathbf{r}(s,t)\right) ds \\
&=&\frac{1}{2}\frac{\mathbf{r}_{s}(-\varepsilon ,t)\times \mathbf{r}%
(-\varepsilon ,t)}{\varepsilon ^{2}}-\int_{-\infty }^{-\varepsilon }\sigma
_{s}(s)\left( \kappa (s,t)\mathbf{n}(s,t)\times \mathbf{r}(s,t)\right) ds,
\end{eqnarray*}%
where%
\begin{eqnarray*}
\int_{-\infty }^{-\varepsilon }(s\sigma _{s}(s))\frac{\kappa (s,t)\mathbf{n}%
(s,t)\times \mathbf{r}(s,t)}{s}ds &=&-\mu (-\varepsilon )\frac{\kappa
(-\varepsilon ,t)\mathbf{n}(-\varepsilon ,t)\times \mathbf{r}(-\varepsilon
,t)}{\varepsilon } \\
&&-\int_{-\infty }^{-\varepsilon }\mu (s)\frac{d}{ds}\frac{\kappa (s,t)%
\mathbf{n}(s,t)\times \mathbf{r}(s,t)}{s}ds,
\end{eqnarray*}%
with%
\begin{eqnarray*}
\mu _{s}(s) &=&s\sigma _{s}(s), \\
\mu (s) &\sim &\frac{1}{2}\log (-s)\text{, as }\left\vert s\right\vert
\rightarrow 0.
\end{eqnarray*}%
Since $\mathbf{r}(-\varepsilon ,t)\simeq -\mathbf{r}(\varepsilon ,t)$, $%
\mathbf{t}(-\varepsilon ,t)\simeq \mathbf{t}(\varepsilon ,t)$ and\textbf{\ }$%
\mathbf{n}(-\varepsilon ,t)\simeq \mathbf{n}(\varepsilon ,t)$ for $%
\varepsilon $ sufficiently small and $\mathbf{r}(\varepsilon ,t)=\mathbf{t}%
(0,t)\varepsilon +O(\varepsilon ^{2})$ we get%
\begin{eqnarray*}
I_{1}+I_{2} &=&-\left( \kappa (0,t)\mathbf{n}(0,t)\times \mathbf{t}%
(0,t)\right) \log \varepsilon +\int_{%
\mathbb{R}
\diagdown \left[ -\varepsilon ,\varepsilon \right] }\mu (s)\frac{d}{ds}\frac{%
\kappa (s,t)\mathbf{n}(s,t)\times \mathbf{r}(s,t)}{s}ds \\
&&+O(\varepsilon ),
\end{eqnarray*}%
where we can see clearly separated the $O(\log \varepsilon )$, $O(1)$ and $%
O(\varepsilon )$ contributions. Note that%
\begin{equation*}
\frac{d}{ds}\frac{\kappa (s,t)\mathbf{n}(s,t)\times \mathbf{r}(s,t)}{s}%
=\left( -\kappa _{s}(0,t)\mathbf{b}(0,t)-\kappa (0,t)\tau (0,t)\mathbf{n}%
(0,t)\right) +O(s),
\end{equation*}%
so that the integrand is $O(\log \left\vert s\right\vert )$ and the integral
coincides with the integral over $%
\mathbb{R}
$ plus an $O(\varepsilon \log \varepsilon )$ correction.

We conclude then that%
\begin{equation*}
\lim_{\varepsilon \rightarrow 0}(\mathbf{I}-\kappa (0,t)\mathbf{b}(0,t)\log
\varepsilon )=\lim_{\varepsilon \rightarrow 0}\left( \int_{%
\mathbb{R}
\diagdown \left[ -\varepsilon ,\varepsilon \right] }\frac{\mathbf{r}%
_{s}(s,t)\times \mathbf{r}(s,t)}{\left\vert \mathbf{r}(s,t)\right\vert ^{3}}%
ds-\kappa (0,t)\mathbf{b}(0,t)\log \varepsilon \right) ,
\end{equation*}%
is well defined and bounded. Analogously, for any point $\mathbf{r}(s,t)$
along the filament we have%
\begin{equation*}
\mathbf{A}(\mathbf{r}(s,t))\equiv \lim_{\varepsilon \rightarrow 0}\left(
\int_{%
\mathbb{R}
\diagdown \left[ -\varepsilon ,\varepsilon \right] }\frac{\mathbf{r}%
_{s}(s^{\prime },t)\times (\mathbf{r}(s,t)-\mathbf{r}(s^{\prime },t))}{%
\left\vert \mathbf{r}(s^{\prime },t)-\mathbf{r}(s^{\prime },t)\right\vert
^{3}}ds^{\prime }+\kappa (0,t)\mathbf{b}(s,t)\log \varepsilon \right) ,
\end{equation*}%
is well defined and bounded. It is simple to show that $\mathbf{A}(\mathbf{r}%
(s,t))$, as a function of $s$, is a smooth vector field provided the curve $%
\mathbf{r}(s,t)$ is sufficiently smooth.

\end{document}